\documentclass[11pt,draft,namelimits,sumlimits]{amsart}
\usepackage{citesort}
\usepackage{amssymb,amsmath}
\usepackage[mathscr]{eucal}

\newcounter{mtheorem}
\newtheorem{mtheorem}[mtheorem]{Theorem}

\newtheorem{theorem}{Theorem}[section]
\newtheorem{lemma}[theorem]{Lemma}
\newtheorem{prop}[theorem]{Proposition}
\newtheorem{corollary}[theorem]{Corollary}

\newtheorem{definition}[theorem]{Definition}

\theoremstyle{remark}
\newtheorem{remark}[theorem]{Remark}

\newtheorem{convention}[theorem]{Convention}

\numberwithin{equation}{section}

\newcommand{\pp}{\mathsf{p}}
\newcommand{\tp}{\tau}
\newcommand{\tb}{{\underline{\tau}}}
\newcommand{\mb}{{\underline{m}}}
\newcommand{\pt}{\pp_\tau}

\newcommand{\ptb}{\pp_{\tb}}
\newcommand{\pthat}{\widehat{\pp}_\tau}
\newcommand{\ptphat}{\widehat{\pp}_{\tp}}
\newcommand{\ptbhat}{\widehat{\pp}_{\tb}}

\newcommand{\RRR}{\mathsf{R}}
\newcommand{\SSS}{\mathsf{S}}
\newcommand{\TTT}{\mathsf{T}}
\newcommand{\QQQ}{\mathsf{Q}}
\newcommand{\PPP}{\mathsf{P}}

\newcommand{\sbar}{\underline{\SSS}}
\newcommand{\tbar}{\underline{\TTT}}
\newcommand{\sbartilde}{\widetilde{\underline{\SSS}}}
\newcommand{\tbartilde}{\widetilde{\underline{\TTT}}}

\newcommand{\qtilde}{\widetilde{\QQQ}}
\newcommand{\rtilde}{\widetilde{\RRR}}
\newcommand{\stilde}{\widetilde{\SSS}}
\newcommand{\ttilde}{\widetilde{\TTT}}
\newcommand{\ptilde}{\widetilde{\PPP}}

\newcommand{\group}{\mathscr{G}}  
\newcommand{\skernel}{\mathscr{K}}

\newcommand{\vspan}{\operatorname{span}}
\newcommand{\join}{\operatorname{Join}}
\newcommand{\beq}{\begin{equation}}
\newcommand{\eeq}{\end{equation}}
\newcommand{\bea}{\begin{eqnarray}}
\newcommand{\eea}{\end{eqnarray}}
\newcommand{\cy}{{C}alabi-{Y}au\ }
\newcommand{\ka}{{K}\"ahler\ }

\newcommand{\slg}{special {L}agrangian\ }
\newcommand{\Slg}{{S}pecial {L}agrangian\ }

\newcommand{\C}{\mathbb{C}}
\newcommand{\R}{\mathbb{R}}
\newcommand{\Q}{\mathbb{Q}}
\newcommand{\Z}{\mathbb{Z}}
\newcommand{\N}{\mathbb{N}}

\newcommand{\Sph}{\mathbb{S}}
\newcommand{\ra}{\rightarrow}

\newcommand{\sn}{\operatorname{sn}}

\newcommand{\cnd}{\operatorname{cn}}
\newcommand{\dn}{\operatorname{dn}}
\newcommand{\diag}{\operatorname{diag}}
\newcommand{\vol}{\operatorname{Vol}}
\newcommand{\graph}{\operatorname{graph}}

\newcommand{\Real}{\operatorname{Re}}
\newcommand{\Imag}{\operatorname{Im}}

\newcommand{\sech}{\operatorname{sech}}
\newcommand{\Xta}{X_{\tau,\alpha}}
\newcommand{\Xunderta}{\underline{X}_{\tau,\alpha}}
\newcommand{\Xunder}{\underline{X}}
\newcommand{\funder}{\underline{f}}
\newcommand{\fundertilde}{\widetilde{{f}}}
\newcommand{\cunder}{\underline{c}\,}
\newcommand{\rhounder}{\underline{\rho}}
\newcommand{\Ytilde}{\widetilde{Y}}
\newcommand{\Yhat}{\widehat{Y}}
\newcommand{\Stilde}{\widetilde{S}}
\newcommand{\Shat}{\widehat{S}}
\newcommand{\Yzeta}{Y_{\boldsymbol{\zeta}}}

\newcommand{\phiunder}{{\underline{\phi}}}

\newcommand{\xunder}{{\underline{x}}}
\newcommand{\xover}{{\overline{x}}}
\newcommand{\xboth}{{\overline{\underline{x}}}}
\newcommand{\thetadisloc}{{\theta_{dislocation}}}
\newcommand{\thetatwist}{{\theta_{twisting}}}
\newcommand{\thetaglue}{{\theta_{gluing}}}
\newcommand{\Lcal}{{\mathcal{L}}}
\newcommand{\Rcal}{{\mathcal{R}}}
\newcommand{\Jcal}{{\mathcal{J}}}
\newcommand{\Ecal}{{\mathcal{E}}}
\newcommand{\Wcal}{{\mathcal{W}}}
\newcommand{\Lchi}{{\mathcal{L}_\chi}}

\newcommand{\ftbold}{f_{\mathbf{t}}}
\newcommand{\fqbold}{f_{\mathbf{q}}}

\newcommand{\zetabold}{{\boldsymbol{\zeta}}}
\newcommand{\zerobold}{{\boldsymbol{0}}}
\newcommand{\xibold}{{\boldsymbol{\xi}}}
\newcommand{\Ubold}{{\boldsymbol{U}}}
\newcommand{\Uboldunder}{{\underline{\boldsymbol{U}}}}
\newcommand{\wt}{\widetilde{w}}
\newcommand{\psihat}{\widehat{\psi}}
\newcommand{\psicheck}{\check{\psi}}

\newcommand{\sym}{{sym}}
\newcommand{\cyl}{{\Sph^1\times\R}}

\begin{document}

\title[Higher Genus Special Lagrangian Cones]{Special Lagrangian cones with higher genus links}

\author[M.~Haskins]{Mark~Haskins}
\address{Department of Mathematics, South Kensington Campus, Imperial
College London}
\email{m.haskins@imperial.ac.uk}

\author[N.~Kapouleas]{Nikolaos~Kapouleas}
\address{Department of Mathematics, Brown University, Providence,
RI 02912} \email{nicos@math.brown.edu}





\begin{abstract}
For every odd natural number $g=2d+1$ we prove the existence of a countably infinite family of special Lagrangian cones in $\C^3$
over a closed Riemann surface of genus $g$, using a geometric PDE gluing method.
\end{abstract}

\maketitle
\section{Introduction}
\label{S:intro}
\nopagebreak
Let $M$ be a \cy manifold of complex dimension $n$ with \ka form $\omega$ and non-zero parallel
holomorphic $n$-form $\Omega$ satisfying a normalization condition. Then $\Real{\Omega}$ is a calibrated
form whose calibrated submanifolds are called special Lagrangian (SLG) submanifolds \cite{harvey:lawson}.
Moduli spaces of SLG submanifolds have appeared recently in string theory 
\cite{becker,syz,hitchin:moduli,hitchin:slg:lectures,joyce:syz}. 
On physical grounds, Strominger, Yau and Zaslow argued that a \cy manifold $M$
with its mirror partner $\widehat{M}$ admits a (singular) fibration by SLG tori, 
and that $\widehat{M}$ should be obtained by compactifying the dual fibration \cite{syz}.
To make their ideas rigorous one needs control over the singularities and compactness 
properties of families of SLG submanifolds. 
In dimensions three and higher these properties
are not well understood. As a result there has been considerable recent interest
in singular SLG subvarieties 
\cite{butscher,gross:slg:examples,gross:slgfib1,gross:slgfib2,haskins:thesis,haskins:slgcones,haskins:complexity,joyce:sing:survey,y:lee,schoen:wolfson,schoen:wolfson:1996}.
In particular some gluing constructions \cite{butscher:thesis,butscher,y:lee,joyce:sing:survey}
using Lawlor necks \cite{lawlor}
have been successfully carried out.

One natural class of singular SLG varieties is the class of SLG varieties with isolated 
conical singularities \cite{joyce:sing:survey}. 
Loosely speaking, these are compact SLG varieties 
of \cy manifolds which are singular at a finite number of points, near each of which
they resemble asymptotically some SLG cone $C$ in $\C^n$ with the origin as the only singular
point of $C$. 
This motivates the recent interest in constructing SLG cones in $\C^n$ with an
isolated singularity at the origin 
\cite{haskins:thesis,haskins:slgcones,mcintosh:slg,carberry:mcintosh,joyce:symmetries}.
Until recently few examples of such SLG cones in $\C^n$ were known. 
Recently, using techniques from equivariant differential geometry, symplectic geometry and integrable systems
many new families of examples have been constructed 
\cite{haskins:thesis,haskins:slgcones,joyce:symmetries,mcintosh:slg,carberry:mcintosh}.

Very interesting results have been obtained for SLG cones in $\C^3$.
If $C$ is a SLG cone in $\C^3$ whose intersection with $\Sph^5$ is topologically a $2$-sphere
then $C$ must be a plane \cite{yau} \cite[Thm. 2.7]{haskins:slgcones}. 
The oldest (circa 1974) and simplest singular SLG cone in $\C^3$ is the cone over
a flat $2$-torus in $\Sph^5$ (the so-called Clifford torus) \cite[Ex. 3.18]{haskins:complexity}.
It is invariant under the diagonal subgroup $T^2 \subset \text{SU(3)}$. 
Recently, countably infinite families of SLG cones over $2$-tori invariant under
$\text{U(1)}$-subgroups of $\text{SU(3)}$ were constructed in \cite{haskins:slgcones,haskins:thesis}
using integrable systems techniques. These $\text{U(1)}$-invariant examples can be described explicitly
in terms of elliptic functions and integrals.
A special case of this construction is crucial to this paper.
Very recently, the existence of many more SLG cones over $2$-tori in $\Sph^5$ 
was proved using more sophisticated
algebro-geometric tools from integrable systems \cite{carberry:mcintosh,mcintosh:slg}.
Recent results \cite[Thms. C \& D]{haskins:complexity} prove that these
examples rapidly become geometrically very complex.

In stark contrast to the case of SLG cones over $2$-spheres or $2$-tori 
practically nothing is known about SLG cones over surfaces of higher genus.
The integrable systems methods effective in the genus one case have so far yielded no insight into the
higher genus case.
The intersection of a SLG cone in $\C^3$ with $\Sph^5$ is called a special Legendrian surface in $\Sph^5$.
In this paper we use a geometric PDE gluing method
to construct special Legendrian closed surfaces in $\Sph^5$ of any odd genus $g>1$---their
associated cones are SLG.
The particular kind of methods we use relates more
closely to the methods developed in \cite{schoen,kapouleas:annals},
especially as they evolved and were systematized in \cite{kapouleas:wente}.
The basic building blocks of our construction are provided by
the $\text{U(1)}$-invariant SLG $T^2$-cones constructed in \cite{haskins:slgcones,haskins:thesis}
and described in Section \ref{S:u1:tori}.
Our main result is the following existence result.

\begin{mtheorem}
For every positive odd integer $g=2d+1$ there exist infinitely many different
special Lagrangian cones in $\C^3$ whose intersection with $\Sph^5$ is
a closed oriented surface of genus $g$.
\end{mtheorem}

We now give a more detailed but nontechnical description of our construction.
We start by describing the basic building blocks we use.
In \cite{haskins:slgcones,haskins:thesis}
a one-parameter family of 
$\text{U(1)}$-invariant
special Legendrian cylinders is constructed.
We call the parameter of the family $\tau$ and its range is $[0,1/3\sqrt{3})$ (see \ref{P:conformal}).
The corresponding immersions
$X_\tau: \Sph^1 \times \R \ra \Sph^5\subset \C^3$---see \ref{D:X:tau}, \ref{E:X:tau}, and \ref{L:Per}---are
periodic for $\tau\ne0$ with the period controlled by a number $\pthat$.
The fundamental domain approaches an equatorial $\Sph^2$ as $\tau\to0$,
where at two antipodal points the surface diverges from the equatorial $\Sph^2$,
and approaches a Lagrangian catenoid of size of order $\sqrt{\tau\,}$
at each of the two antipodal points---see \ref{R:catenoid}.
Through these Lagrangian catenoids the spherical region in consideration
connects to the spherical regions of the adjacent fundamental domains.
This behavior is reminiscent of the Delaunay cylinders---see for example \cite{kapouleas:annals},
although the positivity of the parameter and the rotational character of the period
are more in analogy with the Wente cylinders---see for example \cite{kapouleas:wente}.
Moreover, when $\pthat$, which tends to $\pi/2$ as $\tau\to0$,
is a rational multiple of $\pi$,
$X_\tau$ factors through a special Legendrian embedding of the torus.
For simplicity we concentrate on the case where the torus contains $4\mb-1$ fundamental
regions, where $\mb$ is large and (equivalently) the corresponding $\tau$
is small---see \ref{E:p:hat} and \ref{L:Per}.

We want to combine $g$ copies of the tori we just described to obtain a closed
special Legendrian surface of genus $g$.
The technical aspects of constructions of this type are significantly simplified
when the maximal possible symmetry is imposed.
To achieve such symmetry we use an $\text{SU(3)}$-rotation of order $g$---see \ref{E:rotn}---to
obtain $g$ copies of the torus in consideration.
There is an equatorial $\Sph^2$ invariant under the rotation,
which is close to one spherical region in each of the copies.
By appropriately fusing these spherical regions with the equatorial $\Sph^2$
we obtain a closed surface of genus $g$.
More precisely we consider an equatorial regular $2g$-gon on the equatorial $\Sph^2$
at the vertices of which we attach the $2g$ Lagrangian catenoidal necks through which 
our $g$ tori connect to the common $\Sph^2$.
By carrying out the fusion carefully, following the idea in
\cite{butscher:thesis,butscher,y:lee},
we ensure that the new surface,
which we call the ``initial surface'',
is Legendrian.
Moreover its Lagrangian angle,
which is supported on the $2g$ annuli where the surface transits from the equatorial $\Sph^2$ to the tori,
is of order $\tau$.
The idea is to correct the constructed surface to being special Legendrian when $\tau$ is small enough,
that is the number of fundamental domains in each torus is large enough.
This set-up with the maximal symmetry turns out to work only when $g$ is odd,
otherwise we end up trying to glue two different tori to $\Sph^2$ at the same antipodal points.
This is the reason that our theorem assumes $g$ odd.

There is a simple construction which allows us to obtain Legendrian perturbations
of a given Legendrian surface by using a given real function $f$ on the surface (see Appendix
\ref{A:perturbation}).
From the viewpoint of the Lagrangian cone on the surface this is the same
as the construction used in 
\cite{butscher:thesis,butscher,y:lee}.
The Lagrangian angle (see \ref{S:slg} for the definition) $\theta_f$
of the perturbed surface satisfies
$$
\theta_f=\theta+(\Delta+6)f+Q_f,
$$
where $\theta$ is the Lagrangian angle of the original surface,
$\Delta$ the induced Laplacian on the original surface,
and $Q_f$ is quadratic and higher in $f$ and its derivatives.
To find a special Legendrian perturbation,
that is one with $\theta_f\equiv0$,
we have to find $f$ such that 
\addtocounter{theorem}{1}
\begin{equation}
\label{E:PDE}
(\Delta+6)f=-\theta-Q_f.
\end{equation}
This is similar to the situation in 
\cite{kapouleas:annals,kapouleas:finite,kapouleas:minimal,kapouleas:imc,kapouleas:wente,kapouleas:wente:announce,kapouleas:jdg,kapouleas:cmp,kapouleas:annals,kapouleas:bulletin}:
The Lagrangian angle here plays the role of the mean curvature there,
and the perturbation is constructed by using the gradient of $f$ instead of $f$ times the unit normal.
In some sense---see \ref{mc:lagn:angle:grad}---the gradient of the quantities $\theta$ and $f$ used here,
corresponds to the quantities used in the other constructions.
Anyway, the problem is reduced to finding a solution $f$ of \ref{E:PDE},
where because we have prepared a small $\theta$,
we are expecting a small $f$ as well.
We expect then that $Q_f$ is not a dominant term.
The main difficulty is to understand for given $E$
the linearized equation
\addtocounter{theorem}{1}
\begin{equation}
\label{E:LPDE}
(\Delta+6)f=E.
\end{equation}

This is somewhat simpler than in the cases in
\cite{kapouleas:annals,kapouleas:finite,kapouleas:minimal,kapouleas:imc,kapouleas:wente,kapouleas:wente:announce,kapouleas:jdg,kapouleas:cmp,kapouleas:annals,kapouleas:bulletin}
where the corresponding equation is 
$(\Delta+|A|^2)f=E$,
where $|A|^2$ is the square of the length of the second fundamental form.
In those cases the term $|A|^2$ makes regions of high curvature important spectrally,
something which does not happen in the current case---see \ref{R:catenoid}.
Nevertheless the usual difficulties persist:
The regions on the initial surface which are close to equatorial $\Sph^2$'s,
induce the linearized operator $\Delta+6$ to inherit small eigenvalues
from the five-dimensional kernel it has on an equatorial $\Sph^2$.
Actually, the symmetry we have imposed on the construction (see 
\ref{E:symmetries}),
induces symmetries on the functions $f$ (see 
\ref{D:f:symmetries}),
which simplify the approximate kernel on each spherical region to a two-dimensional one,
except for the spherical region where the fusion occurs which has trivial approximate kernel
(see \ref{P:app:ker}).

Dealing with this high-dimensional approximate kernel---in our case the dimension is $4\mb-2$,
is a major aspect of these constructions
\cite{schoen,kapouleas:annals,kapouleas:finite,kapouleas:minimal,kapouleas:imc,kapouleas:wente,kapouleas:wente:announce,kapouleas:jdg,kapouleas:cmp,kapouleas:annals,kapouleas:bulletin,kapouleas:survey}.
We refer the reader to \cite{kapouleas:survey} for a discussion of the general systematic approach.
Another aspect is that we have to understand carefully the interaction
between the spherical regions.
These connect through regions which are conformal to long cylinders.
As in \cite{kapouleas:wente} we use 
a Fourier decomposition along the meridians to understand the interaction.
In our case the five lowest harmonics on the meridians are important in ensuring the
decay we need, as opposed to the lowest three in \cite{kapouleas:wente}.
Actually because of the symmetries $\sin s$, $\cos s$, and $\sin (2s)$ are not allowed,
and we have only two left,
namely the constants and $\cos (2s)$.
All this leads to Proposition
\ref{P:solution},
where the linear equation \ref{E:LPDE} is solved with appropriate decay estimates,
modulo a $4\mb$-dimensional ``extended substitute kernel'' $\skernel$ (see \ref{D:win}).

Since we only solve \ref{E:LPDE} modulo elements of $\skernel$,
we need to find a way of prescribing small elements of $\skernel$ as part of the initial
Lagrangian angle $\theta$.
This is achieved by applying the so-called geometric principle
\cite{kapouleas:finite,kapouleas:minimal,kapouleas:imc,kapouleas:wente,kapouleas:wente:announce,kapouleas:survey}:
We dislocate the initial surface by moving the spherical regions relative to each other
by using the isometries of the ambient $\Sph^5$ responsible for the existence
of the approximate kernel in the first place.
In our case this means using
$\ttilde_x$ and $\qtilde_x$ which define one-parameter subgroups of $\text{SU(3)}$
(see \ref{D:sym:tilde}).
The generating Killing fields of these one-parameter subgroups are the ones which induce
the approximate kernel on the spherical regions (see \ref{E:f:tq}).
$\ttilde_x$ acts by introducing some ``sliding'' along the ``axis'' of the torus---the points
fixed by the $\text{U(1)}$ action leaving the torus invariant.
It resembles translations along the axis of a Delaunay surface.
$\qtilde_x$ acts by introducing some ``twisting'' along the ``axis''.
It exists because the codimension is higher than in the Delaunay surface case
where no analogous Euclidean motions in $\mathbb{E}^3$ exist.

Technically, we work as follows to apply the geometric principle:
To create extended substitute kernel corresponding to the central spherical region (see \ref{P:phiunder0}),
where the fusion of the tori occurs, we introduce two parameters $\zeta_1$ and $\zeta_2$
in the construction of the initial surface.
The original construction corresponds to $\zeta_1=\zeta_2=0$.
$\zeta_1$ controls a change of the parameter $\tau$ of the original torus,
which induces a change in the period $\pthat$,
which induces a ``mismatch'' between the central spherical region and the torus
by ``sliding'' one relative to the other.
To introduce now a ``mismatch'' which corresponds to a relative ``twisting'' controlled by $\zeta_2$,
we had to perturb the original special Legendrian cylinders of Haskins
to special Legendrian cylinders whose period besides ``sliding'' had some ``twisting'' as well
(see  \ref{L:a-sym}, \ref{L:sliding:twisting:X}, and \ref{T:Xta}).

To create substitute kernel corresponding to the other spherical regions (see \ref{L:phi:inp}),
we can fortunately work at the linear level.
This means that instead of introducing more parameters controlling perturbations
of the initial surface, we construct instead functions on the surface which can be included in the solution $f$
in order to induce the desired effect.
Checking this, that is estimating the amount of substitute kernel created,
is facilitated as usual by using a balancing type formula
and calculating the forces involved
\cite{schoen,kapouleas:annals,kapouleas:finite,kapouleas:minimal,kapouleas:imc,kapouleas:wente,kapouleas:wente:announce,kapouleas:jdg,kapouleas:cmp,kapouleas:annals,kapouleas:bulletin,kapouleas:survey}.
A small technical innovation is that this is done here entirely at the linear level.
This way the balancing formula reduces to Green's second identity \cite{gilbarg}
and the calculation is considerably simplified.
All this leads to Proposition
\ref{P:Phixi},
where the (approximate) prescription of small elements of $\skernel$ is clearly given
in a way easy to use later.

Another small technical innovation introduced in this paper is using a scaling
argument to simplify the estimation of the nonlinear terms (see \ref{L:quadratic}).
Finally we use Schauder's fixed point theorem
\cite[Theorem 11.1]{gilbarg} 
to complete the proof of the main result of the paper (see \ref{T:main}).

The paper is organized in seven sections and three appendices.
Section 1 consists of the introduction and some remarks on the notation.
In Section \ref{S:slg} we recall basic facts and definitions from special Lagrangian geometry.
In Section \ref{S:u1:tori} we describe
in detail the $\text{U}(1)$-invariant SLG $T^2$-cones which form basic building blocks
of our construction.
Section \ref{S:init:surf} describes the two-parameter family we discussed above
of Legendrian immersions $Y_\zetabold: M \ra \Sph^5$, 
where $M$ is a compact oriented surface of genus $g$.
Section \ref{S:linear} discusses the linear theory culminating in stating and proving
\ref{P:solution}.
In Section 6 the geometric principle is applied to carry out the creation
of extended substitute kernel as formalized in \ref{P:Phixi}.
In Section 7 the nonlinear terms are discussed and the main results established.

Appendix \ref{A:elliptic} contains the 
basic definitions and facts concerning Jacobi elliptic functions and elliptic integrals
that are needed in Section \ref{S:u1:tori}.
Appendix \ref{A:interpolation} 
describes how to transit from one Legendrian immersion in $\Sph^5$ to another.
This material is used in the construction of the initial surfaces of Section \ref{S:init:surf}.
Finally, appendix C describes how to use a small function $f$ to perturb a Legendrian immersion $X$
into $\Sph^5$ into a nearby Legendrian immersion $X_f$.

\subsection*{Notation and conventions}
$\phantom{ab}$
\nopagebreak

In this paper we use weighted H\"{o}lder norms.
A definition which works for our purposes in this paper is the following: 
\addtocounter{theorem}{1}
\begin{equation}
\label{E:weightedHolder}
\|\phi: C^{k,\beta}(\Omega,g,f)\|:=
\sup_{x\in\Omega}\frac{\,\|\phi:C^{k,\beta}(\Omega\cap B_x,g)\|\,}{f(x)},
\end{equation}
where $\Omega$ is a domain inside a Riemannian manifold $(M,g)$,
$f$ is a weight function on $\Omega$,
$B_x$ is a geodesic ball centered at $x$ and of radius the minimum of
$1$ and half the injectivity radius at $x$.

We will be using extensively cut-off functions and for this reason we adopt the
following notation:
We fix a smooth function $\Psi:\R\to[0,1]$ with the following properties:
\newline
(i).
$\Psi$ is nondecreasing.
\newline
(ii).
$\Psi\equiv1$ on $[1,\infty]$ and $\Psi\equiv0$ on $(-\infty,-1]$.
\newline
(iii).
$\Psi-\frac12$ is an odd function.
\newline
Given then $a,b\in \R$ with $a\ne b$ we define a smooth function $\psi[a,b]:\R\to[0,1]$
by
\addtocounter{theorem}{1}
\begin{equation}
\label{E:psiab}
\psi[a,b]=\Psi\circ L_{a,b},
\end{equation}
where $L_{a,b}:\R\to\R$ is the linear function defined by the requirements $L(a)=-3$ and $L(b)=3$.

Clearly then $\psi[a,b]$ has the following properties:
\newline
(i).
$\psi[a,b]$ is weakly monotone.
\newline
(ii).
$\psi[a,b]=1$ on a neighborhood of $b$ and 
$\psi[a,b]=0$ on a neighborhood of $a$.
\newline
(iii).
$\psi[a,b]+\psi[b,a]=1$ on $\R$.

We will denote the span of vectors $e_1,\dots,e_k$ with coefficients in a field $F$
by $\left< e_1,\dots,e_k\right>_F$.

\subsection*{Acknowledgments}
The authors would like to thank Rick Schoen and Karen Uhlenbeck
for their interest and support.
We also thank ETH for its hospitality and support.
M.H. would like to thank IHES for its hospitality and EPSRC
for funding his research.

\section{Special Lagrangian cones and special Legendrian submanifolds of $\Sph^{2n-1}$}
\label{S:slg}
\nopagebreak

In this section we recall basic facts about \slg geometry in $\C^n$, \slg cones in $\C^n$
and their connection to minimal Legendrian submanifolds of $\Sph^{2n-1}$.
\Slg geometry is an example of a \textit{calibrated geometry} \cite{harvey:lawson}. We begin by defining
\textit{calibrations} and \textit{calibrated submanifolds}.

\subsection*{Calibrations and Special Lagrangian geometry in $\C^n$}
\label{SS:slg}
$\phantom{ab}$
\nopagebreak

Let $(M,g)$ be a Riemannian manifold. Let $V$ be an oriented tangent $p$-plane on $M$, \textit{i.e.}
a $p$-dimensional oriented vector subspace of some tangent plane $T_xM$ to $M$. The restriction
of the Riemannian metric to $V$, $g|_V$, is a Euclidean metric on $V$ which together with the
orientation on $V$ determines a natural $p$-form on $V$, the volume form $\vol_{V}$.
A closed $p$-form $\phi$ on $M$ is a \textit{calibration} on $M$ if for every oriented tangent $p$-plane $V$ on $M$
we have $\phi|_V \le \vol_V$. Let $L$ be an oriented submanifold of $M$ with dimension $p$. $L$ is a
\textit{$\phi$-calibrated submanifold} if $\phi|_{T_xL} = \vol_{T_xL}$ for all $x\in L$.
There is a natural extension of this definition to singular calibrated submanifolds using the
language of Geometric Measure Theory and rectifiable currents \cite[\S II.1]{harvey:lawson}.
The key property of calibrated submanifolds (even singular ones) is that they are \textit{homologically volume minimizing} \cite[Thm. II.4.2]{harvey:lawson}.
In particular, any calibrated submanifold is automatically \textit{minimal}, \textit{i.e.} has vanishing mean curvature.

Let $z_1 = x_1 + i y_1, \ldots  ,z_n = x_n + i y_n$ be standard complex coordinates on $\C^n$ equipped with the Euclidean metric.
Let 
$$\omega = \frac{i}{2} \sum_{j=1}^n{dz_j \wedge d\bar{z}_j} = \sum_{j=1}^n{dx_j \wedge dy_j},$$ 
be the standard symplectic $2$-form on $\C^n$.
Define a complex $n$-form $\Omega$ on $\C^n$ by
\begin{equation}
\addtocounter{theorem}{1}
\label{E:slg:form}
\Omega = dz_1 \wedge \ldots \wedge dz_n.
\end{equation}
Then $\Real{\Omega}$ and $\Imag{\Omega}$ are real $n$-forms on $\C^n$.
$\Real{\Omega}$ is a calibration on $\C^n$ whose calibrated submanifolds we call
\textit{\slg submanifolds} of $\C^n$, or SLG $n$-folds for short.
There is a natural extension of \slg geometry to any Calabi-Yau manifold $M$ by replacing
$\Omega$ with the natural parallel holomorphic $(n,0)$-form on $M$. \Slg submanifolds
appear to play an important role in a number of interesting geometric properties
of \cy manifolds, \textit{e.g.} Mirror Symmetry \cite{syz,thomas:yau}.


Let $f:L \ra \C^n$ be a Lagrangian immersion of the oriented $n$-manifold $L$, and $\Omega$ be
the standard holomorphic $(n,0)$-form defined in \ref{E:slg:form}.
Then $f^*\Omega$ is a complex $n$-form on $L$ satisfying $|f^*\Omega| = 1$ \cite[p. 89]{harvey:lawson}. Hence we may write
\addtocounter{theorem}{1}
\begin{equation}
\label{E:lagn:phase}
f^*\Omega = e^{i\theta}\vol_L\quad \text{on \ } L,
\end{equation}
for some \textit{phase function} $e^{i\theta}:L \ra \Sph^1$.
We call $e^{i\theta}$ the \textit{phase of the oriented Lagrangian immersion $f$}.
$L$ is a SLG $n$-fold in $\C^n$ if and only if the phase function $e^{i\theta}\equiv 1$.
Reversing the orientation of $L$ changes the sign of the phase function $e^{i\theta}$.
The differential $d\theta$ is a well-defined closed $1$-form on $L$ which satisfies
\addtocounter{theorem}{1}
\begin{equation}
\label{mc:lagn:angle}
d\theta = \iota_H \omega,
\end{equation}
where $H$ is the mean curvature vector of $L$. In particular, \ref{mc:lagn:angle} implies
that a connected component of $L$ is minimal if and only if the phase function $e^{i\theta}$
is constant. 
\ref{mc:lagn:angle} may also be restated as
\addtocounter{theorem}{1}
\begin{equation}
\label{mc:lagn:angle:grad}
H = -J \nabla \theta,
\end{equation}
where $J$ and $\nabla$ denote the
standard complex structure and gradient on $\C^n$ respectively.
For a general Lagrangian submanifold it is not possible to find a global lift 
of the $\Sph^1$ valued phase function $e^{i\theta}$ to a real function $\theta$, although
of course such a lift always exists locally. However, we will always be interested
in Lagrangian submanifolds for which such a global lift $\theta$ does exist.
In this case we call $\theta:L \ra \R$ the \textit{Lagrangian angle} of $L$.
In particular, any Lagrangian submanifold which is sufficiently close to 
a special Lagrangian submanifold will have a globally well-defined Lagrangian angle $\theta$.

\subsection*{Contact geometry}
$\phantom{ab}$
\nopagebreak

We recall some basic definitions from contact geometry \cite{geiges,arnold,marle,mcduff:salamon}.
Let $M$ be a smooth manifold of dimension $2n+1$, and let $\xi$ be a hyperplane field on $M$.
$\xi$ is a (cooriented) \textit{contact structure} on $M$ if there exists a $1$-form $\gamma$ so that
$\ker{\gamma} = \xi$ and
\addtocounter{theorem}{1}
\begin{equation}
\label{contact:nonint}
\gamma \wedge (d\gamma)^n \neq 0.
\end{equation}
The pair $(M,\xi)$ is called a \textit{contact manifold}, and the $1$-form $\gamma$ a \textit{contact form} defining $\xi$.
Condition \ref{contact:nonint} is equivalent to the condition that $(d\gamma)^n|_{\xi} \neq 0$.
In particular, for each $p\in M$ the $2n$-dimensional subspace $\xi_p \subset T_pM$ endowed with the
$2$-form $d\gamma|{\xi_p}$ is a symplectic vector space.
Given a contact form $\gamma$ on $M$, the \textit{Reeb vector field} $R_\gamma$ is the unique vector field on $M$
satisfying
$$\iota(R_\gamma) d\gamma \equiv  0, \qquad \gamma(R_\gamma) \equiv  1.$$
Let $(M,\xi=\ker{\gamma})$ be a contact manifold.
A submanifold $L$ of $(M,\xi)$ is an \textit{integral submanifold of $\xi$} (also called an isotropic
submanifold) if $T_x L \subset \xi_{x}$ for all $x\in L$. Equivalently $L$
is an integral submanifold of $\xi$ if $\gamma|_L =0$.
If $M$ has dimension $2n+1$ then it is well-known \cite[Prop 2.17]{geiges}
that any integral submanifold $L$ of $(M,\xi)$ has dimension less than or equal to $n$. A submanifold $L$ of $(M^{2n+1},\xi)$ is \textit{Legendrian}
if it is an integral submanifold of maximal dimension, \textit{i.e.} $\dim{L}=n.$

\subsection*{Special Legendrian submanifolds and special Lagrangian cones}
$\phantom{ab}$


For any compact oriented embedded (but not necessarily connected) submanifold $\Sigma
\subset \Sph^{2n-1}(1)\subset \C^n$ define the \textit{cone on
$\Sigma$},
$$ C(\Sigma) = \{ tx: t\in \R^{\ge 0}, x \in \Sigma \}.$$
A cone $C$ in $\C^n$  (that is a subset invariant under dilations) 
is $\textit{regular}$ if there exists
$\Sigma$ as above so that $C=C(\Sigma)$, in which case we call
$\Sigma$ the \textit{link} of the cone $C$. $C'(\Sigma):=C(\Sigma) - \{0\}$ is an
embedded smooth submanifold, but $C(\Sigma)$ has an isolated
singularity at $0$ unless $\Sigma$ is a totally geodesic sphere.
Sometimes it will also be convenient to allow $\Sigma$ to be just immersed not embedded,
in which case $C'(\Sigma)$ is no longer embedded. Then we call $C(\Sigma)$ an \textit{almost regular} cone.

Let $r$ denote the radial coordinate on $\C^n$ and let
$X$ be the Liouville vector field
$$X= \frac{1}{2}r \frac{\partial}{\partial r} = \frac{1}{2}\sum_{j=1}^{n}{x_j \frac{\partial}{\partial x_j} + y_j \frac{\partial}{\partial y_j}}.$$
The unit sphere $\Sph^{2n-1}$ inherits a natural contact form
$$\gamma = \iota_X\omega | _{\Sph^{2n-1}} = {\sum_{j=1}^n{ x_j dy_j - y_j dx_j}}\Bigr|_{\Sph^{2n-1}}$$ from
its embedding in $\C^n$. 

There is a one-to-one correspondence between regular Lagrangian cones in $\C^n$ 
and Legendrian submanifolds of $\Sph^{2n-1}$. 
The Lagrangian angle or the phase of a Lagrangian cone in $\C^n$ is homogeneous of degree $0$.
We define the Lagrangian angle of a Legendrian submanifold $\Sigma$ of $\Sph^{2n-1}$ to
be the restriction to $\Sph^{2n-1}$ of the Lagrangian angle of the Lagrangian cone $C'(\Sigma)$.
We call a submanifold $\Sigma$ of $\Sph^{2n-1}$ \textit{special Legendrian} 
if the cone over $\Sigma$, $C'(\Sigma)$
is special Lagrangian in $\C^n$. In other words, $\Sigma$ is special Legendrian
if and only if its Lagrangian phase is identically $1$ or its Lagrangian angle
is identically $0$ modulo $2\pi$. 

A special Legendrian submanifold of $\Sph^{2n-1}$
is  minimal, that is it has mean curvature $H=0$. Conversely, up to rotation by a constant phase
$e^{i\theta}$ any connected minimal Legendrian submanifold of $\Sph^{2n-1}$ 
is special Legendrian. 
Using this language the goal of our paper is 
to construct special Legendrian immersions of surfaces of odd genus $g>1$
into $\Sph^5$.

\section{The $\text{U(1)}$-invariant minimal Legendrian tori}
\label{S:u1:tori}
\nopagebreak

\subsection*{Introduction}
$\phantom{ab}$
\nopagebreak

This section introduces the one-parameter family of
conformal special Legendrian immersions
$X_\tau:\Sph^1 \times \R \ra \Sph^5$ invariant under the $\text{U(1)}$ action
\addtocounter{theorem}{1}
\begin{equation}
\label{D:u1:act}
e^{is}.(w_1, w_2, w_3) =
\stilde_s(w_1,w_2,w_3) =
 (w_1, \cos{s}\ w_2 + \sin{s}\ w_3, -\sin{s}\ w_2 + \cos{s}\ w_3),
\end{equation}
where
$\tau\in [0,1/3\sqrt{3}]$
is the parameter of the family
and $\stilde_s\in\text{SU(3)}$ (see \ref{D:sym:tilde}).
More precisely
\addtocounter{theorem}{1}
\begin{equation}
\label{E:X:tau}
X_\tau(e^{is},t) = w_1(t) e_1 + w_2(t) \stilde_s(e_2) = w_1(t) e_1 + w_2(t) (\cos{s}\ e_2 - \sin{s}\ e_3),
\end{equation}
where $e_1,e_2,e_3$ is the standard basis of $\C^3$ and $w_1,w_2:\R\to\C$ are defined later in
\ref{D:X:tau}.
For a dense set of $\tau$ these cylinders factor through embedded special Legendrian $2$-tori.
When $\tau$ is sufficiently small any such embedded $2$-torus is composed of a large number of identical almost 
spherical regions connected to each of its two neighboring almost spherical regions by 
a small neck. These surfaces form the basic building blocks of our construction.

\begin{remark}
\addtocounter{equation}{1}
The $\text{U(1)}$-invariant special Legendrian immersions $X_\tau$
studied in this section are special cases of those constructed in \cite[Thm. D]{haskins:slgcones}. 
In the terminology of that paper they are the immersions $u_{0,J}$ where $J=\tau$.
This family is distinguished among all $\text{U(1)}$-invariant special Legendrian $2$-tori
because it is the only one for which the $\text{U(1)}$-action has
nontrivial fixed points. As a result it is the only family which
can limit to a two-sphere (which by \cite{yau} \cite[Thm. 2.7]{haskins:slgcones} is necessarily totally geodesic).
\end{remark}

\subsection*{The induced conformal factor $y_\tau=\rho_\tau^2$.}
$\phantom{ab}$
\nopagebreak

We begin by defining a function $y_\tau=\rho_\tau^2:\R \ra \R$ ($\rho_\tau>0$)
which will prove to be the conformal factor of the
metric induced on $\Sph^1 \times \R$ by the immersion $X_\tau$.

Given any $\tau\in [0,1/3\sqrt{3}]$ define $y_\tau(t)$ to be the unique
solution of the initial value problem
\addtocounter{theorem}{1}
\begin{equation}
\label{y:double:dot}
 \ddot{y} = -2y (3y-2); \quad  y_\tau(0) = y_{\text{max}}, \qquad \dot{y}_\tau (0)= 0,
\end{equation}
where $y_{\text{max}}>0$ denotes the largest solution of the cubic equation
\addtocounter{theorem}{1}
\begin{equation}
\label{E:cubic}
y^3 -y^2 + 4\tau^2 =0.
\end{equation}
It is straightforward to check that the equation admits the first integral
\addtocounter{theorem}{1}
\begin{equation}
\label{E:ydot}
\dot{y}^2 = -4(y^3 - y^2 + 4\tau^2)
\end{equation}
and that $y_{\text{max}}$ is the maximum value attained by the solution $y_\tau$.
For notational convenience we will usually drop the $\tau$ and refer
simply to $y$.
For each $\tau\in [0,1/3\sqrt{3}]$, the cubic
$y^3-y^2+4\tau^2$ has three real roots $y_{-} \le 0 \le  y_{\text{min}} \le y_{\text{max}}$.
All three roots are distinct except for the extreme values $\tau=0$ and
$\tau=1/3\sqrt{3}$ in which case we have $y_{\text{min}}=y_{-}=0$, $y_{\text{max}}=1$, and
$y_{\text{min}}=y_{\text{max}}=2/3$, $y_{-}=-1/3$ respectively. 

\addtocounter{equation}{1}
\begin{prop}
\label{P:conformal}
$y_\tau=\rho_\tau^2:\R \ra [y_{\text{min}},y_{\text{max}}] \subset\R$ is a smooth even function depending 
smoothly on $\tau^2$ for $\tau\in[0,1/3\sqrt{3})$.
For $\tau=0$, $y_\tau = \sech^2{t}$ and for $\tau = 1/3\sqrt{3}$, $y_\tau \equiv \frac{2}{3}$.
Moreover, for $\tau \in (0,1/3\sqrt{3})$
\begin{enumerate}
\item[(i)]  $y_{\tau}(t)=\rho_\tau^2(t)$ is given explicitly in terms of the roots $y_{-}\le 0
\le y_{\text{min}} \le y_{\text{max}}$ of $y^3-y^2+4\tau^2$  and the Jacobi elliptic function $\sn$ by
\addtocounter{theorem}{1}
\begin{equation}
\label{E:y:sn}
y_\tau(t) =\rho_\tau^2(t)= y_{\text{max}} - (y_{\text{max}}-y_{\text{min}})\sn^2{(rt,k)}
\end{equation}
where
$$
r^2 = y_{\text{max}}-y_{-},
\quad \textrm{and\ }\
k^2=\frac{y_{\text{max}}-y_{\text{min}}}{y_{\text{max}}-y_{-}}.
$$ 
\item[(ii)] $y_{\text{min}}$ and $y_{\text{max}}$ satisfy as $\tau\ra0$
\addtocounter{theorem}{1}
\begin{equation}
\label{E:ymin:tau}
y_{\text{min}} = 2\tau + O(\tau^2), \qquad y_{\text{max}} = 1 - 4\tau^2 + O(\tau^3),
\end{equation}
\item[(iii)]  $y_{\tau}$ and $\rho_\tau$ are periodic of
period $2\pt = 2K(k)/r$, where $K(k)$ is the complete elliptic
integral of the first kind and as $\tau\ra0$
\addtocounter{theorem}{1}
\begin{equation}
\label{E:pt} \pt= -\frac{1}{2}{\log{\tau}} + O(1),
\qquad
\frac{d\pt}{d\tau} = -\frac{1}{2\tau} + O(1).
\end{equation}
\item[(iv)]
$\rho_\tau$ satisfies $\rho_{\tau}(t)\ge \tfrac{1}{2}e^{-t}$ for all $\tau\in [0,\pt]$,
and for each $k\in \N$ there exists a positive constant $C(k)$ independent of $\tau$ so that for
all sufficiently small $\tau>0$ we have
\addtocounter{theorem}{1}
\begin{equation}
\label{E:exp:decay:y}
\Vert \rho_{\tau}: C^k(\Sph^1 \times [0,\pt],ds^2+dt^2,e^{-t})\Vert \le C(k).
\end{equation}
Note that $\rho_\tau$ is then controlled everywhere by its evenness and periodicity.
\end{enumerate}
\end{prop}

\begin{proof}
The formulae for $y_\tau$ in the 
limiting cases $\tau=0$ and $\tau = \frac{1}{3\sqrt{3}}$ are verified by an easy computation.
Similarly, in the case $\tau\in (0,\frac{1}{3\sqrt{3}})$ 
it is a straightforward computation using standard properties of the Jacobi elliptic 
function $\sn$ to verify that \ref{E:y:sn} satisfies equation \ref{E:ydot} with the correct
initial condition.
The periodicity and the expression for $\pt$
then follow immediately from the basic properties of $\sn$ and $K$.
The facts that $y_\tau$ is smooth and even also follow immediately from the explicit formulae
for $y_\tau$.
It remains to prove (ii),
(iii), (iv),
and the smooth dependence of $y_\tau$ 
on $\tau^2$ for $\tau\in[0,\frac{1}{3\sqrt{3}}]$.

Smooth dependence on $\tau^2$:
One can check that $y_{\text{max}}$ depends smoothly on $\tau^2 \in [0,\frac{1}{3\sqrt{3}})$.
Smooth dependence of $y_\tau$ on $\tau^2$ then follows from smooth dependence
of solutions of the initial value problem \ref{y:double:dot}.

To prove \ref{E:ymin:tau} and \ref{E:pt}
it is convenient to treat $y_{\text{min}}$ as the parameter instead of $\tau$. 
As $\tau$ increases from $0$ to $1/3\sqrt{3}$, $y_{\text{min}}$ increases
monotonically from $0$ to $\frac{2}{3}$. Specifying a value of $y_{\text{min}} \in (0,\frac{2}{3})$ uniquely determines
$\tau\in (0,1/3\sqrt{3})$ by
\addtocounter{theorem}{1}
\begin{equation}
\label{E:tau:y}
\tau = \frac{1}{2}y_{\text{min}} \sqrt{(1-y_{\text{min}})} = \frac{y_{\text{min}}}{2} - \frac{y_{\text{min}}^2}{4} + O(y_{\text{min}}^3).
\end{equation}
The first half of \ref{E:ymin:tau} follows immediately from \ref{E:tau:y}.

One advantage of using $y_{\text{min}}$ in place of $\tau$ is that since one root $y_{\text{min}}$ of the cubic $y^3-y^2+4\tau^2$
is specified the two remaining roots $y_{-}$ and $y_{\text{max}}$ are determined in terms of $y_{\text{min}}$ by solving a
quadratic equation. Solving this equation leads to
$$y_{\text{max}} = \frac{1}{2}\left( \sqrt{(1-y_{\text{min}})(1+3y_{\text{min}})} + 1-y_{\text{min}} \right) = 1 - y_{\text{min}}^2 + y_{\text{min}}^3 + O(y_{\text{min}}^4).$$
The second half of \ref{E:ymin:tau} now follows by combining this expression with the first half of \ref{E:ymin:tau}.
Similarly we have
$$y_{-} = \frac{1}{2}\left( - \sqrt{(1-y_{\text{min}})(1+3y_{\text{min}})} + 1-y_{\text{min}}\right) = -y_{\text{min}} + y_{\text{min}}^2 -y_{\text{min}}^3 +  O(y_{\text{min}}^4).$$
Hence
\addtocounter{theorem}{1}
\begin{equation}
\label{r:y}
r = \sqrt{y_{\text{max}}-y_{-}} = \left((1-y_{\text{min}})(1+3y_{\text{min}})\right)^{\frac{1}{4}} = 1 + \frac{y_{\text{min}}}{2} - \frac{9}{8}y_{\text{min}}^2 + O(y_{\text{min}}^3),
\end{equation}
$$k = \sqrt{\frac{y_{\text{max}}-y_{\text{min}}}{y_{\text{max}}-y_{-}}}=\sqrt{\frac{1}{2} + \frac{1-3y_{\text{min}}}{2\sqrt{(1-y_{\text{min}})(1+3y_{\text{min}})}}} = 1 - y_{\text{min}} + y_{\text{min}}^2 -3y_{\text{min}}^3 +O(y_{\text{min}}^4)$$
and
$$h:=1-k=y_{\text{min}}-y_{\text{min}}^2+3y_{\text{min}}^3 +O(y_{\text{min}}^4).$$
Combining this final expression with \ref{Kexpand} and \ref{Kpexpand} we have
\addtocounter{theorem}{1}
\begin{equation}
\label{K:y}
K(k)= \frac{1}{2}\log{\frac{8}{h}} + O\left(h\log{\frac{8}{h}}\right) = \frac{1}{2}\log{\frac{8}{y_{\text{min}}}} + O\left(y_{\text{min}}\log{\frac{8}{y_{\text{min}}}}\right)
= - \frac{1}{2}\log{y_{\text{min}}} + O(1).
\end{equation}
Hence combining \ref{r:y} and \ref{K:y}  we obtain
$$ \pt=\frac{K(k)}{r} = -\frac{1}{2}\log{y_{\text{min}}} + O(1) = -\frac{1}{2}\log{\tau} + O(1),$$
as required for \ref{E:pt}.

To prove the asymptotic expansion for the derivative in \ref{E:pt}
we first use \ref{dK:dk} to obtain
\addtocounter{theorem}{1}
\begin{equation}
\label{dp:dy}
\frac{d\pt}{dy_{\text{min}}} = \frac{d}{dy_{\text{min}}}\left(\frac{K(k)}{r}\right) =
\frac{E}{rk{k'}^2} \frac{dk}{dy_{\text{min}}} - \left( \frac{1}{rk}\frac{dk}{dy_{\text{min}}} + \frac{1}{r^2}\frac{dr}{dy_{\text{min}}} \right)K.
\end{equation}
From \ref{Eexpand} and \ref{Jpexpand} we find that the analogue of \ref{K:y} for $E$ is
\addtocounter{theorem}{1}
\begin{equation}
\label{E:y}
E(k)=1 + \frac{h}{2}\log{\frac{8}{h}} + O(h) = 1 - y_{\text{min}} \log{y_{\text{min}}} + O(y_{\text{min}}).
\end{equation}
Using the expressions for $k$, $r$, $E$ and $K$ in terms of $y_{\text{min}}$ and combining them with \ref{dp:dy} one obtains
$$ \frac{d\pt}{dy_{\text{min}}} = - \frac{(1-\tfrac{1}{2}y_{\text{min}} \log{y_{\text{min}}})}{2y_{\text{min}}} + \frac{1}{2} \left(-1 + \frac{1}{2}\right)\log{y_{\text{min}}} + O(1) = -\frac{1}{2y_{\text{min}}} + O(1).$$
The asymptotic expansion for the derivative in \ref{E:pt} 
follows using the relation between $\tau$ and $y_{\text{min}}$ given by \ref{E:tau:y}.

(iv) It follows from \ref{E:ydot} that on the interval $[0,\pt]$, $\dot{y}$
satisfies
\addtocounter{theorem}{1}
\begin{equation}
\label{y:dot:sqrt}
\dot{y} = -2y \sqrt{1-y-4\tfrac{\tau^2}{y^2}}.
\end{equation}
Hence we have 
$$ \frac{d}{dt}(e^{2t} y) = 2e^{2t}y\,\left(1 - \sqrt{1-y-4\tfrac{\tau^2}{y^2}}\right) >0,$$
that is, $e^{2t}y_{\tau}$ is increasing on $[0,\pt]$. 
In particular, putting $t=0$ we obtain 
$$ e^{2t}y_{\tau}(t) \ge y_{\tau}(0), \quad \text{for} \quad t\in [0,\pt].$$
Since for any $\tau \in \bigl[0,\tfrac{1}{3\sqrt{3}}\bigr]$
we have $y_{\text{max}}\ge \frac{2}{3}$ it now follows immediately that 
$$y_\tau(t) \ge \tfrac{2}{3}e^{-2t} \quad \text{holds on} \quad  [0,\pt].$$ 

To obtain the required $C^k$ upper bounds for $y_{\tau}$ we proceed as follows. 
Since $e^{2t}y_{\tau}$ is increasing on $[0,\pt]$, putting $t= \pt$ we obtain
$$ e^{2t}y_{\tau}(t) \le e^{2\pt}y_{\tau}(\pt) = e^{2\pt}y_{\text{min}}, \quad \text{for} \quad t\in [0,\pt].$$
It follows from \ref{E:pt} and \ref{E:tau:y} that there exists a constant $C$ independent of $\tau$
so that $e^{2\pt}y_{\text{min}} \le C$ holds for all sufficiently small $\tau>0$.
The required $C^0$ upper bound for $y_\tau$ now follows immediately.
From \ref{y:dot:sqrt} it follows that $|\dot{y}| < 2y$, and hence the required 
$C^1$ upper bound for $y_\tau$ follows from the $C^0$ upper bound for $y_\tau$.
Similarly, combining equation \ref{y:double:dot} with the $C^0$ upper bound for $y_\tau$ 
yields the required $C^2$ upper bound. By repeated differentiation of \ref{y:double:dot} 
we obtain equations for $y_{\tau}^{(k)}$ in terms of a polynomial in $y_\tau, \dot{y}_\tau, \ldots , y^{(k-2)}_\tau$.
Using these equations we can obtain inductively $C^k$ upper bounds for $y_\tau$ in terms of $C^0, C^1, \ldots C^{k-2}$
upper bounds.
\end{proof}

\subsection*{Discrete and continuous symmetries}
$\phantom{ab}$
\nopagebreak

This subsection defines various symmetries of $\Sph^1\times \R$ and $\C^3$
needed in the discussion of the special Legendrian immersions $X_\tau$.
For $x\in \R$ define the following transformations of the cylinder $\Sph^1\times \R$
\addtocounter{theorem}{1}
\begin{equation}
\label{E:sym:cyl}
\begin{aligned}
\TTT_x &:(e^{is},t) \mapsto (e^{is},t+x),\qquad
&\tbar&:(e^{is},t) \mapsto (e^{is},-t), \qquad 
&\tbar_x&:= \TTT_{2x} \circ \tbar, \\
\SSS_x &:(e^{is},t) \mapsto (e^{i(s+x)},t),\qquad 
&\sbar&: (e^{is},t) \mapsto (e^{-is},t),\qquad
&\sbar_x&:= \SSS_{2x} \circ \sbar.
\end{aligned}
\end{equation}
Note that we underline to denote a reflection.
We use a tilde to denote
isometries of $\Sph^5$,
as in the next definition.
\addtocounter{equation}{1}
\begin{definition}
\label{D:sym:tilde}
For $x\in \R$ we define $\ttilde_x$, $\stilde_x$, $\qtilde_x \in \textnormal{SU(3)}$
by taking their matrices with respect to
the standard basis $e_1, e_2, e_3$ of $\C^3$ to be
\begin{equation*}
\ttilde_x = \left(
\begin{matrix}
e^{ix} & 0 & 0 \\
0 & e^{-ix/2} & 0 \\
0 & 0 & e^{-ix/2}
\end{matrix}
\right)
\!,
\quad
\stilde_x = \left(
\begin{matrix}
1 & 0 & 0 \\
0 & \cos{x} & \sin{x} \\
0 & -\sin{x} & \cos{x}
\end{matrix}
\right)
\!,
\quad
\qtilde_x = \left(
\begin{matrix}
1 & 0 & 0 \\
0 & e^{ix} & 0 \\
0 & 0 & e^{-ix}
\end{matrix}
\right)
\!,
\end{equation*}
respectively.
We also define
$\tbartilde_x,\ \tbartilde,\  \sbartilde:\C^3 \ra \C^3$ to be orthogonal reflections
with respect to 
$\ttilde_x\left(\langle Je_1, e_2, e_3 \rangle_\R\right)=
\langle e^{ix}Je_1,e^{-ix/2}e_2, e^{-ix/2}e_3\rangle_\R$,
$\langle Je_1, e_2, e_3\rangle_\R$,
and $\langle e_1, Je_1, e_2, Je_2\rangle_\R$ respectively.
\end{definition}
The action of $\ttilde_x$ on $\Sph^5$
can be thought of as a ``translation'' along the circle
$\langle e_1, Je_1\rangle_\R\cap\Sph^5$.
The actions of both $\stilde_x$ and $\qtilde_x$ preserve the points
of this circle $\langle e_1, Je_1\rangle_\R\cap\Sph^5$.
We will refer to $\stilde_x$ as a ``rotation'' with axis this circle,
and to $\qtilde_x$ as a ``twisting'' around the same axis.
For future reference note that
$\tbartilde=\tbartilde_0$,
$\sbartilde$, 
and
$\stilde_{\pi}$,
act on
$(w_1,w_2,w_3)\in \C^3$
by 
\addtocounter{theorem}{1}
\begin{equation}
\label{E:sym:bartilde}
\begin{aligned}
\tbartilde(w_1,w_2,w_3)&=(-\overline{w}_1, \overline{w}_2, \overline{w}_3),
\\
\sbartilde(w_1,w_2,w_3)&=(w_1, w_2, -w_3),
\\
\stilde_\pi(w_1,w_2,w_3)&=(w_1, -w_2, -w_3);
\end{aligned}
\end{equation}
in particular $\stilde_{\pi}$ is orthogonal reflection with respect to
$\langle e_1, Je_1\rangle_\R$.
\addtocounter{equation}{1}
\begin{prop}
\label{P:s5:commute}
\begin{itemize}
\item[(i)]$\ttilde_x$, $\qtilde_y$, $\stilde_\pi$, $\sbartilde$ all commute with each other
and with $J$.
\item[(ii)]
$\tbartilde$ commutes with $\stilde_\pi$ and $\sbartilde$,
$\, \tbartilde \circ \ttilde_x \circ \tbartilde  = \ttilde_{-x}$,
$\,\tbartilde \circ \qtilde_x \circ \tbartilde  = \qtilde_{-x}$,
and
$\,\tbartilde\circ J=-J\circ\tbartilde$.
\item[(iii)]
$\ttilde_x$, $\qtilde_y$, $\stilde_x$ all preserve both $\Omega$ and $\omega$.
\item[(iv)]
$\tbartilde^* \Omega  = -\overline{\Omega}, \quad \tbartilde^* \omega = -\omega, \quad 
\sbartilde^* \Omega  = - \Omega, \quad \sbartilde^* \omega = \omega.$
\item[(v)]
$\tbartilde_x = \ttilde_{2x} \circ \tbartilde$.
\end{itemize}
\end{prop}
\begin{proof}
Part (i) follows from the fact that $\sbartilde$, $\ttilde_x$, $\qtilde_y$ and $\stilde_\pi$
are all diagonal matrices with respect to the basis $e_1$, $e_2$ and $e_3$.
Part (iii) is immediate from the fact that $\ttilde_x$, $\qtilde_y$, $\stilde_x$
are all $\text{SU}(3)$ matrices.
Part (ii) follows from the expressions for $\tbartilde$, $\ttilde_x$ and $\qtilde_x$ 
given in \ref{D:sym:tilde} and \ref{E:sym:bartilde}.
Part (iv) follows from the expressions
$$\omega = -\frac{1}{2i} \sum_{j=1}^3{dw_j \wedge d\overline{w}_j}, \qquad 
\Omega = dw_1 \wedge dw_2 \wedge dw_3,$$
and the expressions for $\tbartilde$ and $\sbartilde$ given in \ref{E:sym:bartilde}.
To prove part (v) it suffices to prove that $\ttilde_{2x} \circ \tbartilde$ preserves the 
vectors $e^{ix}Je_1,e^{-ix/2}e_2, e^{-ix/2}e_3$ and sends the vectors
$e^{ix}e_1$, $e^{-ix/2}Je_2$, $e^{-ix/2}Je_3$ 
into their negatives. This is easily verified using the expressions provided by
\ref{D:sym:tilde} and \ref{E:sym:bartilde}.
\end{proof}

\subsection*{The special Legendrian immersions $X_\tau$}
$\phantom{ab}$
\nopagebreak

We now proceed to define a family of special Legendrian immersions
$X_\tau: \Sph^1 \times \R \ra \Sph^5$ which have the function $y_\tau=\rho^2_\tau$
as the conformal factor for the induced metric on the cylinder. 
\begin{definition}
\label{D:X:tau}
\addtocounter{equation}{1}
For $\tau\in [0,1/3\sqrt{3}]$ define an immersion $X_\tau: \Sph^1 \times \R \ra
\Sph^5\subset \C^3$
as in \ref{E:X:tau} by
defining $w_1,w_2:\R\ra\C$
\addtocounter{theorem}{1}
\begin{equation}
\label{E:w1}
w_1(t) = 
\begin{cases}
\tanh{t}, &\text{for $\tau=0$;}\\
-i\sqrt{(1-\rho^2_{\tau})} e^{i\psi_1}, &\text{for $\tau >0$;}
\end{cases}
\qquad
w_2(t) = 
\begin{cases}
\sech{t}, &\text{for $\tau=0$;}\\
\rho_\tau e^{i\psi_2}, &\text{for $\tau >0$;}
\end{cases}
\end{equation}
and where for $\tau\in (0,1/3\sqrt{3}]$, $\psi_1,\  \psi_2:\R \ra \R$ are the odd functions defined by 
\addtocounter{theorem}{1}
\begin{equation}
\label{E:y:psi}
(1-y_{\tau}) \dot{\psi}_1 = 2\tau, \qquad y_{\tau} \dot{\psi}_2 = -2\tau, \qquad \text{with\ } \psi_1(0)=\psi_2(0)=0.
\end{equation}
\end{definition}

Define $\text{Per}(X_\tau)$, the period lattice of $X_\tau$, as follows:
\addtocounter{theorem}{1}
\begin{equation}
\label{E:Per}
 \text{Per}(X_\tau)= \{ (s,t) \in \R^2 \ | \ X_\tau \circ \SSS_{s} \circ \TTT_{t} = X_\tau \}.
\end{equation}

\addtocounter{equation}{1}
\begin{prop}
\label{P:immersion}
\ \\
(i) $X_\tau$: $\Sph^1\times \R \ra \Sph^5$ is a smooth special Legendrian immersion 
which depends smoothly on $\tau\in [0,1/3\sqrt{3})$.\\
(ii) The induced metric on $\Sph^1\times \R$ is
$\rho^2_\tau(ds^2 + dt^2)$ and has Gaussian curvature $\kappa = 1 - 8\tau^2/\rho_{\tau}^6.$
(iii) $X_\tau$ has the following decay properties:
\addtocounter{theorem}{1}
\begin{equation}
\label{E:xtau:decay}
\Vert  \partial X_\tau:\ C^k(\Sph^1 \times [0,\pt],ds^2 + dt^2, e^{-t})\Vert  \le C(k)
\end{equation}
where $\partial X_\tau$ are the partial derivatives of the coordinates of $X_\tau$ 
and $C(k)$ is some positive constant depending on $k\in \N$ but independent of $\tau$.\\
(iv)
The following equalities hold and define $\pthat$:
\addtocounter{theorem}{1}
\begin{equation}
\label{E:psi:period}
\pthat:=-\psi_2(2\pt)=-2\psi_2(\pt) = \psi_1(\pt) = \tfrac{1}{2} \psi_1(2\pt).
\end{equation}
Moreover $\pthat$ is a smooth function of $\tau$ and satisfies as $\tau\ra0$
\addtocounter{theorem}{1}
\begin{equation}
\label{E:p:hat:exp}
\,\pthat = \frac{\pi}{2} - \tau \log{\tau} + O(\tau), \qquad \frac{d\pthat}{d\tau} = -\log{\tau} + O(1).
\end{equation} 
(v) For $\tau>0$, $X_\tau$ enjoys the following symmetries
\addtocounter{theorem}{1}
\begin{subequations}
\label{E:sym}
\begin{align}
\label{E:t:sym}
\tbartilde \circ X_\tau &= X_\tau \circ \tbar,\\
\label{E:s:sym}
\sbartilde \circ X_\tau &= X_\tau \circ \sbar,\\
\label{E:sx:sym}
\stilde_x \circ X_\tau & = X_\tau \circ \SSS_x,\\
\label{E:tkp:sym}
\tbartilde_{k\pthat} \circ X_\tau & = X_\tau \circ \tbar_{k\pt}, 
\quad \text{for \ } k\in \N,\\
\label{E:tphat:sym}
\ttilde_{2k\pthat} \circ X_\tau &= X_\tau \circ \TTT_{2k\pt}, \quad \text{for \ }k\in \N.
\end{align}
\end{subequations}
(vi)
For $\tau>0$
\begin{equation}
\addtocounter{theorem}{1}
\label{E:period:k1k2}
\text{Per}(X_\tau)=
\{ (k_1\pi, k_2 2\pt) :
(k_1,k_2 \pthat )\in  (2\Z \times 2\pi \Z)\, \cup \,(\,( 2\Z+1) \times (2\Z+1)\pi\,)\},
\end{equation}
and hence $X_\tau$ factors through a torus if and only if
$\pthat \in \pi \Q$.
%
%
\end{prop}
\begin{proof} Parts (i) and (ii): It is a routine computation to verify that the immersion $X_\tau$ defined
in \ref{E:X:tau} is special Legendrian 
and that the induced metric on $\Sph^1\times \R$ is \mbox{$y_\tau\ (ds^2 + dt^2)$}.
The expression for the Gaussian curvature of $y_\tau$ given in part (ii) then follows immediately.

Smooth dependence on $\tau$: it suffices to show $w=(w_1,w_2):\R \ra \C^2$ depends 
smoothly on $\tau\in \bigl[0,\tfrac{1}{3\sqrt{3}}\bigr)$.
A calculation shows that $w$ satisfies the following system of first order complex-valued ODEs:
\addtocounter{theorem}{1}
\begin{equation}
\label{E:wdot}
\dot{w}_1 = \overline{w}_2^2, \qquad \dot{w}_2  = -\overline{w_1 w_2},
\end{equation}
and can be characterized as the unique solution to \ref{E:wdot} satisfying the initial conditions:
\addtocounter{theorem}{1}
\begin{equation}
\label{E:w:initial}
w_\tau(0)= (w_1,w_2)_{\tau}(0) =  
\begin{cases}
(0, \ 1), &\text{for $\tau=0$;}\\
(-i\sqrt{1-y_{\text{max}}}, \, \sqrt{y_{\text{max}}}\,) , &\text{for $\tau >0$;}
\end{cases}
\end{equation}
where as previously $y_{\text{max}}$ denotes the largest solution of \ref{E:cubic}.
For $\tau \in \bigl(0,\tfrac{1}{3\sqrt{3}}\bigl)$ we have that 
$y_{\text{max}} \in \bigl(\tfrac{2}{3},1\bigr)$ 
depends smoothly on $\tau^2$ and hence so does $w(0)$.
Using \ref{E:cubic} and \ref{E:w:initial} one can verify that close to $\tau=0$ 
the initial conditions for $w$ satisfy
$$w_\tau (0) = (-2i\tau,\, 1)\  + O(\tau^2).$$
Smooth dependence of $w$ on $\tau \in \bigl[0,\tfrac{1}{3\sqrt{3}}\bigr)$ now follows from the 
initial value problem characterization of $w$.

Part (iii) -- decay. It follows easily from the definitions of $w_1$ and $w_2$ and \ref{E:y:psi} that 
$$ \dot{w}_2 = \left(\frac{\dot{\rho}}{\rho} -i \frac{2\tau}{\rho^2}\right) w_2, \quad \text{and} \quad \dot{w}_1 = \left(\frac{-\rho \dot{\rho} + 2i\tau}{1-\rho^2}\right) w_1.$$
From these formulae and using the exponential decay of $\rho$ established 
in Proposition \ref{P:conformal} (iv), by induction we obtain
$$ \Vert w_2, C^k(e^{-t})\Vert  \le C(k), \quad \text{and} \quad \Vert  \dot{w}_1, C^k(e^{-2t})\Vert  \le C(k),$$
from which the statement follows.

Part (iv) -- proof of \ref{E:psi:period}.
$\psi_i(2\pt) = 2\psi_i(\pt)$ for $i=1,2$ follows directly from 
the definition of $\psi_i$ in terms of $\rho_\tau$ and the fact that $\rho_\tau$ is even about its half-period $\pt$.
It remains only to prove that $\psi_1(\pt)+2\psi(\pt)=0$. Multiplying the first equation of \ref{E:wdot} by $\overline{w}_1$
and comparing the imaginary parts of both sides leads to the equality
\addtocounter{theorem}{1}
\begin{equation}
\label{E:psi:theta}
2\tau = y \sqrt{1-y} \cos{(\psi_1+2\psi_2)}.
\end{equation}
Since $(\psi_1+2\psi_2)(0)=0$ and $y\sqrt{1-y}$ is continuous in $t$ and positive, \ref{E:psi:theta} implies
that $(\psi_1+ 2\psi_2)(t) \in (-\tfrac{\pi}{2},\tfrac{\pi}{2})$ holds for all $t\in \R$.
Then since $\dot{y}(\pt)=0$, from \ref{E:ydot} it follows that $(y \sqrt{1-y})(\pt) = 2\tau$
and hence from \ref{E:psi:theta} that $\cos(\psi_1+2\psi_2)(\pt)=1$ as required.

To prove both parts of \ref{E:p:hat:exp}
we first express $\pthat$ in terms of elliptic integrals and then utilize Appendix \ref{A:elliptic}. 
As a first step we establish that
\addtocounter{theorem}{1}
\begin{equation}
\label{E:theta1}
-\pthat = \psi_2(2\pt) = -2\left[KD(\phi,k')-KF(\phi,k') +EF(\phi,k')\right] = -\pi \Lambda_0(\phi,k) 
\end{equation}
where $\Lambda_0$ is the Heuman Lambda function defined in \ref{hlambda:defn},
$\sin{\phi}= \sqrt{\frac{y_{-}}{y_{-}-y_{\text{min}}}}$,
and $k'=\sqrt{1-k^2}$ is the complementary modulus to $k$.

From \ref{E:y:psi} we have
$\psi_2(t) = - 2 \tau\int_0^{t}{\frac{1}{y(t)}dt}$.
Using the explicit expression for $y$ in terms of the Jacobi elliptic function $\sn{}$ we obtain
\addtocounter{theorem}{1}
\begin{equation}
\label{E:theta1:gamma}
\psi_2(t) = -\frac{2\tau}{ry_{\text{max}}} \int_0^{rt}{\frac{dt}{1-\alpha^2 \sn^2{(t,k)}}} = -\frac{2\tau}{ry_{\text{max}}} \Lambda(rt,\alpha,k)
\end{equation}
where $\alpha^2 = \frac{y_{\text{max}}-y_{\text{min}}}{y_{\text{max}}}$, and $\Lambda$ is the elliptic integral of the third kind defined in \ref{ellipint3}.
Since $y_{-}<0$ for $\tau$ nonzero, we have $0<k<\alpha<1$. Hence from \ref{compellipint3} we obtain
\addtocounter{theorem}{1}
\begin{equation}
\label{compelliptint3b}
\Lambda(K,\alpha,k)= c(\alpha,k)
\left[ K(k)D(\phi,k')-K(k)F(\phi,k') + E(k)F(\phi,k') \right]
\end{equation}
where
$$c(\alpha,k) = \frac{\alpha}{\sqrt{(\alpha^2-k^2)(1-\alpha^2)}} \quad \textrm{and\ } \sin{\phi}=\frac{\sqrt{\alpha^2-k^2}}{\alpha k'}.$$
A calculation starting from the definitions of $\alpha$ and $k$ shows that
\begin{equation}
\nonumber
\frac{2\tau}{ry_{\text{max}}} c(\alpha,k) = 1 \quad \textrm{and\ } \sin{\phi} = \sqrt{\frac{y_{-}}{y_{-}-y_{\text{min}}}}.
\end{equation}
Hence \ref{E:theta1} follows from \ref{E:theta1:gamma} and \ref{compelliptint3b}.

We now derive both parts of \ref{E:p:hat:exp} from \ref{E:theta1} using Appendix \ref{A:elliptic}. 
As in the proof of Proposition \ref{P:conformal}(iii)
it is more convenient to treat the minimum value $y_{\text{min}}$ of $y(t)$ 
as the parameter instead of $\tau$.
Combining the expressions for $K$ and $E$ from \ref{K:y} and \ref{E:y} with the expressions \ref{Fexpand} and \ref{Dexpand} we find that
$$E(k)F(\phi,k') = \phi - \frac{1}{2}y_{\text{min}}\phi \log{y_{\text{min}}} + O(y_{\text{min}})$$
and
$$K(k)D(\phi,k')-K(k)F(\phi,k') =  \frac{1}{2}y_{\text{min}} \left( \phi - \sin{\phi}\cos{\phi} \right)\log{y_{\text{min}}}  + O(y_{\text{min}}).$$
\ref{E:theta1} together with the previous two expressions implies that
\addtocounter{theorem}{1}
\begin{equation}
\label{theta1expand}
-\psi_2(2\pt) = 2\phi - \sin{\phi}\cos{\phi} \ y_{\text{min}}\log{y_{\text{min}}}+ O(y_{\text{min}}).
\end{equation}
Since $\sin^2{\phi} = \frac{y_{-}}{y_{-}-y_{\text{min}}}$ using the expression for $y_{-}$ in terms of $y_{\text{min}}$ given in the proof of
Proposition \ref{P:conformal}
we have
$$\sin^2{\phi}= \frac{1}{2} \left(1-\frac{y_{\text{min}}}{2}\right) + O(y_{\text{min}}^2).$$
From this it follows that
$$\phi = \frac{\pi}{4} - \frac{1}{4}y_{\text{min}} + O(y_{\text{min}}^2) \quad \text{and} \quad {k'}^2\sin{\phi}\cos{\phi} = y_{\text{min}} + O(y_{\text{min}}^2).$$
Inserting the last two expressions into \ref{theta1expand} we obtain
\addtocounter{theorem}{1}
\begin{equation}
-\psi_2(2\pt) = \frac{\pi}{2} - \frac{1}{2}y_{\text{min}}\log{y_{\text{min}}} + O(y_{\text{min}}).
\end{equation}
The first half of \ref{E:p:hat:exp} now follows from the expression for $\tau$ in terms of $y_{\text{min}}$ given in \ref{E:tau:y}.

Similarly, to prove the second half of \ref{E:p:hat:exp} it suffices to prove
\begin{equation*}
-\frac{d\pt}{dy_{\text{min}}} = \frac{1}{2}\log{y_{\text{min}}} + O(1).
\end{equation*}
From \ref{E:theta1}, $\pthat =  \pi\Lambda_0(\phi,k)$, where $\Lambda_0$ is the Heuman Lambda function.
Applying the Chain Rule and using the expressions for the derivatives of $\Lambda_0$
given in \ref{dk:hlambda} and \ref{dphi:hlambda} we get
\addtocounter{theorem}{1}
\begin{equation}
\label{dTheta:dy}
\frac{1}{2} \frac{d\pthat}{dy_{\text{min}}} = \frac{1}{\sqrt{1-{k'}^2\sin^2{\phi}}} \left( (E - {k'}^2 \sin^2{\phi}K) \frac{d\phi}{dy_{\text{min}}} +
(E-K)\sin{\phi}\cos{\phi} \frac{1}{k}\frac{dk}{dy_{\text{min}}} \right).
\end{equation}
A short calculation shows that ${k'}^2 \sin^2{\phi} = \frac{-y_{-}}{y_{\text{max}}-y_{-}}$, and $1 - k{'}^2\sin^2{\phi} = \frac{y_{\text{max}}}{y_{\text{max}}-y_{-}}$.
Hence from our previous expansions for $y_{-}$ and $y_{\text{max}}$ in terms
of $y_{\text{min}}$ we obtain
$$ {k'}^2\sin^2{\phi} = y_{\text{min}} + O(y_{\text{min}}^2), \quad \text{and} \quad \frac{1}{\sqrt{1 - k{'}^2\sin^2{\phi}}} = 1 + \frac{1}{2}y_{\text{min}} + O(y_{\text{min}}^2).$$
Similarly from our previous expansions for $\phi$ and $k$ in terms of $y_{\text{min}}$ we have
$$ \frac{d\phi}{dy_{\text{min}}} = - \frac{1}{4} + O(y_{\text{min}}), \quad \text{and} \ \frac{dk}{dy_{\text{min}}} = -1 + O(y_{\text{min}}).$$
From \ref{E:y} we see that $E$ remains bounded as $y_{\text{min}} \ra 0$.
It follows easily that  both terms in \ref{dTheta:dy} involving $E$ remain bounded as $y_{\text{min}}\ra 0$.
Using the expansions for ${k'}^2\sin^2{\phi}$ and for $K(k)$ in terms of $y_{\text{min}}$ one sees that the first term involving $K$ in \ref{dTheta:dy}
is also bounded. Hence from \ref{dTheta:dy} we have
$$ -\frac{d\pthat}{dy_{\text{min}}} = \frac{2}{k\sqrt{1-{k'}^2\sin^2{\phi}}} \frac{dk}{dy_{\text{min}}} K \sin{\phi}\cos{\phi} + O(1) = -K + O(1)
= \frac{1}{2} \log{y_{\text{min}}} + O(1),$$
as required.

Part (v) -- symmetries. \ref{E:t:sym}, \ref{E:s:sym} and \ref{E:sx:sym} all follow 
from the basic definitions and from the facts that $\rho(t)$, $\psi_i(t)$ ($i=1,2$) are even and odd functions of $t$ respectively.\\
Proof of \ref{E:tphat:sym}: from the definitions of $X_\tau$ and $\TTT_x$ and the fact that $\rho$ has period $2\pt$ we have
$$X_\tau \circ \TTT_{2k\pt}(s,t) = X_\tau(s,t+2k\pt) = \diag(e^{ik\psi_1(2\pt)},e^{ik\psi_2(2\pt)},e^{ik\psi_2(2\pt)}) \circ X_\tau(s,t).$$
Hence \ref{E:tphat:sym} holds with $\pthat = -\psi_2(2\pt)$ using the fact from part (iv) that 
$ \psi_1(2\pt) + 2\psi_2(2\pt) = 0.$\\
Proof of \ref{E:tkp:sym}: this follows by using \ref{E:t:sym}, \ref{E:tphat:sym} and Proposition \ref{P:s5:commute}.(v).

Part (vi) -- the period lattice, $\text{Per}(X_\tau)$: it follows directly from the definition 
of the period lattice by taking $s=0$ that any $(s_0,t_0) \in \text{Per}(X_\tau)$ must satisfy 
$X_\tau(e^{is_0},t+t_0)=X_\tau(e^0,t)$ for all $t\in \R$. From \ref{E:X:tau} this is equivalent to demanding that
$$ w_1(t+t_0)=w_1(t), \quad w_2(t+t_0) = w_2(t), \quad \text{and}\ -w_2(t+t_0) \sin{s_0} = 0, \quad \text{for all } t\in \R.$$
Hence $s_0 \in \pi \Z$ and $t_0 \in 2\pt \Z$ as required. If $(s_0,t_0) = (k_1\pi,k_2 2\pt)$ for some $k_1, k_2 \in \Z$, then if follows
from \ref{E:sx:sym} and \ref{E:tphat:sym} that $(s_0,t_0) \in \text{Per}(X_\tau)$ if and only if 
$$ \ttilde_{2k_2\pthat} \circ \stilde_{k_1 \pi} \circ X_\tau = X_\tau.$$
It is easy to check that this is equivalent to the condition that $\ttilde_{2k_2\pthat} \circ \stilde_{k_1 \pi} =Id$. 
Using \ref{D:sym:tilde} it is also straightforward to verify that the previous condition 
leads to \ref{E:period:k1k2}a if $k_1$ is even and to \ref{E:period:k1k2}b if $k_1$ is odd. 
\end{proof}

\addtocounter{equation}{1}
\begin{remark}
\label{R:negative:tau}
The parameter $\tau$ can be extended to negative values to give a family
which depends smoothly on $\tau\in(-1/3\sqrt{3},1/3\sqrt{3})$
by defining for $\tau\in(-1/3\sqrt{3},0)$
$X_\tau:=\overline{X_{-\tau}}$.
Since though the new immersions are simply the complex conjugates of the ones we already have
this is of no use.
The proof of the smoothness is by using \ref{E:wdot} and \ref{E:w:initial},
expressing $\tau$ as a smooth function of $w_1(0)$ explicitly and then inverting
to express the initial data as smooth functions of $\tau$.
\end{remark}

\addtocounter{equation}{1}
\begin{remark}
\label{R:catenoid}
It is easy to calculate that the negatively curved regions of $X_\tau$
magnified by a factor of order $1/\sqrt{\tau}$---see \ref{E:ymin:tau}---and
appropriately translated,
tend as $\tau\to0$ to the so-called Lagrangian catenoid
which is a special Lagrangian surface in $\C^2$ given by 
$$
|w_1|^2-|w_2|^2=0,
\qquad
\Imag(w_1w_2)=1.
$$
Note the analogy with the negatively curved regions of the Delaunay surfaces
which tend to catenoids
\cite{kapouleas:annals,kapouleas:imc,kapouleas:jdg,kapouleas:annals,kapouleas:bulletin,kapouleas:survey}.
For the current construction though we do not need this result and so we will not present its proof here.
\end{remark}

\subsection*{Twisting $X_\tau$ to $\Xta$}
$\phantom{ab}$
\nopagebreak

Recall that by \ref{E:tphat:sym} $X_\tau$ has a translational period.
As we mentioned in the introduction we need to use
modified versions of $X_{\tau}$ so that the period involves a prescribed amount
of twisting $\qtilde_x$ as well.
In what follows we introduce the desired twisting to convert the special Legendrian
immersion $X_\tau$ to a Legendrian immersion $\Xta$.
$\Xta$ is corrected to a special Legendrian immersion
$(X_{\tau,\alpha})_{\fundertilde}$
in \ref{T:Xta}.
$(X_{\tau,\alpha})_{\fundertilde}$
is no longer  $\text{U}(1)$-invariant,
and hence is not in the family of immersions constructed in \cite{haskins:thesis,haskins:slgcones}.

We define now---for $\tau$ and $|\alpha|$ small---$X_{\tau,\alpha}: \Sph^1 \times [-\pt,\pt] \ra \Sph^5$  by requiring 
\addtocounter{theorem}{1}
\begin{equation}
\label{E:xta:0}
\begin{aligned}
\tbartilde \circ X_{\tau,\alpha} = X_{\tau,\alpha} \circ \tbar, \qquad & \text{on\quad}
\Sph^1 \times [-\pt,\pt],\\
\Xta=
\join{[X_{\tau}, \qtilde_\alpha \circ X_{\tau};1,2;e_1,e_2,e_3]} \qquad & \text{on\quad } \Sph^1 \times [-1,\pt].
\end{aligned}
\end{equation}
$\Xta$ then agrees with 
$\qtilde_\alpha \circ X_{\tau}$ on $\Sph^1 \times [2,\pt]$
and by Proposition \ref{P:s5:commute}(ii), with 
$\qtilde_{-\alpha} \circ X_{\tau}$ on $\Sph^1 \times [-\pt,-2]$.
We extend the map $X_{\tau,\alpha}$ to $\Sph^1 \times \R$ by defining
\addtocounter{theorem}{1}
\begin{equation}
\label{E:xta}
X_{\tau,\alpha}=  \ptilde_{\tau,\alpha}^j \circ X_{\tau,\alpha} \circ \TTT_{-2j\pt} \ \text{on} \quad \Sph^1\times[(2j-1)\pt,(2j+1)\pt] \ \text{for\ } j\in\Z,
\end{equation}
where
(recall \ref{P:s5:commute}.i)
\begin{equation}
\addtocounter{theorem}{1}
\label{E:xta:period}
\ptilde_{\tau,\alpha}
:=
\ttilde_{2\ptphat} \circ \qtilde_{2 \alpha}
=
\qtilde_{2 \alpha} \circ \ttilde_{2\ptphat} 
\end{equation}
is the new ``translating-twisting'' period
which reduces to the period $\ttilde_{2\ptphat}$ of $X_\tau$
when $\alpha=0$.
Most but not all of the symmetries of $X_\tau$
(recall \ref{E:sym})
generalize to $X_{\tau,\alpha}$ as follows:

\addtocounter{equation}{1}
\begin{lemma}
\label{L:a-sym}
$X_{\tau,\alpha}:\Sph^1 \times \R\ra\Sph^5$ satisfies the symmetries 
\addtocounter{theorem}{1}
\begin{subequations}
\begin{align}
\label{E:t:a-sym}
\tbartilde \circ X_{\tau,\alpha} &= X_{\tau,\alpha} \circ \tbar,\\
\label{E:s:a-sym}
\sbartilde \circ X_{\tau,\alpha} &= X_{\tau,\alpha} \circ \sbar,\\
\label{E:sx:a-sym}
\stilde_\pi \circ X_{\tau,\alpha} &= X_{\tau,\alpha} \circ \SSS_\pi,\\
\label{E:tkp:a-sym}
\ptilde_{\tau,\alpha}^k \circ\tbartilde \circ \Xta & = \Xta \circ \tbar_{k\pt}, 
\quad \text{for \ } k\in \N,\\
\label{E:tphat:a-sym}
\ptilde_{\tau,\alpha}^k \circ X_{\tau,\alpha} &= X_{\tau,\alpha} \circ \TTT_{2k\pt},
\quad \text{for \ }k\in \N.
\end{align}
\end{subequations}
\end{lemma}

\begin{proof}
Using \ref{E:sym}, the definitions, 
proposition \ref{P:s5:commute}(i,ii,iv),
and $\tbar_{k\pt}=\TTT_{2k\pt}\circ\tbar$,
it is straightforward to check the symmetries above.
\end{proof}

\section{Construction of the initial surfaces}
\label{S:init:surf}
\nopagebreak

\subsection*{Introduction}
$\phantom{ab}$
\nopagebreak

The idea of the construction of the initial immersions is to consider
$g$ copies of $X_\tb$ suitably repositioned relative to each other,
and fuse them together to obtain a Legendrian immersion $Y_\zerobold:M\ra \Sph^5$,
where $M$ is a closed surface of genus $g$.
$\tau=\tb$ is chosen in \ref{E:p:hat}
to ensure that $X_\tb$ factors through a torus
containing a prescribed number of fundamental regions
(see \ref{L:Per}).
In order to achieve the required decay estimates in the construction later,
a two-parameter family of Legendrian perturbations
$Y_{\zetabold}:M\ra \Sph^5$ of $Y_{\zerobold}$ is needed.
In the construction of $Y_{\zetabold}$, $X_\tb$ is substituted
by $\Xta$ suitably placed,
where $\tau$ and $\alpha$ are determined in \ref{E:sigma}
in terms of $\tb$ and the parameters $\zetabold=(\zeta_1,\zeta_2)\in\R^2$.
 
The repositioning uses $\rtilde\in\text{SU(3)}$
(note that ${{{\rtilde}}}^g$ is the identity),
defined by taking its matrix with respect to the
unitary basis $e_1, e_2, e_3$ to be
\addtocounter{theorem}{1}
\begin{equation} 
\label{E:rotn}
\rtilde:= \left(
  \begin{matrix}
    \cos{\frac{2\pi}{g}} & \sin{\frac{2\pi}{g}} & 0 \\
    -\sin{\frac{2\pi}{g}} & \cos{\frac{2\pi}{g}} & 0 \\
    0 & 0 & 1
  \end{matrix}
\right).
\end{equation}

$\tb$ is defined to be the unique small number (by \ref{E:p:hat:exp})
which satisfies 
\addtocounter{theorem}{1}
\begin{equation}
\label{E:p:hat}
\ptbhat = \frac{\pi}{2}\left(1+ \frac{1}{m}\right),
\quad\qquad
\text{
where $m=4\mb-1$,}
\end{equation}
where from now on $\mb\in\N$ is assumed 
fixed and as large as it may be needed.
$\tb$ is then fixed and as small as it may be needed.
We have then 
\addtocounter{equation}{1}
\begin{lemma}
\label{L:Per}
$\text{Per}(X_\tb) = \vspan_{\Z}{\langle(2\pi,0),(0,2m\ptb)\rangle}$,
$\ttilde_{m\ptbhat}=\ttilde_{-m\ptbhat}$,
$\tbartilde_{m\ptbhat}=\tbartilde_{-m\ptbhat}= \tbartilde$,
and $X_\tb$ factors through an embedding of the torus
$\Sph^1 \times (\R / 2m\ptb\Z)$.
\end{lemma}
\begin{proof} From Proposition \ref{P:immersion}(vi) any period is of the form $(k_1\pi,2k_2\ptb)$ for some $k_1, k_2 \in \Z$.
First we show that there are no periods of the form $(k_1\pi,2k_2\ptb)$, when $k_1$ is odd.
In this case, according to \ref{E:period:k1k2}b, $(k_1\pi,2k_2\ptb)$ is a period if and only if 
$$ 2\mb k_1 = (4\mb-1)(2l+1) \quad \text{for some\ } l\in \Z.$$
However, the left-hand side of the previous equation is even while the right-hand side is odd.
Hence there are no solutions to this equation.

If $k_1$ is even, then according to \ref{E:period:k1k2}a, $(k_1\pi,2k_2\ptb)$ is a period if and only if
$$ k_1 \mb = l(4\mb-1), \quad \text{for some\ } l\in \Z.$$
Hence $l=0 \mod{\mb}$, from which we see that the previous equation has solutions exactly when $k_1 \in (4\mb-1)\Z = m\Z$.


\ref{E:p:hat} implies that $m\ptbhat=2\mb\pi$ and the equalities follow from
\ref{D:sym:tilde} and  \ref{P:s5:commute}(v).

Embeddedness: we want to show that $X(s_1,t_1) = X(s_2, t_2)$ implies that
$(s_1-s_2,t_1-t_2) \in \textrm{Per}(X)$. 
Using the definition of $X$ and equation \ref{E:wdot} one can see that the
immersions $X$ satisfy first-order systems of ODEs in both $s$ and $t$, namely
\begin{align*}
X_s & = AX,\\
X_t & = X \times AX,
\end{align*}
where $A=\tfrac{d}{dx}\bigr|_{x=0}\stilde_x$ and the cross-product $u\times v$ 
of two vectors in $\C^3$ is defined to be
$$ u\times v = (\bar{u}_2 \bar{v}_3 - \bar{u}_3 \bar{v}_2, \bar{u}_3 \bar{v}_1 - \bar{u}_1 \bar{v}_3, 
\bar{u}_1 \bar{v}_2 - \bar{u}_2 \bar{v}_1 ).$$
Suppose $X(s_1,t_1) = X(s_2, t_2)$. Then by uniqueness of solutions to these first-order ODEs, we have that 
$X(s+s_1,t+t_1) = X(s+s_2,t+t_2)$ for all $s$ and $t$. Hence $(s_1-s_2,t_1-t_2) \in \textrm{Per}(X)$
as required.
\end{proof}

The range of the parameters $\zetabold=(\zeta_1,\zeta_2)\in\R^2$ is determined by
\begin{equation}
\addtocounter{theorem}{1}
\label{E:sigma:range}
|\zeta_1|,|\zeta_2|\le\cunder\,\tb,
\end{equation}
where $\cunder>0$ will be chosen later
in the proof of
\ref{T:main}.
$\zeta'_1$, $\zeta'_2$, $\tau$, and $\alpha$ are determined by 
\addtocounter{theorem}{1}
\begin{equation}
\label{E:sigma}
\zeta'_1:=c'_1\zeta_1= m(\ptbhat - \ptphat),
\qquad
\zeta'_2:=c'_2\zeta_2 = (1-m)\alpha,
\end{equation}
where $c'_1,c'_2$ are normalization constants determined in \ref{E:c1c2}.
Note that (see \ref{L:sliding:twisting:X}) the total `sliding' introduced by using $\Xta$
instead of $X_\tb$ is $\zeta'_1$ and the total twisting $\zeta'_2$.
Using now \ref{E:p:hat:exp} and \ref{P:conformal}(iii)
it clearly follows that $\tau$ and $\alpha$ depend smoothly on $\zetabold$ and
\begin{equation}
\addtocounter{theorem}{1}
\label{E:tau:alp:interval}
|\tau-\tb|\le C \cunder \tb^2, \qquad
\left| 1-\frac{\pt}{\ptb}\right|\le \frac{C\cunder\tb}{-\log\tb},\qquad
|\alpha| \le C \cunder c'_2 \tb^2 \,|\log\tb|.
\end{equation}

\subsection*{The parametrizing surface $M$ and symmetries}
$\phantom{ab}$
\nopagebreak

In this subsection we describe and define the abstract surface $M$.
$M$ is independent of the parameters $\zetabold$ and actually is also
independent of $\mb$ and $\tau$ since it is simply a closed surface of genus $g$.
We define however $M$ together with convenient coordinate patches
useful in the definition of $Y_{\zetabold}$ later.
These coordinates depend on $\mb$ and $\tb$.
$M$ is defined as follows:
Consider $M_0$, a two-sphere with $2g$ discs removed,
and $g$ copies $M_1,...,M_g$
of a cylinder parametrizing the appropriate part of $X_\tb$ (see \ref{D:Mprime}).
$M$ is obtained by identifying
neighborhoods of the $2g$ boundary circles of $M_0$ with
neighborhoods of the $2g$ boundary circles of the cylinders (see \ref{D:M}).

It is convenient to identify $\R^3$ with the
real $3$-plane
$\langle e_1, e_2, e_3\rangle_\R$.
$\Sph^2$ is then identified with the intersection of this 3-plane with $\Sph^5\subset\C^3$.
We can think then of $X_0$ defined in \ref{D:X:tau} as a diffeomorphism
of the cylinder $\Sph^1\times\R$ onto $\Sph^2\setminus\{e_1,-e_1\}$,
and hence we can use it to describe the identifications of annuli on the cylinders
with annuli on $M_0$.
To specify the annuli of identification in \ref{D:M}
we define constants $\delta, a$ by
\addtocounter{theorem}{1}
\begin{equation}
\label{E:delta:a}
\delta = \frac{\pi}{100g},  \qquad \sech{(a+1)}=\sin{\delta},
\end{equation}
so that 
$X_0(\Sph^1 \times \{a+1\})$ is
the circle on $\Sph^2$ with center $e_1$ and  geodesic radius $\delta$.






The centers of the discs excised from $\Sph^2$ to define $M_0$
form a canonical $2g$-gon on
$\langle e_1,e_2\rangle_\R$.
The set of its vertices is:
\addtocounter{theorem}{1}
\begin{equation} 
\label{E:V}
V: = \left\{\left(\cos{\frac{j\pi}{g}},\sin{\frac{j\pi}{g}},0\right): j\in \Z\right\}
\end{equation}
We define then $V_\delta$ to be 
the set of points of $\Sph^2$ whose geodesic distance from 
$V$ is less than $\delta$,
so 
$V_\delta$ is the union of $2g$ discs of radius $\delta$ and centers
the vertices of the regular polygon.
Note that
$\rtilde$, $\tbartilde$, $\sbartilde$, and $ \stilde_\pi$,
restrict to isometries of $\Sph^2$ which also preserve $V$ and $V_\delta$.
Since $\Sph^2$ is used in the construction of $M$,
we drop the tilde and denote these restrictions to $\Sph^2$ by 
$\RRR$, $\tbar$, $\sbar$, and $\SSS_\pi$ respectively.
Note that they are a rotation by an angle $2\pi/g$ on the real $e_1e_2$-plane,
reflection with respect to the real $e_2e_3$-plane,
reflection with respect to the real $e_1e_2$-plane,
and reflection with respect to the real $e_1$-line respectively,
therefore they satisfy
\addtocounter{theorem}{1}
\begin{equation}
\label{E:commute}
 \tbar \circ \RRR \circ \tbar = \RRR^{-1},
\qquad  \sbar \circ \RRR \circ \sbar = \RRR,
\qquad \SSS_\pi \circ \RRR \circ \SSS_\pi = \RRR^{-1}.
\end{equation}


\addtocounter{equation}{1}
\begin{definition}
\label{D:Mprime}
We define $M' = \coprod_{j=0}^g  M_j$
where 
$$ M_j = 
\begin{cases}
\{j\} \times \Sph^1 \times [a,2m\ptb-a], & \text{for $j=1, \ldots g$;}\\
 \Sph^2 \setminus V_{\delta}, & \text{for $j=0$.}
\end{cases}
$$
\end{definition}

%

The relations \ref{E:commute} imply that on $M_0$
\addtocounter{theorem}{1}
\begin{equation}
\label{E:commute2}
\tbar = \RRR^{-j} \circ \tbar \circ \RRR^{-j}, \qquad
\sbar = \RRR^j \circ \sbar \circ \RRR^{-j}, \qquad
\SSS_\pi = \RRR^{-j} \circ \SSS_\pi \circ \RRR^{-j}.
\end{equation}
We extend
$\RRR$, $\tbar$, $\sbar$, and $\SSS_\pi$
to diffeomorphisms of $M'$ by requiring
$$\RRR(j,q) := (j',q) \text{\ for\ } q\in \Sph^1 \times [a,2m\ptb-a]$$ 
for any 
$j,j'=1, \ldots ,g \text{\ with\ } j'=j+1$ $(\!\!\mod g)$,
by defining for any 
$(1,q) \in M_1$ (recall \ref{E:sym:cyl}),
$$
\tbar(1,q) := (1,\tbar_{m\ptb}(q)), \quad
\sbar(1,q) := (1, \sbar(q)), \quad
\SSS_\pi(1,q) := (1,\SSS_\pi(q)),
$$
and then on $M_j$ ($j=2,...,g$)
$$
\tbar := \RRR^{1-j} \circ \tbar \circ \RRR^{1-j}, \qquad
\sbar := \RRR^{j-1} \circ \sbar \circ \RRR^{1-j}, \qquad
\SSS_\pi := \RRR^{1-j} \circ \SSS_\pi \circ \RRR^{1-j}.
$$

\addtocounter{equation}{1}
\begin{lemma}
\label{L:Mp:sym}
$\tbar$, $\sbar$, and $\SSS_\pi$ are of order $2$ on $M'$,
they commute with each other,
and \ref{E:commute} and \ref{E:commute2} are still valid on $M'$.
\end{lemma}

\begin{proof}
It is straightforward to check these relations by using the definitions.
\end{proof}

To simplify the notation we identify $M_1$ with a subset of $\Sph^1\times\R$
by identifying $(1,e^{is},t)$ with $(e^{is},t)$.
$X_0$ then maps $M_1$ to $\Sph^2$ and is given (recall \ref{D:X:tau}) by
$$
X_0(1,e^{is},t):=\tanh t\,\,e_1
+
\sech t\,\, \stilde_s(e_2).
$$
Clearly then on $M_1$
\addtocounter{theorem}{1}
\begin{equation}
\label{E:X0:sym}
\sbar\circ X_0=X_0\circ \sbar,
\qquad
\SSS_\pi\circ X_0=X_0\circ\SSS_\pi.
\end{equation}
Note also that by \ref{E:delta:a} $X_0$ maps
$ \{1\} \times \Sph^1 \times [a,a+1]$
to a neighborhood of a boundary circle of $M_0$
as needed in the following definition.

\addtocounter{equation}{1}
\begin{definition}
\label{D:M}
We define 
$M:= M' / \sim$ where
$\sim$ is induced by the identifications given by
$$
\RRR^{j}\circ X_0(p)\sim \RRR^j(p),
\qquad
\RRR^{j}\circ \tbar \circ X_0(p)\sim \RRR^j \circ \tbar(p),
$$
for $j=0,...,g-1$ and $p\in
 \{1\} \times \Sph^1 \times [a,a+1]$.
By an abuse of notation we will identify $M_j$ with its image in $M$ under the identification.
The standard coordinates on each $M_j$
($j\in \{1, \ldots, g\}$)
will be denoted by $s$ (modulo $2\pi$) and $t$.
\end{definition}

Clearly with these identifications $g$ handles are attached to $M_0$
and $M$ is an oriented closed surface of genus $g$.


Using \ref{L:Mp:sym} and  \ref{E:X0:sym},
it is straightforward to check that 
$\RRR$, $\tbar$, $\sbar$, and $\SSS_\pi$,
respect the identifications and induce diffeomorphisms of $M$ which we will denote 
by the same symbols.
Moreover, we have the following:

\addtocounter{equation}{1}
\begin{lemma}
\label{L:M:sym}
$\tbar$, $\sbar$, and $\SSS_\pi$ are of order $2$ on $M$,
they commute with each other,
and \ref{E:commute} and \ref{E:commute2} are still valid on $M$.
\end{lemma}

\begin{proof}
It is straightforward to check these relations by using the definitions and \ref{L:Mp:sym}.
\end{proof}

\addtocounter{equation}{1}
\begin{definition}
\label{D:group}
We denote by $\group$ the group of diffeomorphisms
of $M$ generated by $\tbar$, $\sbar$, $\SSS_\pi$ and $\RRR$.
We call $\group$ the group of symmetries of $M$. 
\end{definition}


\subsection*{Standard and transition regions}
$\phantom{ab}$
\nopagebreak

Motivated by the definition of $Y_{\zerobold}$ later in which $X_{\tb}$ is used (see \ref{D:Yzeta}),
we identify various regions on $M$ in the usual fashion of 
\cite{kapouleas:annals,kapouleas:finite,kapouleas:minimal,kapouleas:imc,kapouleas:wente,kapouleas:wente:announce,kapouleas:jdg,kapouleas:cmp,kapouleas:annals,kapouleas:bulletin}
as follows:
For
\addtocounter{theorem}{1}
\begin{equation}
\label{E:nnprime}
n\in \{0, 1, \ldots ,2\mb-1\},
\quad n'\in \{1, 2, \ldots ,2\mb-1\},
\quad\text{ and }\quad n''\in \{1, 2, \ldots ,2\mb\},
\end{equation}
we define
\addtocounter{theorem}{1}
\begin{subequations}
\begin{align} 
S_x[0] &:= M_0 \cup \group (\{1\}\times \Sph^1 \times [a, b+x]),\\
\widetilde{S}_x[0] &:= M_0 \cup \group (\{1\} \times \Sph^1 \times [a,2\ptb -b-x]),\\
S_x[n'] &:= \{1\} \times \Sph^1 \times [2n'\ptb - b - x, 2n'\ptb +b +x],\\
\widetilde{S}_x[n'] &:= \{1\} \times \Sph^1 \times [(2n'-2)\ptb + b + x, (2n'+2)\ptb -b -x],\\
\label{E:neck}
\Lambda_{x,y}[n''] & := \{1\} \times \Sph^1 \times [(2n''-2)\ptb + b + x, 2n''\ptb -b -y],\\
C^+_x[n] &:= \{1\} \times \Sph^1 \times \{2n\ptb+b+x\},\\
C^-_x[n''] &:= \{1\} \times \Sph^1 \times \{2n''\ptb-b-x\},\\
\label{E:Ccentral}
\widetilde{C}[n''] &:= \{1\} \times \Sph^1 \times \{(2n''-1)\ptb\},
\end{align}
\end{subequations}
where $b>a+5$ is a constant independent of $\tau$ chosen in \ref{E:conformaldecay},
$x,y\ge0$ and $\tb$ is assumed small enough so that $x,y< \ptb -b$.
When $x=y=0$ we will drop the subscripts and refer simply to $S[0], \ldots , \Lambda[n'']$.
We also write $\Lambda_x[n'']$ for $\Lambda_{x,x}[n'']$.

Note that
\addtocounter{theorem}{1}
\begin{equation}
\label{E:SL}
\begin{gathered}
\widetilde{S}_x[0]=S_x[0]\cup\group\Lambda_x[1],
\qquad
\widetilde{S}_x[n']=\Lambda_x[n']\cup S_x[n']\cup\Lambda_x[n'+1],
\\
C^+_x[n]=S_x[n]\cap\Lambda_x[n+1],
\qquad
C^-_x[n']=S_x[n']\cap\Lambda_x[n'],
\qquad
C^-_x[2\mb]= \tbar C_x^+[2\mb-1],
\\
\partial S_x[0]=\group C^+_x[0],
\qquad
\partial S_x[n']=C^-_x[n']\cup C^+_x[n'],
\qquad
\partial \Lambda_{x,y}[n''] = C_x^+[n''-1] \cup \, C^-_y[n''].
\end{gathered}
\end{equation}

$S_x[0]$ and $\widetilde{S}_x[0]$ are invariant under $\group$, while the other regions
are invariant only under $\sbar$ and $\SSS_\pi$, except for 
$\Lambda_x[2\mb]$ which is also invariant under $\tbar$.
The regions 
$$S_x[0], \ \RRR^j S_x[n'], \ \RRR^j \tbar S_x[n'], \ \RRR^j \Lambda_x[n'], \ \RRR^j \tbar \Lambda_x[n'] \text{\ and\ } \RRR^j \Lambda_x[2\mb]$$  where $j\in \{1, \ldots ,g\}$ and
$n'\in \{1,\ldots ,2\mb-1\}$
provide a decomposition of $M$ with overlaps only on their boundary circles.

$S[n]$ (and their images under elements of $\group$)
are \textit{standard regions} \cite{kapouleas:survey} and 
\textit{almost spherical regions} in the terminology of
\cite{kapouleas:annals,kapouleas:finite,kapouleas:minimal,kapouleas:imc,kapouleas:wente,kapouleas:wente:announce,kapouleas:jdg,kapouleas:cmp,kapouleas:annals,kapouleas:bulletin}.
Similarly, 
$\Lambda[n]$ (and their images under $\group$) are \textit{transition regions} or alternatively 
\textit{necks}, while 
$\widetilde{S}[n]$ are \textit{extended standard regions} or \textit{extended almost spherical regions} 
and are the union of $S[n]$ with the adjacent
transition regions. $S[0]$ is called the \textit{central almost spherical region}
and is the region where the fusion of the constituent tori occurs.
We call $\Lambda[2\mb]$ the \textit{opposing transition region} or \textit{opposing neck}. 
It is the neck which is farthest away from the central almost spherical region of $M$.
Each $S_x[n]$ is a neighborhood of $S[n]$, while $\widetilde{S}_x[n]$ is
$\widetilde{S}[n]$ with an appropriate neighborhood of its boundary excised.

\subsection*{The initial immersions}
$\phantom{ab}$
\nopagebreak

We will only consider Legendrian immersions $Y:M \ra \Sph^5$ which satisfy the symmetries
\addtocounter{theorem}{1}
\begin{equation}
\label{E:symmetries}
 \tbartilde \circ Y  = Y \circ \tbar, \qquad \sbartilde \circ Y = Y \circ \sbar, \qquad \stilde_\pi \circ Y = Y \circ \SSS_\pi, \qquad \rtilde \circ Y = Y \circ \RRR.
\end{equation}

We first observe that given suitable immersions of $M_0$ and the cylinder
we can use the symmetries to generate an immersion of $M$:
\addtocounter{equation}{1}
\begin{lemma}
\label{L:Y} 
Given an immersion $\check{X}_0:M_0 \ra \Sph^5$ satisfying \ref{E:symmetries} on $M_0$
and an immersion $\check{X}: \Sph^1 \times [a, m\ptb +1] \ra \Sph^5$ such that
\begin{enumerate}
\item[(i)]
$\check{X} = \check{X}_0 \circ X_0\ $ on $\,\Sph^1 \times [a,a+1]$,
\item[(ii)]
$\tbartilde \circ \check{X} = \check{X} \circ \tbar\ $ on $\,\Sph^1 \times [m\ptb-1,m\ptb+1]$,
\item[(iii)]
$\sbartilde \circ \check{X} = \check{X} \circ \sbar$
and
$\stilde_\pi \circ \check{X} = \check{X} \circ \SSS_\pi$ on $\,\Sph^1 \times [a, m\ptb +1]$,
\end{enumerate}
there is a unique immersion $\check{Y}:M \ra \Sph^5$ satisfying \ref{E:symmetries}
and such that 
\begin{equation*}
\check{Y}=
\check{X}_0  \quad \text{on} \quad  M_0, \quad 
\check{Y}(1,q) = \check{X}(q) \quad \text{for} \quad q\in \Sph^1 \times [a, m\ptb +1].
\end{equation*}
\end{lemma}
\begin{proof}
We first extend $\check{X}$ to $\Sph^1\times[a,2m\ptb-a]$
by defining
$$
\check{X}=\tbartilde \circ \check{X}\circ\tbar
\qquad
\text{on}
\qquad
\Sph^1\times[m\ptb,2m\ptb-a],
$$
which is consistent by (ii).
We then define $\check{Y'}:M'\to\Sph^5$ by (recall the identification of $M_1$ with
$\Sph^1\times[a,2m\ptb-a]$):
$$
\begin{aligned}
\check{Y'}&=\check{X}_0 \qquad &\text{on}\qquad &M_0
\\
\check{Y'}&=\rtilde^{j-1}\circ\check{X}\circ\RRR^{1-j} \qquad &\text{on}\qquad &M_j
\quad
(j\in\N).
\end{aligned}
$$
It is straightforward then to check that $\check{Y'}$ satisfies the symmetries in \ref{E:symmetries}
and then factors through $\check{Y}:M\to\Sph^5$ which also satisfies \ref{E:symmetries}
and the other required conditions.
\end{proof}

To construct the desired immersions $Y_{\zetabold}$
we intend to apply \ref{L:Y} with
$\check{X}_0:M_0\ra\Sph^2\subset\Sph^5$
simply the inclusion map,
and $\check{X}$ an appropriate modification of $X_{\tau,\alpha}$.
We discuss now the various modifications of $X_{\tau,\alpha}$
so the conditions required in \ref{L:Y} are satisfied.
First, we need to reparametrize $\Xta$ to accommodate for the dependence of $\pt$ on $\tau$:
We define $\Xunderta:\Sph^1\times\R\to\Sph^5$ by 
\addtocounter{theorem}{1}
\begin{equation}
\label{E:Xunderta}
\Xunderta(e^{is},t)=
\Xta(e^{is},\underline{t}(t)),
\qquad
\text{ where }
\qquad
\underline{t}(t)=\frac{\pt}{\ptb} t.
\end{equation}

The following lemma specifies the symmetries of $\Xunderta$
and the appropriate sliding and twisting required by \ref{L:Y}.(ii):

\addtocounter{equation}{1}
\begin{lemma}
\label{L:sliding:twisting:X}
(i). 
$\Xunderta:\Sph^1 \times \R\ra\Sph^5$ satisfies the symmetries 
\addtocounter{theorem}{1}
\begin{subequations}
\begin{align}
\label{E:t:a:under-sym}
\tbartilde \circ \Xunderta &= \Xunderta \circ \tbar,\\
\label{E:s:a:under-sym}
\sbartilde \circ \Xunderta &= \Xunderta \circ \sbar,\\
\label{E:sx:a:under-sym}
\stilde_\pi \circ \Xunderta &= \Xunderta \circ \SSS_\pi,\\
\label{E:tkp:a:under-sym}
\ptilde_{\tau,\alpha}^k \circ\tbartilde \circ \Xunderta & = \Xunderta \circ \tbar_{k\ptb}, 
\quad \text{for \ } k\in \N,\\
\label{E:tphat:a:under-sym}
\ptilde_{\tau,\alpha}^k \circ \Xunderta &= \Xunderta \circ \TTT_{2k\ptb},
\quad \text{for \ }k\in \N.
\end{align}
\end{subequations}
(ii). 
$
\tbartilde
\circ
\ttilde_{\zeta'_1}
\circ
\qtilde_{\zeta'_2-\alpha}
\circ
\Xunderta 
=
\ttilde_{\zeta'_1}
\circ
\qtilde_{\zeta'_2-\alpha}
\circ
\Xunderta 
\circ
\tbar_{m\ptb}
$
on
$\Sph^1\times\R$.
\newline
(iii). 
$ \qtilde_{-\alpha} \circ X_{\tau,\alpha}=X_\tau$
on $\Sph^1\times[2,2\pt-2]$.
\end{lemma}

\begin{proof}
(i). The symmetries follow from
definition \ref{E:Xunderta} and
the symmetries of $\Xta$ in \ref{L:a-sym}.
(ii). Follows from \ref{E:tkp:a:under-sym}, \ref{L:Per}, and the commutation laws
in \ref{P:s5:commute}.
\newline
(iii). Follows from the definition of $X_{\tau,\alpha}$.
\end{proof}

The last lemma motivates us to define
$\widetilde{X}_{\tau,\alpha}: \Sph^1 \times \R \ra \Sph^5$ by 
\addtocounter{theorem}{1}
\begin{equation}
\widetilde{X}_{\tau,\alpha} = 
\begin{cases}
\join{[X_0,\ttilde_{\zeta'_1} \circ\qtilde_{\zeta'_2} \circ X_0;a+1,a+2;e_1,e_2,e_3]}
& \text{on} \quad \Sph^1 \times (-\infty,a+2],\\
\ttilde_{\zeta'_1} \circ\qtilde_{\zeta'_2} \circ
\join{[X_0,\qtilde_{-\alpha}\circ \Xunderta;a+2,a+3;e_1,e_2,e_3]}
& \text{on} \quad \Sph^1 \times [a+2,\infty).
\end{cases}
\end{equation}

Note that
$\widetilde{X}_{\tau,\alpha}$ transits from $X_0$ to a translated (by  $\zeta'_1$)
and twisted (by $\zeta'_2$)
$X_0$,
and then transits to the appropriately 
(by \ref{L:sliding:twisting:X}.i)
translated and twisted $\Xunderta$,
that is
$\widetilde{X}_{\tau,\alpha} =
\ttilde_{\zeta'_1}
\circ
\qtilde_{\zeta'_2-\alpha}
\circ
\Xunderta$
on $\Sph^1 \times [a+3,\infty)$.


\addtocounter{equation}{1}
\begin{definition}
\label{D:Yzeta}
We now define $Y_{\zetabold}=\check{Y}:M \ra \Sph^5$ by applying Lemma \ref{L:Y}
with
\begin{enumerate}
\item[(a)]
$\check{X}_0:M_0\ra\Sph^2\subset\Sph^5$
the inclusion map.
\item[(b)]
$\check{X}: \Sph^1 \times [a_1-1,m\pt+1] \ra \Sph^5$ to be 
the restriction of
$ \widetilde{X}_{\tau,\alpha}$.
\end{enumerate}
\end{definition}


\subsection*{Remarks on the geometry and the Lagrangian angle of $Y_{\zetabold}$}
$\phantom{ab}$
\nopagebreak

We define maps $\Ytilde[n]$ on $\widetilde{S}[n]$,
which can be interpreted---in a sense made precise in
\ref{L:Ytilde}---as limits of $Y_\zetabold$ as $\tb\to0$:


\addtocounter{equation}{1}
\begin{definition}
\label{D:Ytilde}
Arguing as in the proof of \ref{L:Y} we define a map $\Ytilde[0]:\widetilde{S}[0]\to\Sph^2$
by requiring it satisfies the symmetries in \ref{E:symmetries},
it is the inclusion map on $M_0$,
and it equals $X_0$ on
$\Lambda[1]$.
For $n'$ as in \ref{E:nnprime} we define 
$\Ytilde[n']:\widetilde{S}[n']\to\Sph^2$
by
$
\Ytilde[n']:=X_0\circ \TTT_{-2n'\ptb}.
$
\end{definition}

\addtocounter{equation}{1}
\begin{lemma}
\label{L:Ytilde}
For $n$ and $n'$ as in \ref{E:nnprime} we have the following:
$\Ytilde[n]$ is a diffeomorphism from $\widetilde{S}[n]$
onto a domain of $\Sph^2$ and
for $\mb$ large enough in terms of a given $x>0$ we have
$$
\begin{aligned}
\|Y_\zetabold-\Ytilde[0]:C^{2,\beta}(S_x[0],(\Ytilde[0])^*g_{\Sph^2})\|
\le&
C(\tb+|\zeta'_1|+|\zeta'_2|)
\le
C\cunder c'_2\tb,
\\
\|Y_\zetabold-\ptilde_{\tau,\alpha}^{n'} \circ \Ytilde[n']:C^{2,\beta}(S_x[n'],(\Ytilde[n'])^*g_{\Sph^2})\|
\le&
C(\tau+|\alpha|)
\le
C\tb,
\end{aligned}
$$
where the maps are considered as $\mathbb{C}^3$-valued.
\end{lemma}

\begin{proof}
The restriction of $Y_{\zetabold}$ to $S_x[0]$ depends smoothly on the parameters $\tau$,
$\zeta_1$, and $\zeta_2$,
and gives $\Ytilde$ when all the parameters vanish,
therefore the result follows when $n=0$.
By 
\ref{E:tphat:a:under-sym}
we have
$$
\Yzeta=
\ptilde_{\tau,\alpha}^{n'} \circ \Xunderta \circ \TTT_{-2n'\ptb}
\qquad
\text{ on $\widetilde{S}[n']$}.
$$
Since $\Xunderta$ depends smoothly on $\tau$ and $\alpha$
and it gives $X_0$ when $\tau=\alpha=0$, the result follows.
\end{proof}

Before we discuss the metrics we will be using 
we define a cut-off function we will need as follows:

\addtocounter{equation}{1}
\begin{definition}
\label{D:psihat}
$\psihat:M\to[0,1]$ is defined by 
$\psihat=\psi[a+1,a]\circ t$ on each $\{j\}\times\Sph^1\times[a,a+1]$,
$\psihat=\psi[2m\ptb-a- 1,2m\ptb-a]\circ t$ on each
$\{j\} \times \Sph^1 \times [2m\ptb-a- 1,2m\ptb-a]$,
$\psihat$ vanishes on each $M_j\setminus M_0$,
and $\psihat\equiv1$ on $M_0\setminus\bigcup_{j=1}^gM_j$.
\end{definition}

Note that $\psihat$ and $1-\psihat$ form a partition of unity
subordinate to
$\{M_0,\bigcup_{j=1}^gM_j\}$.
Also that $\psihat$ is invariant under the action of $\group$
(in the usual sense or see \ref{D:psi:symmetries}).
In the next definition we define metrics $g$ and $\chi$ on $M$
and $\Sph^1\times\R$ we will be using in the estimates of the paper.

\addtocounter{equation}{1}
\begin{definition}
\label{D:chi}
We define on $\Sph^1\times\R$ $\underline{\rho}_\tau:=\rho_\tau\circ\underline{t}$,
$g=\Xunderta^*\,g_{\Sph^5}$, and $\chi:=\underline{\rho}_\tau^{-2} g$.
On $M$ we define 
$g:=g_\zetabold:=Y_\zetabold^*\,g_{\Sph^5}$,
and $\chi:=\chi_\zetabold:=\rho_\zetabold^{-2}g$,
where
$$
\rho=\rho_\zetabold:=\psihat + \,(1-\psihat)\,\underline{\rho}_\tau
\quad
\text{ where }  \quad \underline{\rho}_\tau:=\rho_\tau\circ\underline{t}\text{ on }M_j.
$$
\end{definition}

Clearly then $g$ and $\chi$ on $M$ depend smoothly on $\zetabold$,
while $g$ and $\chi$ on $\Sph^1\times\R$ depend smoothly on $\tau$ and $\alpha$.
$\group$ acts by isometries on $(M,g)$ and $(M,\chi)$.
On each neck $\Lambda[n']\subset M$
and each $\Sph^1\times[2+2n\ptb,2(n+1)\ptb-2]\subset\Sph^1\times\R$ we have
\addtocounter{theorem}{1}
\begin{equation}
\label{E:chi:necks}
\chi=ds^2+\frac{\pt^2}{\ptb^2}\,dt^2=ds^2+d\underline{t}^2.
\end{equation}
Finally on the spherical regions $S_x[n]$ all the metrics in consideration
($g$ and $\chi$ for the various $\zetabold$)
are uniformly equivalent as can be seen by applying
\ref{L:Ytilde} and \ref{P:conformal}.
By uniformly equivalent here we mean that if $h_1$ and $h_2$ are any two such metrics,
then
\addtocounter{theorem}{1}
\begin{equation}
\label{E:equivalent}
\|h_1:C^k(S_x[n],h_2)\|
\le
C(k,x),
\end{equation}
where $C(k,x)$ depends only on $k$ and $x$.

%

We study now the Lagrangian angle $\theta$ induced on $M$ by $Y_\zetabold$.
We first decompose the Lagrangian angle $\theta$ on $M$ into three parts by writing 
\addtocounter{theorem}{1}
\begin{equation}
\label{E:decomposition}
\theta=\thetaglue + \thetadisloc + \thetatwist,
\end{equation}
where $\thetadisloc$ is supported on $\group (\Sph^1 \times [a+1,a+2])$,
$\thetaglue$ is supported on $\group (\Sph^1 \times [a+2,a+3])$ and 
$\thetatwist$ is supported on the complement of $S[0]$ and the necks.
$\thetaglue$ is created by the transition from the round equatorial $2$-sphere to a shifted $X_{\tp}$,
and depends only on $\tp$.
$\thetadisloc$ is due to the mismatch induced by the sliding and twisting controlled by the 
parameters $\zetabold$.
It vanishes when $\zetabold$ vanishes.
Finally $\thetatwist$ is due to the twisting introduced periodically on the spherical regions of $X_{\tp}$
to convert it to $X_{\tp,\alpha}$.
It vanishes when $\alpha$---equivalently $\zeta_2$---vanishes.
It is easy to estimate these as follows, where $g$ could have been used instead of $\chi$
as well:

\begin{lemma}
\addtocounter{equation}{1}
\label{L:theta}
(i).
$\Vert \thetaglue:C^k(M,\chi)\Vert  \le C(k)\, \tp\le
C(k)\,\,\tb$.
\newline
(ii).
$\Vert \thetadisloc:C^k(M,\chi)\Vert  \le C(k)\, |(\zeta'_1,\zeta'_2)|
\le C(k)\,\,\cunder\,c'_2\,\tb$.
\newline
(iii).
$\Vert \thetatwist:C^k(M,\chi)\Vert  \le C(k)\, |\alpha|
\le
C(k)\,\cunder \,c'_2\,\tb^2 \,|\log\tb|$.
\end{lemma}

\begin{proof}
These estimates follow from the smooth dependence on the parameters, 
\ref{E:sigma:range}, and \ref{E:tau:alp:interval}.
\end{proof}

Note that \ref{L:theta}.(iii) is equivalent to
\addtocounter{theorem}{1}
\begin{equation}
\label{E:theta}
\Vert \theta:C^k(\Sph^1\times\R,\chi)\Vert
\le C(k)\, |\alpha|
\le
C(k)\,\cunder \,\tb^2 \,|\log\tb|,
\end{equation}
where $\theta$ is the Lagrangian angle on $\Sph^1\times\R$ induced by $\Xunderta$.

\section{The Linearized Equation}
\label{S:linear}
\nopagebreak

\subsection*{Introduction}
$\phantom{ab}$
\nopagebreak

In this section we study the linearized equation and the linearized operator
for functions which should correct the Legendrian immersions we are working
with to special Legendrian immersions.
Our purpose is to state and prove Propositions \ref{P:solution} and \ref{P:solution:Xta}.
We start by defining the linear operators $\Lcal$ and $\Lchi$ by (recall \ref{D:chi}):
\addtocounter{theorem}{1}
\begin{equation}
\label{E:L}
 \Lcal := \Delta_g + 6=\rho^{-2}\Lchi,
\qquad\text{where}\qquad
\Lchi:=
\Delta_\chi+6\rho^2
=\rho^2\Lcal.
\end{equation}
By proposition \ref{P:Leg:pert}
$\Lcal$ is the linearized operator we are interested in.
$\Lchi$ is a conformally modified version of $\Lcal$
which we will find useful later.
The inhomogeneous linear equation we will be studying is
\addtocounter{theorem}{1}
\begin{equation}
\label{E:Lequation}
 \Lcal u= \rho^{-2} E,
\qquad\text{or, equivalently,}\qquad
\Lchi u= E.
\end{equation}

\subsection*{Symmetries}
$\phantom{ab}$

The symmetries involved play an important role.
The next definition is motivated by proposition \ref{P:f:symmetries}.
Note the odd symmetry with respect to $\tbar$.

\addtocounter{equation}{1}
\begin{definition}
\label{D:f:symmetries}
We will call a function $f$ defined on $M$ appropriately symmetric if and only if
the following hold:
\begin{equation*}
f = -f \circ \tbar, \qquad f = f \circ \sbar, \qquad f = f \circ \SSS_\pi, \qquad f = f \circ \RRR.
\end{equation*} 
We use the subscript ``$\sym$'' to denote subspaces of appropriately symmetric functions,
for example
$C^{k,\beta}_\sym(M)=\{f\in C^{k,\beta}(M):
f = -f \circ \tbar, \ f = f \circ \sbar, \ f = f \circ \SSS_\pi, \ f = f \circ \RRR\}$.

More generally, a function defined on $\Omega\subset M$
will be called appropriately symmetric if it can be
extended to a function satisfying the above symmetries on $\group\Omega$.
\end{definition}

We will be using various cut-off functions which
need to have even symmetry with respect to $\tbar$,
so that they preserve the odd symmetry of a function in multiplication.
An example of such a cut-off function is $\psihat$ defined in \ref{D:psihat}.
To avoid confusion we have the following:

\addtocounter{equation}{1}
\begin{definition}
\label{D:psi:symmetries}
We will call a function $\psi$ defined on a $\group$-inva\-ri\-ant
domain of $M$ $\group$-inva\-ri\-ant if and only if the following hold:
\begin{equation*}
\psi = \psi \circ \tbar,
\qquad
\psi = \psi \circ \sbar,
\qquad
\psi = \psi \circ \SSS_\pi,
\qquad
\psi = \psi \circ \RRR.
\end{equation*}
\end{definition}

Clearly by this definition
if $f$ is appropriately symmetric in the sense of \ref{D:f:symmetries}
and $\psi$ $\group$-invariant,
then $\psi f$ is appropriately symmetric.

The functions on $M$ controlling the perturbations of $Y_\zetabold$
will be appropriately symmetric as in
\ref{D:f:symmetries}.
The following proposition allows us to make use of these symmetries:

\addtocounter{equation}{1}
\begin{prop}
\label{P:f:symmetries}
Let $\,Y:M \ra \Sph^5$ be a smooth Legendrian immersion
and assume $Y$ satisfies the symmetries \ref{E:symmetries}.
\newline
(i).
The Lagrangian angle $\theta$ of $Y$ is appropriately symmetric in the sense
of \ref{D:f:symmetries}.
\newline
(ii).
If $f\in C^{2,\beta}_\sym(M)$
is appropriately small so that
the Legendrian perturbation $Y_f$ of $Y$ is well defined
as in Appendix \ref{A:perturbation},
then $Y_f$ satisfies the same symmetries \ref{E:symmetries} as $Y$,
and therefore $\Lcal f,\Lchi f, \theta_f\in C^{0,\beta}_\sym(M)$.
\end{prop}

\begin{proof}
(i).
We extend $Y:M\to\Sph^5$ to $Y_{cone}:M\times\R^+\to\C^3$
by
$Y_{cone}(p,R)=R\,Y(p)$.
We also extend $\tbar$ to a symmetry of $M\times\R^+$
by $\tbar(p,R)=(\tbar(p),R)$ and similarly for the other symmetries.
We have then by \ref{E:symmetries} that
$$
Y_{cone}^*\Omega=\tbar^*\circ  Y_{cone}^* \circ \tbartilde^* \Omega.
$$
Using \ref{E:lagn:phase} and \ref{P:s5:commute}.(iv)
we conclude
$$
e^{i\theta}\vol=\tbar^*(-e^{-i\theta}\vol)=e^{-i\theta\circ\tbar}\vol,
$$
where the last equality follows because $\tbar$ is an orientation-reserving isometry.
$\theta=-\theta\circ\tbar$ then follows.
The other symmetries follow by a similar argument using
$\sbartilde^* \Omega  = - \Omega$ (see \ref{P:s5:commute}.(iv))
and the invariance of $\Omega$ by $\stilde_\pi, \rtilde\in\text{SU(3)}$.
\newline
(ii).
The vector field $V$ as defined on a tubular neighborhood of $Y$ in Appendix \ref{A:perturbation}
satisfies
$$
\tbartilde\circ V=\tbartilde(-J\nabla f)=J\tbartilde(\nabla f)=
V\circ\tbartilde,
$$
since $\tbartilde\circ J=-J\circ\tbartilde$ by \ref{P:s5:commute}.ii
and $f=-f\circ\tbar$.
$\tbartilde\circ Y_f=Y_f\circ\tbar$ follows.
The other symmetries follow by a similar argument using
the commutation properties with $J$ in \ref{P:s5:commute}.i.
\end{proof}

We have also to perturb $\Xta$ to a special Legendrian immersion satisfying the symmetries.
For this we need the following analogues of \ref{D:f:symmetries}
and \ref{P:f:symmetries}:

\addtocounter{equation}{1}
\begin{definition}
\label{D:f:symmetries:Xta}
We will call a function $f$ defined on $\Sph^1\times\R$
appropriately symmetric
if and only if
the following hold:
\begin{equation*}
f = -f \circ \tbar, \quad f = f \circ \TTT_{2\ptb}, \quad f = f \circ \sbar, \quad f=f\circ \SSS_{\pi}.
\end{equation*} 
We use the subscript ``$\sym$'' to denote subspaces of appropriately symmetric functions
as for example in
$C^{k,\beta}_\sym(\Sph^1\times\R)$.
\end{definition}

\addtocounter{equation}{1}
\begin{prop}
\label{P:f:symmetries:Xta}
(i).
The Lagrangian angle $\theta$ of any immersion
satisfying the same symmetries
$\Xunderta:\Sph^1 \times \R\ra\Sph^5$
satisfies in
\ref{L:sliding:twisting:X},
is appropriately symmetric in the sense of
\ref{D:f:symmetries:Xta}.
\newline
(ii).
If $f\in C^{2,\beta}_\sym(\Sph^1\times \R)$
is appropriately small so that
the Legendrian perturbation $\Xunder_{\tau,\alpha,f}$ of $\Xunderta$ is well defined
as in Appendix \ref{A:perturbation},
then $\Xunder_{\tau,\alpha,f}$ satisfies the same symmetries in
\ref{L:sliding:twisting:X} as $\Xunderta$,
and therefore $\Lcal f,\Lchi f, \theta_f\in C^{0,\beta}_\sym(\Sph^1\times \R)$.
\end{prop}

\begin{proof}
(i).
The proof for $\theta=-\theta\circ \tbar$ is similar to the one in 
\ref{P:f:symmetries}(i).
The other symmetries follow 
also from \ref{L:a-sym} 
and the invariance of $\Omega$ by $\ptilde_{\tau,\alpha}, \stilde_\pi, \rtilde\in\text{SU(3)}$.
\newline
(ii).
This is also similar to the proof of 
\ref{P:f:symmetries}(ii).
\end{proof}

\subsection*{The linearized equation on the necks}
$\phantom{ab}$
\nopagebreak

In this subsection we consider the linearized equation on a neck $\Lambda_{x,y}[n'']$
defined as in \ref{E:neck}, where we assume that $x,y\in[0,5]$.
For simplicity in this subsection we will denote the neck under consideration by $\Lambda$, 
and its boundary circles $C_x^+[n''-1]$ and $C^-_y[n'']$ by $C^+$ and $C^-$ respectively.
We next define $\xunder,\xover,\xboth:\Lambda\ra\R$
to measure the $t$-coordinate distance from $C^+$, $C^-$, and $\partial\Lambda=C^+\cup C^-$
respectively (recall \ref{E:neck}):
\addtocounter{theorem}{1}
\begin{equation}
\label{E:xdistance}
(2n'-2)\ptb+b+x+\xunder
=t,
\quad
2n'\ptb-b-y-\xover
=t,
\quad
\xboth:=\min(\xunder,\xover).
\end{equation}
We define $\ell$ to be the $t$-coordinate length of the cylinder,
so that
\addtocounter{theorem}{1}
\begin{equation}
\label{E:ell}
\ell=2\ptb-2b-x-y,
\qquad
\xunder + \xover\equiv\ell.
\end{equation}
Recall (see \ref{D:chi} and \ref{E:chi:necks}) that the metrics on $\Lambda$ are given by
\addtocounter{theorem}{1}
\begin{equation}
\label{E:Lambdametrics}
g=\rho_\tau^2\circ\underline{t}\,\chi = \rho^2\,\chi,
\qquad
\chi=ds^2+d\underline{t}^2=ds^2+ \frac{\pt^2}{\ptb^2}\,dt^2=\rho^{-2}\,g.
\end{equation}

The functions we will consider on $\Lambda$ are required to satisfy the symmetries
\addtocounter{theorem}{1}
\begin{equation}
\label{E:f:Lsymmetries}
f = f \circ \sbar, \qquad f = f \circ \SSS_\pi.
\end{equation} 
This is consistent with \ref{D:f:symmetries} because $\sbar$ and $\SSS_\pi$
preserve $\Lambda$ while $\tbar$ and $\RRR$ do not.
We use the subscript ``$\SSS$'' to denote subspaces of functions
which satisfy \ref{E:f:Lsymmetries},
as for example 
$C^{k,\beta}_\SSS(\Lambda)=\{f\in C^{k,\beta}:f = f \circ \sbar, \ f = f \circ \SSS_\pi\}$.
Note that exponential decays with respect to the $t$ and $\underline{t}$ coordinates differ
(recall \ref{E:tau:alp:interval})
by factors of order $e^\tau$ and so we can ignore the difference from now on.
In particular by choosing $b$ large enough in terms of $\varepsilon_1$
and using \ref{E:exp:decay:y}
we have that
\addtocounter{theorem}{1}
\begin{equation}
\label{E:conformaldecay}
\Vert \rho^2: C^3(\Lambda,\chi,e^{-2\xboth})\Vert \le \varepsilon_1,
\end{equation}
where $\varepsilon_1$ is a small positive constant to be determined later.
Using \ref{E:pt} and \ref{E:tau:alp:interval}
we have
\addtocounter{theorem}{1}
\begin{equation}
\label{E:ell:bound}
|\ell+\log\tb|\le C(\varepsilon_1).
\end{equation}
We are considering therefore the linear operator $\Lchi$ acting on functions on a long
cylinder by \ref{E:ell:bound} and 
\ref{E:Lambdametrics},
and where by \ref{E:conformaldecay} $\Lchi$ is a small perturbation of $\Delta_\chi$.
This leads to the following:

\addtocounter{equation}{1}
\begin{prop}
\label{P:Lambda:low}
The lowest eigenvalue of the Dirichlet problem for $\Lchi$ on $\Lambda$ is $>C\ell^{-2}$.
\end{prop}

\begin{proof}
The proof is similar to the arguments leading to Proposition 2.28 in \cite{kapouleas:wente}.
We omit the details.
\end{proof}

\addtocounter{equation}{1}
\begin{corollary}
\label{C:L:unique}
(i).
The Dirichlet problem for $\Lchi$ on $\Lambda$ for given $C^{2,\beta}$ Dirichlet data has a unique solution.
\newline
(ii). 
For $E\in C^{0,\beta}(\Lambda)$ there is a unique $\varphi\in C^{2,\beta}(\Lambda)$
such that $\Lchi \varphi=E$ on $\Lambda$ and $\varphi=0$ on $\partial \Lambda$.
Moreover
$
\|\varphi :C^{2,\beta}(\Lambda,\chi)\|
\le
C(\beta)\,\ell^2\,
\|E:C^{0,\beta}(\Lambda,\chi)\|.
$
\end{corollary}

\begin{proof}
(i) follows trivially and (ii) by using standard linear theory.
\end{proof}


In the next Proposition and its Corollary we study the Dirichlet problem
when we are allowed to modify the lower harmonics on the boundary data
so that we can get decay estimates appropriate for our purposes:

\addtocounter{equation}{1}
\begin{prop}
\label{P:Lambda:solution}
Assuming $\varepsilon_1$ small enough in terms of given
$\beta\in(0,1)$ and $\gamma\in(2,3)$,
there is a linear map $\Rcal_\Lambda:C^{0,\beta}_\SSS(\Lambda)\to C^{2,\beta}_\SSS(\Lambda)$
such that the following hold for
$E\in C^{0,\beta}_\SSS(\Lambda)$ and $V:=\Rcal_\Lambda\, E$:
\newline
(i).
$\Lchi V= E$ on $\Lambda$.
\newline
(ii).
$\left.V\right|_{C^+}\in \left<1,\cos 2s\right>_\R$
and $V$ vanishes on $C^-$.
\newline
(iii).
$\|V:C^{2,\beta}_\SSS(\Lambda,\chi,e^{-\gamma\xunder})\|
\le
C(\beta,\gamma)\,
\|E:C^{0,\beta}_\SSS(\Lambda,\chi,e^{-\gamma\xunder})\|$.
\newline
(iv).
$\Rcal_\Lambda$ depends continuously on $\tau$.

The theorem still holds if the roles of $C^+$ and $C^-$ are exchanged in (ii)
and $\xunder$ is replaced by $\xover$ in (iii).
Another possibility is to allow both
$\left.V\right|_{C^+}$ and
$\left.V\right|_{C^-}$
to be in $\left<1,\cos 2s\right>_\R$
in (ii)
while $\xunder$ is replaced by $\xboth$ in (iii).
%
\end{prop}

\begin{proof}
By 
\ref{E:f:Lsymmetries}
the only harmonics allowed on the meridians of order up to two
are the constants and $\cos 2s$.
The proposition then follows by standard theory if $\Lchi$ is replaced by $\Delta_\chi$.
We denote the corresponding linear map and solution in the $\Delta_\chi$ case
by $\widetilde{\Rcal}_\Lambda$ and $\widetilde{V}$ respectively.
Using then \ref{E:conformaldecay} we have
$$
\|\Lchi\,\widetilde{V}:C^{0,\beta}(\Lambda,\chi,e^{-\gamma\xunder})\|
\le C(\beta,\gamma)\,\varepsilon_1
\|E:C^{0,\beta}(\Lambda,\chi,e^{-\gamma\xunder})\|,
$$
and the proposition then follows by an iteration where we treat $\Lchi$ and $\Rcal_\Lambda$
as small perturbations of $\Delta_\chi$ and $\widetilde{\Rcal}_\Lambda$ by assuming
$\varepsilon_1$ small enough.
\end{proof}

%
%

\addtocounter{equation}{1}
\begin{corollary}
\label{C:Lambda:solution}
Assuming $\varepsilon_1$ small enough in terms of given
$\varepsilon_2>0$, $\beta\in(0,1)$, and $\gamma\in(2,3)$,
there is a linear map 
$$
\Rcal_\partial:\{u\in C_\SSS^{2,\beta}(C^+):u \text{\ is $L^2(C^+,ds^2)$-orthogonal to }\left<1,\cos 2s\right>_\R\}
\to C_\SSS^{2,\beta}(\Lambda)
$$
such that the following hold for $u$ in the domain of $\Rcal_\partial$ and $V:=\Rcal_\partial u$:
\newline
(i).
$\Lchi V= 0$ on $\Lambda$.
\newline
(ii).
$\left.V\right|_{C^+}-u\in \left<1,\cos 2s\right>_\R$
and $V$ vanishes on $C^-$.
\newline
(iii).
$\|\left.V\right|_{C^+}-u:C_\SSS^{2,\beta}(C^+,ds^2)\|
\le
\varepsilon_2
\,\|u:C_\SSS^{2,\beta}(C^+,ds^2)\|$.
\newline
(iv).
$\|V:C_\SSS^{2,\beta}(\Lambda,\chi,e^{-\gamma\xunder})\|
\le
C(\beta,\gamma)
\,\|u:C_\SSS^{2,\beta}(C^+,ds^2)\|$.
\newline
(iv).
$\Rcal_\partial$ depends continuously on $\tau$.

The Proposition still holds if the roles of $C^+$ and $C^-$ are exchanged 
and $\xunder$ is replaced by $\xover$.
\end{corollary}

\begin{proof}
By standard theory there is a linear map
$$
\widetilde{\Rcal}_\partial:
\{u\in C_\SSS^{2,\beta}(C^+):u \text{\ is $L^2(C^+,ds^2)$-orthogonal to }\left<1,\cos 2s\right>_\R\}
\to C_\SSS^{2,\beta}(\Lambda)
$$
such that for $u$ in the domain and $\widetilde{V}=\widetilde{\Rcal}_\partial u$
the following hold:
\newline
(a).
$\Delta_\chi \widetilde{V}= 0$ on $\Lambda$.
\newline
(b).
$\left.\widetilde{V}\right|_{C^+}=u$ on $C^+$
and $\widetilde{V}$ vanishes on $C^-$.
\newline
(c).
$\|\widetilde{V}:C_\SSS^{2,\beta}(\Lambda,\chi,e^{-\gamma\xunder})\|
\le
C(\beta,\gamma)
\, \|u:C_\SSS^{2,\beta}(C^+,ds^2)\|$.

The corollary then follows by defining
$$
\Rcal_\partial u:=
\widetilde{\Rcal}_\partial u
-
\Rcal_\Lambda\,
\Lchi
\widetilde{\Rcal}_\partial u,
$$
applying the Proposition,
and
using
\ref{E:conformaldecay}.
\end{proof}

In the next proposition \ref{P:Lambda:V}
we estimate the solutions to the 
Dirichlet problem for $\Lchi$ on $\Lambda$,
with Dirichlet data the two lowest harmonics allowed by the symmetries
\ref{E:f:Lsymmetries}.
The estimates in \ref{P:Lambda:V} compare these solutions with
the corresponding solutions for $\Delta_\chi$
which are given explicitly in \ref{E:Lambda:Vtilde}.
To facilitate reference to these solutions we have the following:

\addtocounter{equation}{1}
\begin{definition}
\label{D:Lambda:V}
For $i=1,2$ we denote by $V_i[\Lambda,a_1,a_2]$
and $\widetilde{V}_i[\Lambda,a_1,a_2]$
the solutions to the Dirichlet problems
on $\Lambda$ given by
$$
\Lchi V_i[\Lambda,a_1,a_2]=0,
\qquad
\Delta_\chi \widetilde{V}_i[\Lambda,a_1,a_2]=0,
$$
with boundary data
\begin{align*}
V_1[\Lambda,a_1,a_2]=
\widetilde{V}_1[\Lambda,a_1,a_2]=&a_1
& &
\text{on}
\quad
C^+,
\\
V_1[\Lambda,a_1,a_2]=
\widetilde{V}_1[\Lambda,a_1,a_2]=&a_2
& &
\text{on}
\quad
C^-,
\\
V_2[\Lambda,a_1,a_2]=
\widetilde{V}_2[\Lambda,a_1,a_2]=&a_1\cos 2s
& &
\text{on}
\quad
C^+,
\\
V_2[\Lambda,a_1,a_2]=
\widetilde{V}_2[\Lambda,a_1,a_2]=&a_2\cos 2s
& &
\text{on}
\quad
C^-.
\end{align*}
\end{definition}

Note that $V_i[\Lambda,a_1,a_2]$ and $\widetilde{V}_i[\Lambda,a_1,a_2]$
are linear on each of $a_1$ and $a_2$
and the roles of $C^+$ and $C^-$ can be exchanged.
Therefore it is enough to understand
$V_i[\Lambda,1,0]$ and $\widetilde{V}_i[\Lambda,1,0]$.
It is straightforward to check
\addtocounter{theorem}{1}
\begin{equation}
\label{E:Lambda:Vtilde}
\widetilde{V}_1[\Lambda,1,0]=\frac{\xover}{\ell},
\qquad
\widetilde{V}_2[\Lambda,1,0]=\frac{\,\sinh(2\pt\xover/\ptb)\,}{\sinh(2\pt\ell/\ptb)}\,\cos2s.
\end{equation}


\addtocounter{equation}{1}
\begin{prop}
\label{P:Lambda:V}
$V_1[\Lambda,a_1,a_2]$ is constant on the meridians,
and 
$V_2[\Lambda,a_1,a_2]$ 
on each me\-ri\-di\-an is a multiple of $\cos 2s$.
Moreover
by assuming $\varepsilon_1$ small enough
in terms of a given $\varepsilon_3>0$,
there are constants $A_1$, $A_1^-$, and $A_2$ such that the following hold:
\newline
(i). $|A_1-1|\le\varepsilon_3$ and $|A_1^-|\le \varepsilon_3/\ell$.
\newline
(ii).
$\|V_1[\Lambda,1,0]-\widetilde{V}_1[\Lambda,A_1,A_1^-]:C^{2,\beta}(\Lambda,\chi,e^{(1-\gamma)\xunder}
+\ell^{-1}e^{(1-\gamma)\xover})\|\le\varepsilon_3$.
\newline
(iii). $|A_2-1|\le\varepsilon_3$.
\newline
(iv).
$\|V_2[\Lambda,1,0]-\widetilde{V}_2[\Lambda,A_2,0]:C^{2,\beta}(\Lambda,\chi,e^{-\gamma\xunder}
+e^{-3\ell/2})\|\le\varepsilon_3$.
\end{prop}

\begin{proof}
The rotational invariance of $\Lchi$ on $\Lambda$ implies the first part of the Proposition.
Using \ref{E:Lambda:Vtilde} and \ref{E:conformaldecay}
we conclude
$\|\Lchi\widetilde{V}_1[\Lambda,1,0]:C^{2,\beta}(\Lambda,\chi,e^{(1-\gamma)\xunder}
+\ell^{-1}e^{(1-\gamma)\xover})\|\le\varepsilon_1\,C(\beta,\gamma)$
\newline
and
$\|\Lchi\widetilde{V}_1[\Lambda,0,1]:C^{2,\beta}(\Lambda,\chi,e^{(1-\gamma)\xover}
+\ell^{-1}e^{(1-\gamma)\xunder})\|\le\varepsilon_1\,C(\beta,\gamma)$.
Using then an appropriately  modified version of \ref{P:Lambda:solution}
to account for the weaker decay available,
we find $\widehat{V}_1[\Lambda,1,0]$, $\widehat{V}_1[\Lambda,0,1]\in C^{2,\beta}(\Lambda)$
such that the following hold:
\newline
(a).
$\Lchi\widehat{V}_1[\Lambda,1,0]=-\Lchi\widetilde{V}_1[\Lambda,1,0]$
and
$\Lchi\widehat{V}_1[\Lambda,0,1]=-\Lchi\widetilde{V}_1[\Lambda,0,1]$.
\newline
(b).
$\widehat{V}_1[\Lambda,1,0]$ and $\widehat{V}_1[\Lambda,0,1]$
are constant on $C^+$ and on $C^-$ by the rotational invariance.
\newline
(c).
$\|\widehat{V}_1[\Lambda,1,0]:C^{2,\beta}(\Lambda,\chi,e^{(1-\gamma)\xunder}
+\ell^{-1}e^{(1-\gamma)\xover})\|\le\varepsilon_1\,C(\beta,\gamma)$
\newline
and
$\|\widehat{V}_1[\Lambda,0,1]:C^{2,\beta}(\Lambda,\chi,e^{(1-\gamma)\xover}
+\ell^{-1}e^{(1-\gamma)\xunder})\|\le\varepsilon_1\,C(\beta,\gamma)$.

$A_1$ and $A_1^-$ are determined then uniquely by requiring that
\addtocounter{theorem}{1}
\begin{equation}
\label{E:estimate:V}
V_1[\Lambda,1,0]=\widetilde{V}_1[\Lambda,A_1,A_1^-]+
A_1\,\widehat{V}_1[\Lambda,1,0]
+
A_1^-\,\widehat{V}_1[\Lambda,0,1]
\end{equation}
holds
on $\partial\Lambda=C^+\cup C^-$.
Since $\Lchi$ kills both sides we conclude by \ref{P:Lambda:low}
that the equality holds on $\Lambda$.
(i) and (ii) follow then from the available estimates.

Using again \ref{E:Lambda:Vtilde} and \ref{E:conformaldecay}
we conclude
$\|\Lchi\widetilde{V}_2[\Lambda,1,0]:
C^{2,\beta}(\Lambda,\chi,e^{-\gamma\xunder} +e^{-2\ell})\|
\le\varepsilon_1\,C(\beta,\gamma)$.
Arguing then in a similar way as for (i) and (ii) we prove (iii) and (iv)
completing the proof.
\end{proof}

\addtocounter{equation}{1}
\begin{corollary}
\label{C:Lambda:u}
If $u\in C_\SSS^{2,\beta}(\Lambda)$, $\mathcal{L}u=0$ on $\Lambda$,
and $u=0$ on $C^-$,
then
$$
\|u:C_\SSS^{2,\beta}(\Lambda,\chi,(\xover+1)/\ell)\|
\le
C(\beta)\,
\|u:C_\SSS^{2,\beta}(C^+,\chi)\|.
$$
\end{corollary}

\begin{proof}
We decompose on $C^+$ $u=u_1+u_2+u_3$
where $u_1$ is constant,
$u_2$ is a multiple of $\cos 2s$,
and $u_3$ is in the domain of $\Rcal_\partial$
defined in 
\ref{C:Lambda:solution}.
Applying then \ref{P:Lambda:V} and 
\ref{C:Lambda:solution}
we obtain estimates which imply the result.
\end{proof}

\subsection*{The approximate kernel}
$\phantom{ab}$
\nopagebreak

We proceed now to discuss the approximate kernel of $\mathcal{L}$ on the various extended 
standard regions, \textit{cf.} \cite[Prop. 2.22]{kapouleas:wente}.
By approximate kernel we mean the span of eigenfunctions whose eigenvalues are close to $0$.
We understand the approximate kernel in the next proposition by comparing it to
\addtocounter{theorem}{1}
\begin{equation}
\label{E:f1f2}
\widehat{f}_1=
\ftbold\circ \Yzeta \,/\,
\|\ftbold:L^2(\Sph^2)\|,
\qquad
\widehat{f}_2=
\fqbold\circ \Yzeta \,/\,
\|\fqbold:L^2(\Sph^2)\|,
\end{equation}
where
as in section \ref{S:init:surf}
$\Sph^2=\left<e_1,e_2,e_3\right>_\R\cap\Sph^5$,
and $\ftbold,\fqbold:\Sph^5\to\R$
are defined by 
\addtocounter{theorem}{1}
\begin{equation}
\label{E:f:tq}
\ftbold = -\tfrac12 |z_1|^2 + \tfrac14 |z_2|^2 + \tfrac14 |z_3|^2, \qquad \fqbold = -\tfrac12 |z_2|^2 + \tfrac12 |z_3|^2,
\end{equation}
where $z_1$, $z_2$, $z_3$ are the standard coordinates in $\C^3$.
Note that $\mathbf{t} = -J\nabla \ftbold$ and $\mathbf{q} = -J\nabla \fqbold$
are the two Killing vector fields which generate the $1$-parameter groups 
$\{\ttilde_x\}_{x\in\R}$ and $\{\qtilde_x\}_{x\in \R}$ respectively.
Note that by the definitions of $\ttilde_x$ and $\qtilde_x$ in \ref{D:sym:tilde},
$\ftbold$ and $\fqbold$ are preserved by the action of $\ttilde_x$ and $\qtilde_x$,
and so by 
\ref{E:tphat:sym}
their restrictions to $\widetilde{S}_x[n']$ do not depend on $n'$.
Their restrictions on the necks depend only on $\tau$.
Using $\Yzeta$ instead of $\Ytilde[n']$ in \ref{E:f1f2} is helpful in the proof of \ref{L:phi:inp}.

\addtocounter{equation}{1}
\begin{prop}
\label{P:app:ker}
Assuming $b$ large enough in absolute terms,
and $\tb$ small enough in terms of a given $\varepsilon > 0$,
the following hold:
\newline
(i).
$\mathcal{L}$ acting
on appropriately symmetric (recall \ref{D:f:symmetries}) functions on $\widetilde{S}[0]$
with vanishing Dirichlet conditions,
has no eigenvalues
in $[-1,1]$.
\newline
(ii).
$\mathcal{L}$ acting
on appropriately symmetric
(recall \ref{D:f:symmetries})
functions
on $\group\widetilde{S}[n']$
($n'$ as in \ref{E:nnprime}),
and with vanishing Dirichlet conditions,
has exactly $2$ eigenvalues
in $[-\varepsilon,\varepsilon]$,
and no other eigenvalues in $[-1,1]$.
We will refer to the 
$2$-dimensional vector space spanned by the corresponding eigenfunctions
as the approximate kernel of $\mathcal{L}$ on $\widetilde{S}[n']$.
Moreover the approximate kernel has an orthonormal basis $\{f_{1,n'},f_{2,n'}\}$
where for $i=1,2$,
$f_{i,n'}\in C^{2,\beta}(\widetilde{S}[n],\chi)$
depends continuously on $\zetabold$ and satisfies
$$
\Vert  f_{i,n'} - \widehat{f}_i \,: C^{2,\beta}(S_5[n']) \Vert  < \varepsilon.
$$
\end{prop}

\begin{proof}
The proof is based on the results of \cite[Appendix B]{kapouleas:annals}
which are based on basic facts about eigenvalues and eigenfunctions \cite{chavel}.
Before using those results we remark the following:
First, the first inequality in \cite[B.1.6]{kapouleas:annals}
should read
$$
\|F_i f\|_\infty \le2\|f\|_\infty
$$
instead.
Second, the spaces of functions can be constrained to satisfy appropriate symmetries,
as indeed was the case in some of the constructions in \cite{kapouleas:annals},
and will be the case here.
Third, the only use of the Sobolev inequality
\cite[B.1.5]{kapouleas:annals}
is to establish supremum bounds for the eigenfunctions.
These in our case can be alternatively established by using
the uniformity of geometry of $S[n]$
to obtain interior estimates on $S_1[n]$,
and then using
\ref{C:Lambda:u}
to obtain estimates on the necks.
We proceed to discuss the proof in the two cases under consideration.

(i).
Consider $\Stilde[0]\setminus\group\widetilde{C}[1]$.
It consists of $2g+1$ connected components.
We denote by $\Stilde^+[0]$ the closure of the connected component which contains $M_0$,
and by $\Lambda^-[1]$ the connected component which contains $C^-[1]$.
Note then that $\Lambda^-[1]\subset\Lambda[1]$,
and that $\Stilde[0]$ is the disjoint union of
$\Stilde^+[0]$, $\RRR^j\Lambda^-[1]$, and $\RRR^j \tbar \Lambda^-[1]$, where $j\in \{1, \ldots ,g\}$.

Consider now two abstract copies of $\Sph^2$,
$\{1\}\times\Sph^2$ and $\{2\}\times\Sph^2$,
and their disjoint union $\{1,2\}\times\Sph^2$.
We define the action of $\group$ on $\{1,2\}\times\Sph^2$ by requiring that it acts on each copy
as on $\Sph^2$.
We also define a metric $g$ on $\{1,2\}\times\Sph^2$ by taking it to restrict to the usual
$g_{\Sph^2}$ on each copy.
$\group$ acts then on $\{1,2\}\times\Sph^2$ by isometries.

Recall \ref{D:Ytilde}.
Let $D$ be the smallest geodesic disc in $\Sph^2$ which contains $\Ytilde[1](\Lambda^-[1])$.
For illustration purposes note that $\partial D=\Ytilde[1](C^-[1])$.
We define
$$
\Shat[0]:=(\{1\}\times\Sph^2)\,
\bigcup\,(\{2\}\times\bigcup_{j=1}^g \RRR^j D)
\bigcup\,(\{2\}\times\bigcup_{j=1}^g \RRR^j \tbar D).
$$
Note that all the unions are disjoint and $\Shat[0]$ is invariant under the action of $\group$.
Moreover $\Shat[0]$ is the disjoint union of $\{1\}\times\Sph^2$ and $\group(\{2\}\times D)$,
the latter consisting of $2g$ isometric connected components.
$\partial\Shat[0]$ consists of $2g$ circles.

Recall that $\group$ acts by isometries on
$\left(\widetilde{S}[0],g=\Yzeta^*g_{\Sph^5}\right)$,
and we are considering the eigenfunctions of $\mathcal{L}=\Delta_g+6$
acting on functions satisfying the symmetries in \ref{D:f:symmetries}
with Dirichlet boundary data.
To prove (i) we compare with eigenfunctions of $\Delta_g+6$
acting on functions on
$(\Shat[0],g)$ 
satisfying also  \ref{D:f:symmetries},
and with Dirichlet boundary data.

In order to apply the results of  \cite[Appendix B]{kapouleas:annals}
we define maps $F_1$ and $F_2$ as follows:
We start by defining a map $\Yhat[0]:\Stilde[0]\to\Shat[0]$ by requesting that
it is equivariant under the action of $\group$ on $\Stilde[0]$ and $\Shat[0]$,
and that for $p\in \Stilde^+[0]$ we have $\Yhat[0](p):=(1,\Ytilde[0](p))$
and for $p\in\Lambda^-[1]$ we have $\Yhat[0](p):=(2,\Ytilde[1](p))$.
Recall \ref{L:Ytilde} and note that $\Yhat[0]$ is injective.
We fix a $d>0$ which we assume large enough in terms of $\varepsilon$.
Given 
$f\in C^\infty_0(\Stilde[0])$,
we define $F_1(f)\in C^\infty_0(\Shat[0])$
by requiring that it vanishes on the complement of the image of $\Yhat[0]$
and that on the image of $\Yhat[0]$ it satisfies
$$
F_1(f)\circ\Yhat[0]=\psicheck\,f,
$$
where $\psicheck$ is a cut-off function defined on $\widetilde{S}[0]$
as follows:
$\psicheck$ satisfies \ref{D:psi:symmetries},
equals $1$ on $S[0]$,
and $\psi[2d,d]\circ\xboth$ on $\Lambda[1]$.
Note that this provides a logarithmic cut-off as in  \cite[IV.2.4]{kapouleas:annals}.

Conversely, given 
$f\in C^\infty_0(\Shat[0])$,
we define $F_2(f)\in C^\infty_0(\Stilde[0])$,
by requiring 
$$
F_2(f)=\psicheck\,f\circ\Ytilde[0].
$$
By assuming $d$ large enough in terms of $\varepsilon$,
is is then straightforward to check the hypotheses needed so that the results of
\cite[Appendix B]{kapouleas:annals} apply.

It remains to check that $\Delta_g+6$ 
acting on functions on
$(\Shat[0],g)$ 
satisfying \ref{D:f:symmetries},
has no eigenvalues in $[-2,2]$.
By assuming $b$ large enough we can ensure that the radius of the disc $D$
is small enough so that the smallest eigenvalue of $\Delta_g$ on $(D,g_{\Sph^2})$
with vanishing Dirichlet data is larger than $8$.
We need therefore to be concerned only with the spectrum of the Laplacian $\Delta$ on $\Sph^2$
which is well known \cite{chavel}.
In particular the only eigenvalue in $[4,8]$ is $6$ with corresponding eigenfunctions
given by
$$
f=
\mu_1 x_1^2+
\mu_2 x_2^2+
\mu_3 x_3^2+
\mu'_1 x_2 x_3 +
\mu'_2 x_3 x_1+
\mu'_3 x_1 x_2,
$$
where $\mu_1+\mu_2+\mu_3=0$,
$\mu_1,\mu_2,\mu_3,\mu'_1,\mu'_2,\mu'_3\in\R$,
and $x_1,x_2,x_3$ are the standard coordinates on $\Sph^2$.
Using \ref{E:sym:bartilde} we see that the action of $\sbar$ and $\SSS_\pi$ gives
any of
$(x_1,x_2,x_3)\mapsto (x_1,\pm x_2,\pm x_3)$,
and hence for $f$ to satisfy \ref{D:f:symmetries}
we need $\mu'_1=\mu'_2=\mu'_3=0$.
Similarly by \ref{E:sym:bartilde} the action of $\tbar$ gives
$(x_1,x_2,x_3)\mapsto (-x_1,x_2,x_3)$.
$f=-f\circ\tbar$ implies then that
$\mu_1=\mu_2=\mu_3=0$.
Therefore there are no nontrivial eigenfunctions respecting the symmetries and the result
follows.

(ii).
In this case the stabilizer of $\widetilde{S}[n]$
is generated by $\sbar$ and $\SSS_\pi$,
therefore we are interested in the eigenfunctions of $\mathcal{L}$
acting on functions on $\widetilde{S}[n]$
satisfying 
\addtocounter{theorem}{1}
\begin{equation}
\label{E:L:symmetries}
f=f\circ\sbar=f\circ\SSS_\pi
\end{equation}
and Dirichlet boundary data.
To prove (ii) we compare as for (i) with
eigenfunctions of $\Delta_g+6$
acting on functions on
$(\Shat[n'],g)$ 
satisfying also
\ref{E:L:symmetries}
and with Dirichlet boundary data.
Here $\Shat[n']\subset\{1,2\}\times\Sph^2$
is defined in analogy with $\Shat[0]$ as follows:
$$
\Shat[n']:=\{1\}\times\Sph^2
\cup(\{2\}\times D^-)
\cup(\{2\}\times D^+),
$$
where $D^-$ and $D^+$ are the smallest geodesic discs in $\Sph^2$ containing
$\Ytilde[n'+1](\Lambda^-[n'+1])$ and 
$\Ytilde[n'](\Lambda^+[n'])$,
where
$\Lambda^-[n'+1]$ and $\Lambda^+[n']$ are the components of
$\Stilde[n']\setminus(\widetilde{C}[n'+1]\cup\widetilde{C}[n'])$
containing $C^-[n'+1]$ and $C^+[n'-1]$ respectively.
We also define $\Stilde^+[n']$ to be the closure of the third
remaining component.
Therefore $\Stilde[n']$ is the disjoint union of $\Lambda^+[n']$, $\Stilde^+[n']$,
and $\Lambda^-[n'+1]$.
Note that $\Shat[n']$ is invariant under the action of $\sbar$ and $\SSS_\pi$
on $\{1,2\}\times\Sph^2$,
and therefore it makes sense to consider functions
on $\Shat[n']$ satisfying \ref{E:L:symmetries}.

The definition of the maps $F_1$ and $F_2$ is analogous to the one in (i)
with some obvious modifications:
First,
instead of $\Yhat[0]$ we use $\Yhat[n']:\Stilde[n']\to\Shat[n']$
which is defined as follows:
For $p\in\Lambda^+[n']$ we have $\Yhat[n'](p)=(2,\Ytilde[n'-1](p))$,
for $p\in\Stilde^+[n']$ we have $\Yhat[n'](p)=(1,\Ytilde[n'](p))$,
and for $p\in\Lambda^-[n'+1]$ we have $\Yhat[n'](p)=(2,\Ytilde[n'+1](p))$.
Second,
$\psicheck$ is now defined on $\widetilde{S}[n']$
by requiring that $\psicheck\equiv1$ on $S[n']$
and 
$\psicheck=\psi[2d,d]\circ\xboth$ on $\Lambda[n'+1]\cup\Lambda[n']$.

By arguing as in (i),
and since the odd symmetry with respect to $\tbar$ is lacking,
it is clear that $\Delta+6$ acting on functions on $\Sph^2$
satisfying \ref{E:L:symmetries},
has a two-dimensional kernel with its eigenfunctions given by
$\mu_1 x_1^2+ \mu_2 x_2^2+ \mu_3 x_3^2$,
where $\mu_1,\mu_2,\mu_3\in\R$
and $\mu_1+\mu_2+\mu_3=0$.
An orthonormal basis of the kernel is given then by 
$$
\widehat{f}_{1,\Sph^2}=
\ftbold\,/\,
\|\ftbold:L^2(\Sph^2)\|,
\qquad
\widehat{f}_{2,\Sph^2}=
\fqbold\,/\,
\|\fqbold:L^2(\Sph^2)\|,
$$
and therefore on $S_d[n']$ we have
\begin{equation*}
F_2(\widehat{f}_{1,\Sph^2})=
\ftbold\circ \Ytilde[n'] \,/\,
\|\ftbold:L^2(\Sph^2)\|,
\qquad
F_2(\widehat{f}_{2,\Sph^2})=
\fqbold\circ \Ytilde[n'] \,/\,
\|\fqbold:L^2(\Sph^2)\|.
\end{equation*}
%
Since $\ftbold$ and $\fqbold$ are invariant under the action of
$\ptilde_{\tau,\alpha}$,
we conclude by assuming $\tb$ small enough and using \ref{L:Ytilde},
that 
$$
\Vert  F_2(f_{i,\Sph^2}) - \widehat{f}_i \,: C^{2,\beta}(S_6[n']) \Vert  < \varepsilon/2,
\qquad
(i=1,2).
$$

The proof of the estimate is then completed by applying the results of 
\cite[Appendix B]{kapouleas:annals} and upgrading the $L^2$ to $C^{2,\beta}$ estimates 
by using the uniformity of the geometry of $S_6[n']$ and standard linear theory interior estimates.
To ensure the continuous dependence of $f_{i,n'}$ on $\zetabold$ we choose $f_{1,n'}$ to be the closest
element of the approximate kernel to $\widehat{f}_i$ in the $L^2(\widetilde{S}[n],g)$ metric.
The continuous dependence follows then by standard arguments.
\end{proof}

We also need to understand the approximate kernel for $\Xta$:

\addtocounter{equation}{1}
\begin{prop}
\label{P:app:ker:Xta}
$\mathcal{L}$ acting
on appropriately symmetric functions
on
$\Sph^1\times\R$
(satisfying \ref{D:f:symmetries:Xta}),
has no eigenvalues in $[-1,1]$.
\end{prop}

\begin{proof}
The proof is similar to the one for \ref{P:app:ker}.(i), only simpler.
The comparison is with the spectrum of $\Delta+6$
acting on functions on $(\Sph^2,g_{\Sph^2})$ satisfying the symmetries in \ref{D:f:symmetries:Xta}
except for the second one,
that is functions $f$ satisfying
$$
f = -f \circ \tbar, \qquad f = f \circ \sbar, \qquad f=f\circ \SSS_{\pi}.
$$
We have already seen in the proof of \ref{P:app:ker}.(i) that the corresponding
kernel is trivial.
\end{proof}

\subsection*{The extended substitute kernel}
$\phantom{ab}$
\nopagebreak

In order to solve the linearized equation globally on the initial surface
following the general methodology of 
\cite{kapouleas:survey,kapouleas:wente},
we will need to modify the inhomogeneous term by elements of the extended substitute kernel
\cite[\S 19]{kapouleas:survey}, which we now proceed to define.

\addtocounter{equation}{1}
\begin{definition}
\label{D:win}
Following \cite{kapouleas:survey,kapouleas:wente}
we define the extended substitute kernel $\skernel$ by
$\skernel:=
\bigoplus_{n=0}^{2\mb-1}\skernel[n]$,
where $\skernel[n]:= \left<w_{1,n}, w_{2,n}\right>_\R$,
where $w_{i,n}$ $(i=1,2)$ are smooth,
appropriately symmetric in the sense of \ref{D:f:symmetries},
functions on $M$,
determined as follows:

$w_{i,0}$ is supported on $\group\Lambda[1]$,
and on $\Lambda[1]$ 
in the notation of 
\ref{D:Lambda:V},
we have
\addtocounter{theorem}{1}
\begin{equation}
\label{E:wi0}
w_{i,0}:=
\widetilde{c}_i\,
\mathcal{L}\left(\psi[0,1]\circ\xunder \,\, V_i[\Lambda[0],1,0]\right),
\end{equation}
where $\widetilde{c}_i$ are constants defined by
$$
\widetilde{c}_1\, V_1[\Lambda[0],1,0] =1
\quad \text{and}\quad
\widetilde{c}_2\, V_2[\Lambda[0],1,0] =\cos 2s
\quad \text{on} \quad
C^+_1[0].
$$

For $n'$ as in \ref{E:nnprime}
$w_{i,n'}$ is supported on $\group S[n']$
and satisfies
$$
w_{i,n'} := c_i \,\psi[2n'\ptb-1,2n'\ptb]\circ t\,\psi[2n'\ptb+1,2n'\ptb]\circ t\,\widehat{f}_i
\qquad
\text{ on }S[n'],
$$
where the coefficients $c_i$
depend on $\zetabold$ and are determined by the requirements
\addtocounter{theorem}{1}
\begin{equation}
\label{E:win:dij}
\int_{S[n']}w_{i,n'}\,\widehat{f}_j\,dg= \delta_{ij} \qquad i,j \in \{1,2\}.
\end{equation}
\end{definition}

To motivate the definition of the extended substitute kernel $\skernel$ just given,
we remark that the linearized equation is solved later modulo the extended substitute kernel (see \ref{P:solution}.(i)).
As we will see in the proof of \ref{L:wS:solution} there are two reasons for this:
First, since the approximate kernel on $\Stilde[n]$ is nontrivial when $n\ne0$ (see \ref{P:app:ker}),
we have to solve modulo $\skernel[n]$ to ensure that the inhomogeneous term is orthogonal to
the approximate kernel (see \ref{L:L2w}.(iii)).
This purpose could be achieved as well by using the substitute kernel 
$\bigoplus_{n=1}^{2\mb-1}\skernel[n]$ instead of $\skernel$.
The second reason is that to ensure appropriately fast exponential decay along the necks
we need to be able to prescribe the low harmonics of the solution on $C^+_1[n]$.
For this we allow ourselves the freedom of modifying the semi-local solutions on $\Stilde[n]$
by elements of $\skernel_v[n]$ defined in \ref{D:vskernel}.
This means that the inhomogeneous term gets modified by $\Lcal v$ for some $v\in\skernel_v[n]$.
Fortunately we can arrange for $\Lcal v\in \skernel[n]$ (see \ref{L:vskernel}.(i)),
so when $n\ne0$ we do not need to extend the substitute kernel $\skernel[n]$ to anything new.
In the case $n=0$ however, where the approximate kernel is trivial,
we need the nontrivial $\skernel[0]$ defined above.
This forces us to define and use the extended substitute kernel $\skernel$,
instead of the substitute kernel,
as is often the case
\cite{kapouleas:survey,kapouleas:wente,kapouleas:finite}.

We record now the following for future reference:

\addtocounter{equation}{1}
\begin{lemma}
\label{L:L2w}
(i).
$w_{i,n}$ is supported on $\group S_1[n]$.
\newline
(ii).
$\Vert w_{i,n}:C^{2,\beta}(M,\chi)\Vert  \le C$.
\newline
(iii).
For $n'$ as in \ref{E:nnprime}
and $E\in C^{0,\beta}(\widetilde{S}[n'],\chi)$,
there is a unique
$\widetilde{w}\in\skernel[n']$ such that $E+\widetilde{w}$
is $L^2(\widetilde{S}[n'],g)$-orthogonal to the approximate kernel on $\widetilde{S}[n']$.
Moreover if $E$ is supported on $S_1[n']$, then
$$
\|\widetilde{w}:C^{2,\beta}(M,\chi)\|
\le
C(b)\,\|E: C^{0,\beta}(\widetilde{S}[n'],\chi)\|.
$$
\end{lemma}

\begin{proof}
(i) and (ii) follow from the definitions and (when $n=0$) \ref{P:Lambda:V}.
(iii) follows then from \ref{P:conformal}.(iv), \ref{P:app:ker}, and \ref{E:win:dij}.
\end{proof}


It is useful to define normalization constants  $c'_1,c''_1,c'_2,c''_2\in\R$
by requesting that on each $C^+_1[n']$ the following hold:

\addtocounter{theorem}{1}
\begin{equation}
\label{E:c1c2}
c'_1 \ftbold \circ \Yzeta = c''_1 \widehat{f}_1=1,
\qquad
c'_2 \fqbold \circ \Yzeta = c''_2 \widehat{f}_2=\cos 2s.
\end{equation}

To arrange for the required decay along the necks we need to be able to
prescribe the low harmonics on $C^+_1[n]$ at the solution level.
For this purpose we need the following:

\begin{definition}
\addtocounter{equation}{1}
\label{D:vskernel}
We define
$\skernel_v[n]:=\left<v_{1,n},v_{2,n}\right>_\R$,
where $v_{i,n}$ $(i=1,2)$ are smooth,
appropriately symmetric in the sense of \ref{D:f:symmetries},
functions on $\group \widetilde{S}[n]$,
determined as follows:

$v_{i,0}$ is supported on $\group\Lambda[1]$,
and on $\Lambda[1]$ 
we have
$
v_{i,0}:=
\widetilde{c}_i\,
\psi[0,1]\circ\xunder \,\, V_i[\Lambda[0],1,0]
$.

For $n'$ as in \ref{E:nnprime}
$v_{i,n'}$ is defined by
\addtocounter{theorem}{1}
\begin{equation}
\label{E:vin}
 v_{i,n'} := c''_i\,(f_{i,n'} + u_{i,n'}),
\end{equation}
where $c''_i$ is as in \ref{E:c1c2},
$ f_{i,n'}$ is as in proposition \ref{P:app:ker},
and $u_{i,n'}$ is the solution of the Dirichlet problem on
$\group\widetilde{S}[n']$,
$$ \mathcal{L} u_{i,n'} = -\mathcal{L} f_{i,n'} + \widetilde{w}_{i,n'},$$
with vanishing boundary data,
and $\widetilde{w}_{i,n'}\in \skernel[n']$
is determined by \ref{L:L2w}(iii) so that the
the inhomogeneous term 
is $L^2(\widetilde{S}[n'],g)$-orthogonal to the approximate kernel on $\widetilde{S}[n']$.
\end{definition}


For future reference we record the following:

\addtocounter{equation}{1}
\begin{lemma}
\label{L:vskernel}
By assuming $b$ large enough in terms of $\varepsilon_3$ we can ensure that 
the $v_{i,n}$'s defined above are smooth on $\widetilde{S}[n]$ and satisfy the following:
\newline
(i).
$\mathcal{L} v_{i,n} \in \skernel[n]$
and therefore
$\Lcal v_{i,n}$
and
$\Lchi v_{i,n}$
are supported on $\group S_1[n]$.
\newline
(ii).
$v_{1,n}=0$ on $\partial \widetilde{S}[n]$.
\newline
(iii).
$\Vert v_{i,n}: C^{2,\beta}(\widetilde{S}[n],\chi)\Vert  \le C(b).$
\newline
(iv).
$\Vert v_{1,n}-1:C^{2,\beta}(C^+_1[n],\chi)\Vert  < \varepsilon_3$
and
$\Vert v_{2,n}-\cos 2s:C^{2,\beta}(C^+_1[n],\chi)\Vert  < \varepsilon_3$.
Moreover
on $C^+_1[0]$ we have $v_{1,0}=1$ and $v_{2,0}=\cos 2s$.
\end{lemma}

\begin{proof}
The case $n=0$ 
follows from the definitions, \ref{L:L2w}, and \ref{P:Lambda:V}.
In the case $n\ne0$ (i) and (ii) follow from the definitions.
Using now \ref{P:app:ker} we have
$
\|\Lcal f_{i,n'}:L^2(\Stilde[n'],g)\|
\le \varepsilon,
$
which implies using \ref{L:L2w} that
$
\|u_{i,n'}:L^2(\Stilde[n'],g)\|
\le C \varepsilon.
$
The $L^2$ estimate gives then interior $C^{2,\beta}$ estimates which together
with \ref{C:Lambda:u} imply (iii) and together with \ref{P:app:ker} imply (iv).
\end{proof}

\subsection*{Solving the linearized equation semi-locally}
$\phantom{ab}$
\nopagebreak

Lemmas \ref{L:L2w} and \ref{L:vskernel} provide us with what we need
to solve with appropriate estimates the linearized problem on extended
spherical regions.
As usual we impose the appropriate symmetry.
To facilitate reference we have the following:

\addtocounter{equation}{1}
\begin{definition}
\label{D:stilde:symmetries}
We use the subscript ``$\ \SSS$'' to denote the subspace of a space of functions
on $\Stilde[n]$ defined by requesting the following:
\newline
If $n=0$ the functions should be appropriately symmetric in the sense of \ref{D:f:symmetries}.
\newline
If $n\ne0$ then the functions should satisfy \ref{E:f:Lsymmetries}.
\end{definition}

The above definition is consistent with the fact that $\Stilde[0]$ is invariant under the action
of $\group$,
while the stabilizer of $\Stilde[n]$ for $n\ne0$ under the action of $\group$ is generated by
$\sbar$ and $\SSS_\pi$.
Because of our earlier lemmas on the linearized equation on the necks it is enough---see the proof
of Proposition \ref{P:solution}---to assume in the next lemma
that the inhomogeneous term is supported on $S_{1}[n]$.
The range of $n$ is as usual as in \ref{E:nnprime}.

\addtocounter{equation}{1}
\begin{lemma}
\label{L:wS:solution}
There is a linear map
$$
\Rcal_{\Stilde[n]}:
\{E\in C^{0,\beta}_\SSS(\Stilde[n]): E \text{ is supported on }S_{1}[n]\}
\to
C^{2,\beta}_\SSS(\Stilde[n])\times\skernel[n],
$$
such that the following hold for
$E$ in the domain of $\Rcal_{\Stilde[n]}$ above
and $(\varphi,w)=\Rcal_{\Stilde[n]}(E)$:
\newline
(i).
$\Lchi \varphi = E + \rho^2{w}$---equivalently
$\Lcal \varphi = \rho^{-2}E + {w}$---on $\widetilde{S}[n]$.
\newline
(ii).
$\varphi$ vanishes on $\partial \widetilde{S}[n]$.
\newline
(iii).
$|\mu_1|+|\mu_2|\le
C(b,\beta)\,\Vert E: C_\SSS^{0,\beta}({S}_1[n],\chi)\Vert$
where $w=\mu_1\,w_{1,n}+\mu_2\,w_{2,n}$.
\newline
(iv).
$\Vert \varphi: C_\SSS^{2,\beta}(\widetilde{S}[n],\chi) \Vert 
\le C(b,\beta)\, \Vert E: C_\SSS^{0,\beta}({S}_1[n],\chi)\Vert .$
\newline
(v).
$\Vert \varphi: C_\SSS^{2,\beta}(\Lambda[n+1],\chi,e^{-\gamma\xunder})\Vert 
\le C(b,\beta,\gamma)\, \Vert E: C_\SSS^{0,\beta}({S}_1[n],\chi)\Vert .$
\newline
(vi).
If $n\ne0$ then 
$\Vert \varphi: C_\SSS^{2,\beta}(\Lambda[n],\chi,(\xover+1)/\ell) \Vert 
\le C(b,\beta)\,\Vert E: C_\SSS^{0,\beta}({S}_1[n],\chi)\Vert .$
\newline
(vii).
$\Rcal_{\Stilde[n]}$ depends continuously on $\zetabold$.
\end{lemma}

\begin{proof}
Let $\wt$ be as in \ref{L:L2w}.iii so that $\rho^{-2}E+\wt$
is orthogonal to the approximate kernel if $n\ne0$,
and $\wt=0$ if $n=0$.
We can solve then uniquely to find 
$\varphi'\in C_\SSS^{2,\beta}(\widetilde{S}[n])$
vanishing on $\partial\widetilde{S}[n]$,
and such that on $\Stilde[n]$
$$
\Lcal \varphi' =\rho^{-2} E + \wt,
\text{ or, equivalently, }
\Lchi \varphi' = E + \rho^{2}\wt.
$$
We define then $\varphi=\varphi'+v$,
where $v\in\skernel_v[n]$ is chosen as follows:
Consider the decomposition $\varphi=\varphi_{low}+\varphi_{\perp}$ on $C^+_1[n]$,
where $\varphi_{low}\in\left< 1, \cos 2s\right>_\R$
and $\varphi_{\perp}$ is an element of the domain of $\Rcal_\partial$
defined in \ref{C:Lambda:solution} with $\Lambda=\Lambda_{1,0}[n+1]$.
$v$ is uniquely determined by requesting
$\Rcal_\partial\varphi_{\perp}-\varphi_{\perp}=\varphi_{low}$
on $C^+_1$ which implies $\Rcal_\partial\varphi_\perp=\varphi$ on $\Lambda$,
providing this way the necessary estimates.
\end{proof}

\subsection*{Solving the linearized equation globally}
$\phantom{ab}$
\nopagebreak

In order to solve the linearized equation \ref{E:Lequation}
globally on $M$ and provide estimates for the solutions,
we paste together the semi-local solutions provided by
\ref{P:Lambda:solution} and \ref{L:wS:solution}
in a construction we proceed to present:
We start by defining various cut-off functions we will need.
In the next definition $n$, $n'$, and $n''$, assume the ranges specified in \ref{E:nnprime}:

\addtocounter{equation}1
\begin{definition}
\label{D:partition}
We define uniquely $\psi_{S[n]}$,
$\psi_{\widetilde{S}[n]}$,
and $\psi_{\Lambda[n'']}$,
smooth functions on $M$,
by requesting the following:
\newline
(i).
They are $\group$-invariant in the sense of
\ref{D:psi:symmetries}.
\newline
(ii).
$\psi_{S[n]}$ is supported on $\group S_{1}[n]$,
$\psi_{\widetilde{S}[n]}$ is supported on $\group \widetilde{S}[n]$,
and $\psi_{\Lambda[n'']}$ is supported on $\group \Lambda[n'']$.
\newline
(iii).
$\psi_{S[n]}\equiv\psi_{\widetilde{S}[n]}\equiv1$ on $S[n]$.
$\psi_{S[n]}=\psi[1,0]\circ\xunder$ and
$\psi_{\widetilde{S}[n]}=\psi[0,1]\circ\xover$ on $\Lambda[n+1]$.
$\psi_{S[n']}=\psi[1,0]\circ\xover$ and
$\psi_{\widetilde{S}[n']}=\psi[0,1]\circ\xunder$ on $\Lambda[n']$.
$\psi_{\Lambda[n'']}=\psi[0,1]\circ\xboth$
on $\Lambda[n'']$.
\end{definition}

Note that the functions
$\psi_{S[n]}$ and $\psi_{\Lambda[n'']}$
form a partition of unity on $M$.
The cut-off function $\psi_{\widetilde{S}[n]}$
is identically one on
$\widetilde{S}[n]$
except close to the boundary
$\partial \widetilde{S}[n]$
where it transits smoothly to $0$.

In order to state the main proposition of this section we need to define appropriate norms.
Note that the fast rate of exponential decay imposed along the necks will be useful later,
especially when we prescribe $\skernel[n]$ for $n\ne0$ in \ref{L:phi:inp}.

\addtocounter{equation}{1}
\begin{definition}
\label{D:globalnorm}
For $k\in\mathbb{N}$,
$\beta\in(0,1)$, and $\gamma\in(2,3)$,
we define a norm $\|\, .\,\|_{k,\beta,\gamma}$
on $C^{k,\beta}_\sym(M)$
by taking 
$\|f\|_{k,\beta,\gamma}$ to be the maximum of
the following semi-norms
with $n$ as in \ref{E:nnprime}:
\newline
(i).
$\tb^{-\gamma n}\|f:C^{k,\beta}(S_1[n],\chi)\|$.
\newline
(ii).
$\tb^{-\gamma n}\|f:C^{k,\beta}(\Lambda[n+1],\chi,e^{-\gamma\xunder})\|$,
except that $\xunder$ is replaced with $\xboth$ when $n+1=2\mb$.
\end{definition}

We also need norms with appropriate decay on $\skernel$ and $\R^{4\mb}$  which we define by
\addtocounter{theorem}{1}
\begin{equation}
\label{E:wnorm}
\|\boldsymbol{\mu}\|_\gamma=\left\|\sum_{i,n}\mu_{i,n}w_{i,n}\right\|_\gamma:=\max_{i,n}\,\tb^{-\gamma n}|\mu_{i,n}|
\end{equation}
for $\boldsymbol{\mu}=\{\mu_{1,n},\mu_{2,n}\}_{n=0}^{2\mb-1}\in\R^{4\mb}$.

To combine naively functions defined on $\widetilde{S}[n]$ to a global function on $M$
we will be using the following:

\addtocounter{equation}{1}
\begin{lemma}
\label{L:comb}
Given $u_n\in C^{k,\beta}_\SSS(\Stilde[n])$
for each $n$ as in \ref{E:nnprime},
with $u_n$ vanishing in a neighborhood of $\partial\Stilde[n]$,
there is a unique
$U=\Ubold(\{u_n\})\in C^{k,\beta}_\sym(M)$
such that the following hold for $n$ and $n'$ as in \ref{E:nnprime}:
\newline
(i). $U=u_n$ on $S[n]$.
\newline
(ii). $U=u_{n'-1}+u_{n'+1}$
on $\Lambda[n']$.
\newline
(iii). $U=u_{2\mb-1}-u_{2\mb-1}\circ\tbar$
on $\Lambda[2\mb]$.

Moreover in the case that $u_{2\mb-1}$ satisfies $u_{2\mb-1}=-u_{2\mb-1}\circ\tbar$ on $\Lambda[2\mb]$,
there is a unique $\underline{U}=\Uboldunder(\{u_n\})\in C^{k,\beta}_\sym(M)$
such that (i) and (ii) hold
and
\newline
(iii\'{}).
$\underline{U}=u_{2\mb-1}$
on $\Lambda[2\mb]$.
\end{lemma}

\begin{proof}
$U$ and $\underline{U}$ are clearly uniquely defined by (i), (ii), and (iii)
on $\bigcup_n\Stilde[n]$.
Using the odd symmetry with respect to $\tbar$ they are then uniquely extended to
$M_0\cup M_1$, and then using $\RRR$ uniquely extended to $M$.
It is straightforward then to check that they satisfy the required symmetries and smoothness.
\end{proof}

\addtocounter{equation}{1}
\begin{prop}
\label{P:solution}
There is a linear map $\Rcal_M:C^{0,\beta}_{sym}(M)\to C^{2,\beta}_{sym}(M)\times\skernel$
such that for $E \in C^{0,\beta}_{sym}(M)$
and $(\varphi,w)=\Rcal_M E$
the following hold:
\newline
(i).
$\Lchi \varphi = E + \rho^2{w}$,
or equivalently
$\Lcal \varphi = \rho^{-2}E + {w}$,
on $M$.
\newline
(ii).
$\|\varphi\|_{2,\beta,\gamma}\le C(b,\beta,\gamma) \, \|E\|_{0,\beta,\gamma}$.
\newline
(iii).
$\|{w}\|_{\gamma}\le C(b,\beta,\gamma) \, \|E\|_{0,\beta,\gamma}$.
\newline
(iv).
$\Rcal_M$ depends continuously on $\zetabold$.
\end{prop}

\begin{proof}
We first prove the proposition under the assumption that $E$ 
is supported on $\bigcup_n\group S_1[n]$:
First we apply \ref{L:wS:solution} to obtain $\varphi_n\in C^{2,\beta}(\widetilde{S}[n])$
and $w_n\in\skernel[n]$
such that
$$
(\varphi_n,w_n)=\Rcal_{\Stilde[n]}(\left. E\right|_{S_1[n]}).
$$
We define then
$$
\begin{gathered}
\Wcal E:=\sum_n w_n\in\skernel,
\qquad
\Rcal E:=\Ubold(\{\psi_{\widetilde{S}[n]}\varphi_n\})\in C^{2,\beta}_{sym}(M),
\\
\Ecal E:=\Ubold(\{[\psi_{\widetilde{S}[n]},\Lchi]\varphi_n\})\in C^{0,\beta}_{sym}(M),
\end{gathered}
$$
where $[\,,\,]$ denotes the commutator, that is
$[\psi_{\widetilde{S}[n]},\Lchi]\varphi_n=\psi_{\widetilde{S}[n]}\,\Lchi \varphi_n-
\Lchi(\psi_{\widetilde{S}[n]}\varphi_n)$.

It is straightforward to check then that
\addtocounter{theorem}{1}
\begin{equation}
\label{E:REWcal:0}
\Lchi\Rcal E+\Ecal E=E+\rho^2\Wcal E\qquad\text{ on } M,
\end{equation}
that $\Ecal E$ is supported on 
$\bigcup_n\group S_1[n]$,
and that we have 
\newline
$\|\Rcal E\|_{2,\beta,\gamma}\le C(b,\beta,\gamma) \, \|E\|_{0,\beta,\gamma}$,
\newline
$\|\Ecal E\|_{0,\beta,\gamma}\le C(b,\beta,\gamma) \, \tau^{\gamma'-\gamma}\,\|E\|_{0,\beta,\gamma}$,
\newline
$\|{w}\|_{\gamma}\le C(b,\beta,\gamma) \, \|E\|_{0,\beta,\gamma}$,
\newline
where for the second estimate we used
\ref{L:wS:solution}
with $\gamma$ replaced by a fixed $\gamma'\in(\gamma,3)$.
By assuming $\tb$ small enough then,
the second estimate gives
\newline
$\|\Ecal E\|_{0,\beta,\gamma}\le \frac12\,\|E\|_{0,\beta,\gamma}$.

By induction then we have for $r=0,1,\ldots$
\addtocounter{theorem}{1}
\begin{equation}
\label{E:REWcal:k}
\Lchi\Rcal \Ecal^r E+\Ecal^{r+1} E=\Ecal^r E+\rho^2\Wcal \Ecal^r E\qquad\text{ on } M,
\end{equation}
and using the estimates above we can define
\addtocounter{theorem}{1}
\begin{equation}
\label{E:RcalWcal}
\varphi:=\sum_{r=0}^\infty\Rcal\Ecal^r E,
\qquad
w:=\sum_{r=0}^\infty\Wcal\Ecal^r E,
\end{equation}
and complete the proof in this case where $E$ is assumed supported on 
$\bigcup_n\group S_1[n]$.

To prove the proposition in general we first apply \ref{P:Lambda:solution}
with $\Lambda=\Lambda[n'']$
to find $V_{n''}=\Rcal_{\Lambda[n'']} \left. E \right|_{\Lambda[n'']}$
where we require
$\left.V_{n''}\right|_{C^+}\in \left<1,\cos 2s\right>_\R$,
and,
when $n''\ne2\mb$,
$V_{n''}=0$ on $C^-$,
while when $n''=2\mb$
we require by appealing to the uniqueness in \ref{C:L:unique}
$V_{2\mb}=-V_{2\mb}\circ\tbar$
on $\Lambda[2\mb]$.
We define then
$$
\widetilde{E}:=\Ubold(\{\psi_{S[n]}\,E\})
+
\Uboldunder(\{[\psi_{\Lambda[n+1]},\Lchi]V_{n+1}\})\in C^{0,\beta}_{sym}(M),
$$
which is clearly supported on
$\bigcup_n\group S_1[n]$.
We apply then the special case of the proposition we have already proven
with $\widetilde{E}$ instead of $E$,
to obtain an element of $C^{2,\beta}_{sym}(M)$
which we will call $\widetilde{\varphi}$ instead of $\varphi$,
and a $w\in\skernel$.
In particular we have $\Lchi\widetilde{\varphi}=\widetilde{E}+w$ on $M$.

We define then $\varphi\in C^{2,\beta}_{sym}(M)$ by
$$
\varphi=\widetilde{\varphi}+\Uboldunder(\{\psi_{\Lambda[n+1]}\, V_{n+1}\}).
$$
It is easy then to check that the $\varphi$ and $w$ we defined
satisfy the required properties and the proof is complete.
\end{proof}

\subsection*{The linearized equation for $\Xunderta$}
$\phantom{ab}$
\nopagebreak


We need to understand the linear equation on $(\Sph^1\times\R,\chi)$ in the fashion of
\ref{P:solution}.
Fortunately,
because of the lack of approximate kernel (see \ref{P:app:ker:Xta}),
the linear theory in this case is much simpler than in the case considered in
\ref{P:solution}.
The following is enough for our purposes:

\addtocounter{equation}{1}
\begin{prop}
\label{P:solution:Xta}
There is a linear map $\Rcal_\cyl:C^{0,\beta}_{sym}(\cyl)\to C^{2,\beta}_{sym}(\cyl)$
such that for $E \in C^{0,\beta}_{sym}(\cyl)$
and $\varphi=\Rcal_\cyl E$
the following hold:
\newline
(i).
$\Lchi \varphi = E $,
or equivalently $\Lcal\varphi=\rho^{-2}E$,
on $\cyl$,
where $\Lchi$ is induced by $\Xunderta$ (see \ref{D:chi} and \ref{E:L}).
\newline
(ii).
$\|\varphi:C^{2,\beta}_\sym(\cyl,\chi)\|
\le C(b,\beta) \,
(\,
\|E:C^{0,\beta}_\sym(\cyl,\chi)\|
+ |\log\tb|^2\,
\|E:C^{0,\beta}_\sym(\widetilde{\Lambda},\chi)\|
\,)$,
where $\widetilde{\Lambda}=\Sph^1\times[b,2\ptb-b]\subset\cyl$.
(Notice that $(\widetilde{\Lambda},\chi)$ is isometric to each $(\Lambda[n''],\chi)$
when $\tau$ and $\alpha$ are determined from $\zetabold$ through \ref{E:sigma}.)
\newline
(iii).
$\Rcal_\cyl$ depends continuously on $\tau$ and $\alpha$.
\end{prop}

\begin{proof}
In analogy with the proof of \ref{P:solution} we first prove the proposition under the assumption
that $E$ vanishes on $\Sph^1\times[b+1,2\ptb-b-1]\subset\cyl$.
By \ref{P:conformal}.iv we have then that on the support of $E$, $\rho^{-2}<C(b)$,
and hence
$$
\|\rho^{-2}E:L^2(\cyl/\TTT_{2\ptb},g)\|\le
C(b)\,\|E:C^{0,\beta}_\sym(\cyl,\chi)\|,
$$
where $\cyl/\TTT_{2\ptb}$ is the quotient of $\cyl$ by the action of the group generated by $\TTT_{2\ptb}$.
Applying then \ref{P:app:ker:Xta} we find by standard linear theory a unique $\varphi\in C^{2,\beta}_{sym}(\cyl)$
which satisfies (i)
and moreover
$$
\|\varphi:L^2(\cyl/\TTT_{2\ptb},g)\|\le
C(b)\,\|E:C^{0,\beta}_\sym(\cyl,\chi)\|.
$$
The uniformity of geometry of $\Sph^1\times[-b-1,b+1]$ allows us to apply standard linear
theory to estimate the $\|\varphi:C^{2,\beta}(\Sph^1\times[-b-1,b+1],\chi)\|$ norm.
Applying then \ref{C:Lambda:u} with $\Lambda=\widetilde{\Lambda}$ to get an estimate on $\widetilde{\Lambda}$,
we conclude that
$$
\|\varphi:C^{2,\beta}_\sym(\cyl,\chi)\|
\le C(b,\beta) \,
\|E:C^{0,\beta}_\sym(\cyl,\chi)\|,
$$
which completes the proof of (ii) in this case.

We deal now with the general case:
By applying \ref{C:L:unique}.ii we obtain $V\in C^{2,\beta}(\widetilde{\Lambda})$
which satisfies (i) on $\widetilde{\Lambda}$ and also
$$
\|V:C^{2,\beta}(\widetilde{\Lambda},\chi)\|
\le
C(\beta)\,|\log\tb|^2\,\|E:C^{0,\beta}(\widetilde{\Lambda},\chi)\|,
$$
where we also used \ref{E:ell:bound}.
Let $\widetilde{E}\in C_\sym^{0,\beta}(\widetilde{\Lambda},\chi)$
be defined by requesting $\widetilde{E}=E$ on $\Sph^1\times[-b,b]$,
and on $\widetilde{\Lambda}$
$$
\widetilde{E}=\psi[1,0]\circ\xboth\, E+\,[\psi[0,1]\circ\xboth,\Lchi]V.
$$
We apply then the special case of the proposition we have already proven
with $\widetilde{E}$ instead of $E$,
to obtain an element of $C^{2,\beta}_{sym}(\cyl)$
which we will call $\widetilde{\varphi}$ instead of $\varphi$.
We define then $\varphi\in C^{2,\beta}_{sym}(M)$ by
by requesting $\varphi=\widetilde{\varphi}$ on $\Sph^1\times[-b,b]$,
and on $\widetilde{\Lambda}$
$$
\varphi=\widetilde{\varphi}+\,
\psi[0,1]\circ\xboth\, V.
$$
It is easy to chck then that (i) and (ii) hold.
(iii) follows easily from the continuous dependence of $\chi$ on the parameters.
\end{proof}

\section{Using the Geometric Principle to prescribe the extended substitute kernel}
\label{S:prescribe}
\nopagebreak

In this section we discuss how to prescribe elements of the extended substitute kernel
as required by our general approach \cite{kapouleas:wente,kapouleas:survey}.
To simplify the notation we adopt from now on the following:

\addtocounter{equation}{1}
\begin{convention}
\label{convention}
We fix some $\beta\in (0,1)$, $\gamma\in(2,3)$, $\beta'\in(0,\beta)$, and $\gamma'\in(\gamma,3)$.
We will be using $C$ to denote positive constants which may depend on $b,\beta,\beta',\gamma,\gamma'$
(not mentioned explicitly)
and any other constants mentioned explicitly.
Recall that 
$\tb$ is always assumed small enough in terms of any other constants---equivalently
by \ref{E:p:hat} $\mb$ large enough---in accordance with the proof of \ref{T:main}.
\end{convention}

\subsection*{Prescribing $\skernel[0]$}
$\phantom{ab}$
\nopagebreak

The prescription of the $w_{i,0}$'s as part of the original Lagrangian angle $\theta$
is done by the introduction of $\thetadisloc$ which is
due to the deformation controlled by the parameters $\zetabold$.
The construction had to be carried out at the nonlinear level,
that is the initial immersion itself had to be modified to $Y_\zetabold$.
This contrasts with the construction in the next subsection where the
deformations are introduced through functions depending on the appropriate parameters
and are studied at the linear level.
The nonlinear terms are included then in the error of linearizing which is dealt with
in the next section.
The way $\zetabold$ controls the prescription of the $w_{i,0}$'s
and related estimates are provided in the following proposition \ref{P:phiunder0},
especially \ref{P:phiunder0}.(iii).
We remark that \ref{P:phiunder0} is similar to part of
\cite[Proposition 6.7]{kapouleas:wente}.

\begin{prop}
\addtocounter{equation}{1}
\label{P:phiunder0}
There exists $\phiunder_\zetabold\in C_\SSS^{2,\beta}(\Stilde[0])$ and
$(\mu_{1,0},\mu_{2,0})\in \R^2$ 
such that the following hold where the initial immersion under consideration is $Y_\zetabold$:
\newline
(i).
$\Lcal \phiunder_\zetabold + \thetadisloc = \mu_{1,0} w_{1,0} + \mu_{2,0} w_{2,0} $
on $\widetilde{S}[0]$,
where $\thetadisloc$ is induced by $Y_\zetabold$.
\newline
(ii).
$\phiunder_\zetabold=0$ on $\partial \widetilde{S}[0].$
\newline
(iii).
$|\zetabold-(\mu_{1,0},\mu_{2,0})| \le
C\,|\zetabold|/|\log\tb|$.
\newline
(iv).
$\Vert \,\phiunder_\zetabold: C^{2,\beta}(\widetilde{S}[0],\chi)\,\Vert  \le
C\,|\zetabold|.$
\newline
(v).
$\Vert \,\phiunder_\zetabold: C^{2,\beta}(\Lambda[1],\chi,e^{-\gamma'\xunder})\,\Vert \le
C\,|\zetabold|.$
\newline
(vi).
$\phiunder_\zetabold$ is unique by its construction and depends continuously on $\zetabold$.
\end{prop}

\begin{proof}
By the construction of $Y_\zetabold$ and \ref{L:Ytilde}
there is $\underline{f}:\Sph^1\times[a,a+3]\to\R$
such that
the Legendrian perturbation of $\Ytilde[0]$ by $\underline{f}$ defined as in Appendix 
\ref{A:perturbation}
satisfies
$$
(\Ytilde[0])_{\underline{f}}=Y_\zetabold\circ D
$$
for some diffeomorphism $D$ from $\Sph^1\times[a,a+3]$ to a subset of $\Sph^1\times\R$,
where $D$ is a small perturbation of the identity map.
Clearly then $\underline{f}$ satisfies the following, where for (c) we use \ref{P:Leg:pert}:
\newline
(a). $\underline{f}=0$ on $\Sph^1\times[a,a+1]$.
\newline
(b). $\|\underline{f}:C^{2,\beta}(\Sph^1\times[a,a+3],\chi)\|\le C\,|\zetabold|^2$.
\newline
(c). $\|\Lcal_{\Ytilde[0]}\,\underline{f}-\thetadisloc\circ D:C^{2,\beta}(\Sph^1\times[a,a+3],\chi)\|\le
C\,|\zetabold|^2$,
where $\Lcal_{\Ytilde[0]}=\Delta_{\Ytilde[0]}+6$,
where $\Delta_{\Ytilde[0]}$ is the Laplacian with respect to the metric induced by $\Ytilde[0]$.
\newline
(d). 
$\|\underline{f}-\zeta'_1\,\ftbold\circ\Ytilde[0]-\zeta'_2\,\fqbold\circ\Ytilde[0]
:
C^{2,\beta}(\Sph^1\times[a+2,a+3],\chi)   \|
\le
C\,|\zetabold|^2$
(recall \ref{E:f:tq}).

Using \ref{L:Ytilde} and \ref{L:theta},
(c) and (d) are modified to the following:
\newline
(e). $\|\Lcal\,\underline{f}-\thetadisloc:C^{2,\beta}(\Sph^1\times[a,a+3],\chi)\|\le
C\,|\zetabold|^2$,
where $\Lcal$ is as in \ref{E:L}.
\newline
(f). 
$\|\underline{f}-\zeta'_1\,\ftbold\circ Y_\zetabold-\zeta'_2\,\fqbold\circ Y_\zetabold
:
C^{2,\beta}(\Sph^1\times[a+2,a+3],\chi)   \|
\le
C\,|\zetabold|^2$.
%

We define now $\phiunder'_\zetabold\in C^{2,\beta}(\widetilde{S}[0])$
by requesting the following:
\newline
(g).
$\phiunder'_\zetabold$ is appropriately symmetric in the sense of \ref{D:f:symmetries}.
\newline
(h).
$\phiunder'_\zetabold=0$ on $M_0$.
\newline
(i).
$\phiunder'_\zetabold=(1-\check{\psi})\underline{f}
+
\check{\psi}\,(\zeta'_1\,\ftbold\circ Y_\zetabold+\zeta'_2\,\fqbold\circ Y_\zetabold)$
on $\Sph^1\times[a,b]$,
where
$\check{\psi}=\psi[a+2,a+3]\circ t$,
where $t$ is the coordinate defined in \ref{D:M}.
\newline
(k).
$\phiunder'_\zetabold=\zeta'_1\,\ftbold\circ Y_\zetabold+\zeta'_2\,\fqbold\circ Y_\zetabold
-(1-\psi_{S[0]})
\underline{V}$
on $\Lambda[1]$,
where $\underline{V}$ is the solution to the Dirichlet problem on $\Lambda[1]$
given by
$\Lcal\underline{V}=0$ on $\Lambda[1]$,
$\underline{V}=0$ on $C^+[0]$,
and
$\underline{V}=\zeta'_1\,\ftbold\circ Y_\zetabold+\zeta'_2\,\fqbold\circ Y_\zetabold$
on $C^-[1]$.
Using \ref{C:Lambda:u} we have 
\newline
(l). $\|\underline{V}:C^{2,\beta}(S_2[0]\setminus S[0],\chi)\|\le
C\,|\zetabold|/ |\log\tb|$.

$\phiunder'_\zetabold$ has the following properties:
\newline
(m).
$\Lcal \phiunder'_\zetabold$ is supported on $S_1[0]$
(by \ref{P:killing})
and
$\phiunder'_\zetabold=0$ on $\partial \Stilde[0]$.
\newline
(n).
$\|\Lcal \phiunder'_\zetabold-\thetadisloc: C^{0,\beta}(S_1[0],\chi)\|\le
C\,|\zetabold|/ |\log\tb|$
follows from (e), (f), and (l).
\newline
(o).
$\|\phiunder'_\zetabold-\zeta'_1\,\ftbold\circ Y_\zetabold-\zeta'_2\,\fqbold\circ Y_\zetabold
:
C^{2,\beta}(C^+_1[0],\chi)   \|
\le
C\,|\zetabold|/ |\log\tb|$
follows from (l).

By applying then \ref{L:vskernel} and \ref{C:Lambda:solution} we determine $(\mu'_1,\mu'_2)\in\R^2$
such that
$|\zetabold-(\mu_1',\mu_2')| \le
C\,|\zetabold|/ |\log\tb|$
and 
$$
\Vert \,\phiunder'_\zetabold-\mu'_1 v_{1,0}-\mu'_2 v_{2,0}:
C^{2,\beta}(\Lambda[1],\chi,e^{-\gamma'\xunder})\,\Vert 
\le C \,|\zetabold|/ |\log\tb|.
$$
Using (m) and (n) above we apply \ref{L:wS:solution}
with $E=\rho^2(\Lcal \phiunder'_\zetabold-\thetadisloc)$
to obtain $\varphi$ and $w=\mu''_1w_{1,0}+ \mu''_2w_{2,0}$.
We define then
$$
\phiunder_\zetabold:=\varphi-\phiunder'_\zetabold+\mu'_1 v_{1,0}+\mu'_2 v_{2,0},
\qquad
(\mu_{1,0},\mu_{2,0}):=(\mu'_1,\mu'_2)+(\mu''_1,\mu''_2),
$$
and we complete the proof by appealing to the estimates above and the conclusions of \ref{L:wS:solution}.
\end{proof}

\subsection*{Prescribing $\skernel[n']$}
$\phantom{ab}$
\nopagebreak

The introduction of the other elements of the extended substitute kernel,
that is the elements of the substitute kernel,
can be done at the linear level.
The amount we have of these elements can be monitored by using
a linearized version of the balancing argument which 
amounts to using Green's second identity \cite{gilbarg}.

We proceed to define now $\phi_{i,n'}\in C^{0,\beta}(\widetilde{S}[n'])$
for $i=1,2$ and $n'$ as in \ref{E:nnprime}.
We first define 
$\phi'_{i,n'}\in C^{0,\beta}(\widetilde{S}[n'])$
by requiring that it is supported on $\Lambda[n']$ where
\addtocounter{theorem}{1}
\begin{equation}
\label{E:phiprime:inp}
\phi'_{1,n'}=\psi_{\widetilde{S}[n'-1]}\,V_1[\Lambda[n'],\widehat{c}_1\ell,0],
\qquad
\phi'_{2,n'}=\psi_{\widetilde{S}[n'-1]}\,V_2[\Lambda[n'],\widehat{c}_2\sinh{2\ell},0],
\end{equation}
where $\widehat{c}_1$ and $\widehat{c}_2$ are determined by requesting that
(recall \ref{E:Ccentral})
\addtocounter{theorem}{1}
\begin{equation}
\label{E:chat}
\int_{\widetilde{C}[n']}\frac{\partial \phi'_{i,n'}}{\partial t} \widehat{f}_i \, ds=-1.
\end{equation}
Using then \ref{E:Lambda:Vtilde}, \ref{P:Lambda:V}, \ref{E:ell:bound}, and \ref{E:ymin:tau},
we conclude that
\addtocounter{theorem}{1}
\begin{equation}
\label{E:chat:bounds}
C^{-1}<\widehat{c}_i< C,
\qquad
(i=1,2).
\end{equation}

We apply now \ref{L:wS:solution} with $E=-\rho^2\,\Lcal\phi'_{i,n'}$
to obtain
$\phi''_{i,n'}\in C^{0,\beta}(\widetilde{S}[n'])$
and $\underline{\widetilde{w}}_{i,n''}\in\skernel[n']$
satisfying (i)-(vii) in \ref{L:wS:solution}.
We define then
\addtocounter{theorem}{1}
\begin{equation}
\label{E:phiprimeprime:inp}
\phi_{i,n'}= \phi'_{i,n'}+ \phi''_{i,n'},
\end{equation}
and we have the following lemma:

\addtocounter{equation}{1}
\begin{lemma}
\label{L:phi:inp}
$\phi_{i,n'}\in C^{0,\beta}(\widetilde{S}[n'])$
and $\underline{\widetilde{w}}_{i,n'}\in\skernel[n']$
as defined above satisfy the following:
\newline
(i).
$\Lcal \phi_{i,n'} = \underline{\widetilde{w}}_{i,n'}$
on $\widetilde{S}[n']$.
\newline
(ii).
$\phi_{i,n'}$ vanishes on $\partial \widetilde{S}[n]$.
\newline
(iii).
$|\underline{\widetilde{w}}_{i,n'}-w_{i,n'}|\le C/|\log\tb|$.
\newline
(iv).
$\Vert \phi_{i,n'}: C^{2,\beta}(\widetilde{S}[n']\setminus\Lambda[n'],\chi) \Vert 
\le C.$
\newline
(v).
$\Vert \phi_{i,n'}: C^{2,\beta}(\Lambda[n'+1],\chi,e^{-\gamma'\xunder})\Vert 
\le C$.
\newline
(vi).
$\Vert \phi_{i,n'}: C^{2,\beta}(\Lambda[n'],\chi,e^{-2\xunder})) \Vert 
\le C \,\tb^{-2}$.
\newline
(vii).
$\phi_{i,n'}$
and $\underline{\widetilde{w}}_{i,n'}$
are unique by their construction.
\end{lemma}

\begin{proof}
(i), (ii), and (iii) follow from the definitions.
Applying \ref{P:Lambda:V}, \ref{E:Lambda:Vtilde}, and \ref{E:ell:bound} we get
\begin{equation*}
\begin{aligned}
\|\phi'_{i,n'}: C^{2,\beta}(\Lambda[n'],\chi,e^{-2\xunder})
\le&e^{2\ell},
\\
\|\Lchi\phi_{i,n'}: C^{2,\beta}(\widetilde{S}[n'],\chi)\|
\le&C.
\end{aligned}
\end{equation*}
Using then the estimates in \ref{L:wS:solution} we establish (iv), (v) and (vi).

It remains to prove (iii).
By Green's second identity \cite{gilbarg} we have
$$
\int_\Omega
(\Lcal \phi_{i,n'}\,\widehat{f}_j - \phi_{i,n'}\,\Lcal\widehat{f}_j) dg
=
\int_{\partial\Omega}
(\vec{N} \phi_{i,n'}\,\widehat{f}_j - \phi_{i,n'}\,\vec{N}\widehat{f}_j) dg,
$$
where $\Omega$ is an appropriate domain in $M$,
and $\vec{N}$ is the unit conormal to $\partial \Omega$ pointing to $M\setminus\Omega$.
We apply this with $\Omega=\Sph^1\times[(2n'-1)\ptb,(2n'+1)\ptb]$:
Using the invariance of the right hand side under conformal changes of the metric,
that $\Lcal \phi_{i,n'}=\underline{\widetilde{w}}_{i,n'}$
which is supported on $S_1[n']$,
that by \ref{P:killing} and \ref{E:f1f2} $\Lcal \widehat{f}_j=0$,
and that by \ref{E:Ccentral}, \ref{D:X:tau}, and \ref{P:conformal}
$\frac\partial{\partial t} \widehat{f}_j=0$ on $\widetilde{C}[n'']$,
we conclude
$$
\int_{S_1[n']}\underline{\widetilde{w}}_{i,n'}\widehat{f}_j\,dg=
\int_{\widetilde{C}[n'+1]}\frac{\partial \phi_{i,n'}}{\partial t} \widehat{f}_j \, ds
-
\int_{\widetilde{C}[n']}\frac{\partial \phi_{i,n'}}{\partial t} \widehat{f}_j \, ds.
$$
Using then \ref{E:chat}
and the estimates on $\phi''_{i,n'}$
we establish (iii).
\end{proof}

\subsection*{Prescribing $\skernel$ globally}
$\phantom{ab}$
\nopagebreak

We assume now $\xibold=\{\xi_{1,n},\xi_{2,n}\}_{n=0}^{2\mb-1}\in\R^{4\mb}$
given and we proceed to construct $\Phi_\xibold\in C^{0,\beta}_{sym}(M)$.
We assume that
\addtocounter{theorem}{1}
\begin{equation}
\label{E:xibold:bound}
\|\xibold\|_\gamma\le\cunder\tb,
\end{equation}
and we take
\addtocounter{theorem}{1}
\begin{equation}
\label{E:zetaxi}
\zetabold=(\xi_{1,0},\xi_{2,0})
\end{equation}
and define $\Phi'_\xibold:=\Ubold(\{\psi_{\widetilde{S}[n]}\,\phiunder[n]\})\in C^{2,\beta}_{sym}(M)$,
where $\phiunder[0]=\phiunder_\zetabold$
and for $n\ne0$
$\phiunder[n]=\xi_{1,n}\phi_{1,n}+ \xi_{2,n}\phi_{2,n}$.
Let
\addtocounter{theorem}{1}
\begin{equation}
\label{E:wunderprime}
\underline{w}':=\mu_{1,0}w_{1,0}+\mu_{2,0}w_{2,0}+\sum_{i,n'} \xi_{i,n'} \underline{\widetilde{w}}_{i,n'}\in\skernel,
\end{equation}
where $\mu_{i,0}$ is as in \ref{P:phiunder0} and $\underline{\widetilde{w}}_{i,n'}$ as in \ref{L:phi:inp}.
We apply then \ref{P:solution} to obtain 
\addtocounter{theorem}{1}
\begin{equation}
\label{E:Phi2p}
(\Phi''_\xibold,\underline{w}'')=\Rcal_M(\rho^2(-\Lcal\Phi'_\xibold+\underline{w}')).
\end{equation}
We define then 
\addtocounter{theorem}{1}
\begin{equation}
\label{E:Phixi}
\Phi_\xibold:=\Phi'_\xibold+\Phi''_\xibold,
\qquad
\underline{w}_\xibold:=\underline{w}'+\underline{w}'',
\end{equation}
and we have the following Proposition:

\addtocounter{equation}{1}
\begin{prop}
\label{P:Phixi}
$\Phi_\xibold\in C^{2,\beta}_{sym}(M)$ and $\underline{w}$,
defined as above for the immersion $Y_\zetabold$ where $\zetabold$ is as in \ref{E:zetaxi},
depend continuously on $\xibold$ and satisfy the following:
\newline
(i).
$\Lcal \Phi_\xibold + \thetadisloc=\underline{w}_\xibold$ on $M$.
\newline
(ii).
$\|\Phi_\xibold\|_{2,\beta,\gamma}
\le
C\,\|\xibold\|_\gamma$.
\newline
(iii).
$\|\underline{w}_\xibold-w_\xibold \|_\gamma
\le
\,C\,\|\xibold\|_\gamma/|\log\tb|$,
where
$w_\xibold:=\sum_{i,n}\xi_{i,n}w_{i,n}$.
\end{prop}

\begin{proof}
Using the definitions, \ref{P:phiunder0}, and \ref{L:phi:inp} we clearly have the following:
\newline
(a). 
$\|\Phi'_\xibold\|_{2,\beta,\gamma}
\le
C\,\|\xibold\|_\gamma$.
\newline
(b). 
$\|\Lchi\Phi'_\xibold\|_{2,\beta,\gamma}
\le
C\,\tau^{\min(\gamma-2,\gamma'-\gamma)}\,\|\xibold\|_\gamma$.
\newline
(c). 
$\|\underline{w}'-\underline{w}_\xibold\|_\gamma
\le
\,C\,\|\xibold\|_\gamma/|\log\tb|$.

Using then the estimates provided in \ref{P:solution} and (b) we complete the proof.
\end{proof}

\section{The main results}
\label{S:results}
\nopagebreak

In this section we prove our main results.

\subsection*{The nonlinear terms}
$\phantom{ab}$
\nopagebreak

Using a rescaling argument we prove now a local estimate for the nonlinear terms:

\addtocounter{equation}{1}
\begin{lemma}
\label{L:quadratic}
Consider $\Xunderta:\cyl\to\Sph^5$
where $\tau$ and $\alpha$ satisfy 
\ref{E:tau:alp:interval}.
Let $D\subset\cyl$ be a disc of radius $1$ with respect to $\chi$
and center some $p\in\cyl$.
If $f,v\in C^{2,\beta}(D,\chi)$ satisfy
$$
\Vert f:C^{2,\beta}(D,\chi)\Vert  < c_2\,\rhounder^2_\tau(p),
\qquad
\Vert v:C^{2,\beta}(D,\chi)\Vert  < c_2\,\rhounder^2_\tau(p),
$$ 
then the Legendrian perturbations
$(\Xunder_{\tau,\alpha)_f},(\Xunder_{\tau,\alpha})_{f+v}:D\to\Sph^5$
are well defined as in Appendix
\ref{A:perturbation}
and satisfy
\begin{multline*}
\Vert \rhounder_\tau^2\,\theta_{f+v} - \rhounder_\tau^2\,\theta_{f} - \Lchi v:C^{0,\beta}(D,\chi)\Vert
 \le
\\
C\,\rhounder^{-2}_\tau(p)\,
(\Vert f:C^{2,\beta}(D,\chi)\Vert +\Vert v:C^{2,\beta}(D,\chi)\Vert )\,
\Vert v:C^{2,\beta}(D,\chi)\Vert ,
\end{multline*}
where $\theta_f$ and $\theta_{f+v}$ are the Lagrangian angles of 
$(\Xunder_{\tau,\alpha})_f$
$(\Xunder_{\tau,\alpha})_{f+v}$
respectively
and $\rhounder_\tau$ is as in \ref{D:chi}.
\end{lemma}

\begin{proof}
By \ref{P:immersion}.(ii-iii) and \ref{P:conformal}.(iv)
we conclude that there is a universal constant $c'>0$ such that
$$
X:=R\,\Xunderta:D\to\Sph^5(R)
$$
satisfies the assumptions of \ref{P:Leg:pert},
where $R:=c'/\rhounder_\tau(p)$.
Since by taking $\phi=R^2f$ and $\varphi=R^2v$ we have 
$$
\rhounder_\tau^2\,\theta_{f+v} - \rhounder_\tau^2\,\theta_{f} - \Lchi v
=\rhounder^2_\tau\,(
\theta_{\phi+\varphi} - \theta_{\phi} - (\Delta + 6R^{-2})\varphi
\,),
$$
where $\theta_\phi$ and $\theta_{\phi+\varphi}$ are the Lagrangian angles of
$X_\phi:D\to\Sph^5(R)$ and $X_{\phi+\varphi}:D\to\Sph^5(R)$ respectively,
and $\Delta$ is the Laplacian induced by $X$ as in \ref{P:Leg:pert},
the lemma follows by applying \ref{P:Leg:pert} and \ref{P:conformal}.(iv).
\end{proof}

Using the local estimate for the nonlinear terms we just proved it is easy
to obtain global estimates suitable for our purposes:

\addtocounter{equation}{1}
\begin{corollary}
\label{C:quadratic}
There is a universal constant $c_3>0$ such that if $f,v\in C^{2,\beta}_\sym(M)$
satisfy
$$
\|f:C^{2,\beta}_\sym(M,\chi)\|\le c_3\,\tb,
\qquad
\|v\|_{2,\beta,\gamma}\le c_3,
$$
then the Legendrian perturbations $(Y_{\zetabold})_f$
and $(Y_{\zetabold})_{f+v}$ of $Y_\zetabold$
are well-defined as in Appendix
\ref{A:perturbation}
and satisfy
$$
\Vert \rho^2\,\theta_{f+v} - \rho^2\,\theta_{f} - \Lchi v\Vert_{0,\beta,\gamma}
 \le
C\,
(   \tb^{-1}\, \Vert f:C^{2,\beta}(D,\chi)\Vert +
\, \Vert v \Vert_{2,\beta,\gamma}      \,   )\,
 \Vert v \Vert_{2,\beta,\gamma}      \,   )\,
$$
where $\theta_f$ and $\theta_{f+v}$ are the Lagrangian angles of 
$(Y_{\zetabold})_f$
and
$(Y_{\zetabold})_{f+v}$
respectively
and $\rho$ is as in \ref{D:chi}.
\end{corollary}

\begin{proof}
Using the definitions of the norms involved, \ref{P:conformal}.(iv),
and \ref{L:quadratic}, the corollary follows.
\end{proof}

\addtocounter{equation}{1}
\begin{corollary}
\label{C:quadratic:Xta}
There is a universal constant $c_4>0$ such that if $f,v\in C^{2,\beta}_\sym(\cyl)$
satisfy
$$
\|f:C^{2,\beta}_\sym(\cyl,\chi)\|\le c_3\,\tb,
\qquad
\|v:C^{2,\beta}_\sym(\cyl,\chi)\|\le c_3\,\tb,
$$
then the Legendrian perturbations
$(\Xunder_{\tau,\alpha})_f,(\Xunder_{\tau,\alpha})_{f+v}:\cyl\to\Sph^5$
are well defined as in Appendix
\ref{A:perturbation}
and satisfy
\begin{multline*}
\Vert \rhounder_\tau^2\,\theta_{f+v} - \rhounder_\tau^2\,\theta_{f} - \Lchi v:C^{0,\beta}_\sym(\cyl,\chi)\Vert
 \le
\\
C\,\tb^{-1}\,
(\Vert f:C_\sym^{2,\beta}(\cyl,\chi)\Vert +\Vert v:C_\sym^{2,\beta}(\cyl,\chi)\Vert )\,
\Vert v:C_\sym^{2,\beta}(\cyl,\chi)\Vert ,
\end{multline*}
where $\theta_f$ and $\theta_{f+v}$ are the Lagrangian angles of 
$(\Xunder_{\tau,\alpha})_f$
and
$(\Xunder_{\tau,\alpha})_{f+v}$
respectively
and $\rhounder_\tau$ is as in \ref{D:chi}.
\end{corollary}

\begin{proof}
Using the definitions of the norms involved and \ref{L:quadratic} the corollary follows.
\end{proof}

\subsection*{Correcting $\Xunderta$ to a special Legendrian immersion}
$\phantom{ab}$
\nopagebreak

\addtocounter{equation}{1}
\begin{theorem}
\label{T:Xta}
There is $\funder=\funder_{\tau,\alpha}\in C^\infty_\sym(\cyl)$
such that $(\Xunder_{\tau,\alpha})_{\funder}:\cyl\to\Sph^5$
is well-defined as in Appendix
\ref{A:perturbation}
and is special Legendrian.
Moreover $\funder$ depends continuously on $\tau,\alpha$ and satisfies
$$
\|\funder_{\tau,\alpha}:C^{2,\beta}_\sym(\cyl,\chi)\|
\le
C\,\cunder\tb^2\,|\log\tb|.
$$
\end{theorem}

\begin{proof}
Applying \ref{P:solution:Xta} we find $f:=-\Rcal_\cyl(\rhounder_\tau^2\theta)$
which satisfies then
\newline
(a). 
$\Lchi f =-\rhounder_\tau^2\theta$.
\newline
(b).
$\|f:C^{2,\beta}_\sym(\cyl,\chi)\|
\le C \,
\cunder \,\tb^2 \,|\log\tb|
$,
where we used \ref{E:theta} also.

Applying then \ref{C:quadratic:Xta} with $0$ and $f$ in place of $f$ and $v$ respectively
we conclude
\newline
(c).
$\Vert \rhounder_\tau^2\,\theta_{f} - \rhounder_\tau^2\,\theta - \Lchi f:C^{0,\beta}_\sym(\cyl,\chi)\Vert
\le C \,
\cunder^2 \,\tb^3 \,|\log\tb|^2
\le \,\tb^2$.

We define now a map $\Jcal:B\to B$,
where
$$
B:=\{v\in C^{2,\beta}_\sym(\cyl):\|v:C^{2,\beta}_\sym(\cyl,\chi)\|\le\tb^2\},
$$
as follows:
By applying \ref{C:quadratic:Xta} again we have that if $v\in B$, then
\newline
(d).
$
\Vert \rhounder_\tau^2\,\theta_{f+v} - \rhounder_\tau^2\,\theta_{f} - \Lchi v:C^{0,\beta}_\sym(\cyl,\chi)\Vert
\le C \,
\cunder \,\tb^3 \,|\log\tb|
\le \,\tb^2$.

Applying then \ref{P:solution:Xta} again we define
$$
\Jcal v:=-\Rcal_\cyl
(\rhounder_\tau^2\,\theta_{f+v} - \rhounder_\tau^2\,\theta_{f} - \Lchi v)
-\Rcal_\cyl
(\rhounder_\tau^2\,\theta_{f} - \rhounder_\tau^2\,\theta - \Lchi f),
$$
which then satisfies by (a), (c), and (d), that
\newline
(e).
$\Lchi\Jcal v=-\rhounder_\tau^2\,\theta_{f+v} + \Lchi v$.
\newline
(f).
$\|\Jcal v :C^{2,\beta}_\sym(\cyl,\chi)\Vert
\le C \,
\cunder^2 \,\tb^3 \,|\log\tb|^4
\le\,\tb^2$,
which implies that $\Jcal v\in B$ as we need.

We prove now that $\Jcal$ is a contraction:
If $v_1,v_2\in B$, we have
$$
\Jcal v_1- \Jcal v_2=
\Rcal_\cyl
(\rhounder_\tau^2\,\theta_{f+v_2} - \rhounder_\tau^2\,\theta_{f+v_1} - \Lchi (v_2-v_1)\,).
$$
Applying then 
\ref{C:quadratic:Xta} again with $f+v_1$ and $v_2-v_1$ instead of $f$ and $v$ respectively,
and \ref{P:solution:Xta} afterwards,
we obtain
\newline
(g).
$\|\Jcal v_2-\Jcal v_1 :C^{2,\beta}_\sym(\cyl,\chi)\Vert
\le C \,
\cunder \,\tb \,|\log\tb|^3
\,
\| v_2- v_1 :C^{2,\beta}_\sym(\cyl,\chi)\Vert$.
\newline
It follows that $\Jcal$ is a contraction and since $B$ is complete there is a unique fixed
point of $\Jcal$ which we call $\underline{v}$.
By defining then $\funder=f+\underline{v}$ we complete the proof.
\end{proof}

\subsection*{The main theorem}
$\phantom{ab}$
\nopagebreak

We are ready to state and prove our main theorem.
While the proof depends on more or less everything in this paper,
the statement uses only some of the definitions,
in particular \ref{E:p:hat}, \ref{D:M}, \ref{D:Yzeta}, \ref{D:chi}, 
\ref{D:f:symmetries}, \ref{D:globalnorm},
and the construction in Appendix 
\ref{A:perturbation}:

\addtocounter{equation}{1}
\begin{theorem}
\label{T:main}
There is a constant $C$ such that if $\mb$ is large enough,
there is $\zetabold\in\R^2$ with $|\zetabold|<C/\mb$
and $f\in C^\infty_\sym(M)$ with $\|f\|_{2,\beta,\gamma}\le C/\mb$
such that $(Y_{\zetabold})_f:M\to \Sph^5$ is well-defined
as in Appendix
\ref{A:perturbation}
and is a special Legendrian immersion satisfying the symmetries in
\ref{E:symmetries}.
\end{theorem}

\begin{proof}
We assume $\tb$ fixed and small enough as in \ref{convention},
which is equivalent to $\mb$ being large enough.
We define
a map $\Jcal : B \to B$ where
$$
B:=\{u\in C^{2,\beta}_\sym(M):\|u\|_{2,\beta,\gamma}\le\tb^{3/2}\}
\times
\{\xibold\in\R^{4\mb}:\|\xibold\|_\gamma\le \cunder\tb\},
$$
as follows:
We assume $(u,\xibold)\in B$ given.
Note then that $\xibold$ determines $\zetabold$ by \ref{E:zetaxi},
which determines then $\tau$ and $\alpha$ by \ref{E:sigma:range} and \ref{E:sigma},
and $Y_\zetabold:M\to\Sph^5$ as in section \ref{S:init:surf}.
We define $\fundertilde\in C^{2,\beta}_\sym(M)$ by requiring
$\fundertilde=\psihat\,\funder_{\tau,\alpha}$ on $M_1$,
and $\fundertilde=0$ on $M_0\setminus\cup_{j=1}^g M_j$---recall
\ref{T:Xta}, \ref{D:psihat}, and \ref{D:M}.
We have then
\newline
(a).
$
\|\fundertilde:C^{2,\beta}_\sym(M,\chi)\|
\le
C\,\cunder\tb^2\,|\log\tb|.
$
\newline
(b). 
$\theta_{\fundertilde}- \thetadisloc - \thetaglue $
is supported on $S[0]$ by using 
\ref{E:decomposition}
and satisfies by \ref{P:Leg:pert}
$$
\|
\theta_{\fundertilde}- \thetadisloc - \thetaglue 
:C^{0,\beta}_\sym(M,\chi)\|
\le
C\,\cunder\tb^2\,|\log\tb|.
$$

We apply now
\ref{P:solution}
to obtain $(\varphi,w)=\Rcal_M (-\rho^2 \theta_{\fundertilde} +\rho^2 \thetadisloc)$
which satisfies
\newline
(c).
$\Lcal \varphi =- \theta_{\fundertilde} + \thetadisloc + w$ on $M$.
\newline
(d).
$\|\varphi\|_{2,\beta,\gamma}\le C\,\tb$,
where we used \ref{L:theta}.(i) and (b) to estimate
$-\theta_{\fundertilde} +\thetadisloc=-\thetaglue-(\theta_{\fundertilde}- \thetadisloc - \thetaglue )$.
\newline
(e)
$\|w\|_\gamma\le C\,\tb$.

We define $v:=\Phi_\xibold + \varphi - u$ and we have then by appealing to \ref{P:Phixi},
(c), (d), and (e), that
\newline
(f).
$\Lcal v = -\theta_{\fundertilde} +  w+\underline{w}_\xibold-\Lcal u$.
\newline
(g).
$\|v\|_{2,\beta,\gamma}\le C\,\cunder\,\tb$.
\newline
(h).
$\|w + \underline{w}_\xibold - w_\xibold \|_\gamma\le \cunder\,\tb/4$,
where for this we have to choose $\cunder$ large enough in terms of $C$.

We apply
\ref{P:solution}
to obtain $(u',w')=\Rcal_M (\rho^2\,\theta_{\fundertilde+v} - \rho^2\,\theta_{\fundertilde} - \Lchi v)$
which satisfies by using (f)
\newline
(i).
$\Lcal u'= \theta_{\fundertilde+v} - w - \underline{w}_\xibold + \Lcal u +w'$.

Using then (a), (g), \ref{C:quadratic}, and \ref{P:solution} we conclude that
\newline
(j).
$\|u'\|_{2,\beta,\gamma}\le C\,\cunder^2\,\tb^2\,|\log\tb|\le\tb^{3/2}$.
\newline
(k).
$\|w'\|_{\gamma}\le C\,\cunder^2\,\tb^2\,|\log\tb|\le\tb^{3/2}$.

We define
$\boldsymbol{\mu}=\{\mu_{1,n},\mu_{2,n}\}_{n=0}^{2\mb-1}\in\R^{4\mb}$
by
\newline
(l).
$
\sum_{i,n}\mu_{i,n}w_{i,n}
=
-w - \underline{w}_\xibold + w_\xibold +w'.
$

We have then by (h) and (k) 
\newline
(m).
$\|\boldsymbol{\mu} \|_\gamma \le \cunder\,\tb$.

By defining then $\Jcal(u,\xibold)=(u',\boldsymbol{\mu})$
we have by (j) and (m) that $\Jcal(u,\xibold)\in B$.
Therefore $\Jcal:B\to B$ is well defined.
$B$ is clearly a compact convex subset of $C^{2,\beta'}_\sym(M)\times\R^{4\mb}$
for some $\beta'\in(0,\beta)$,
and it is easy to check that $\Jcal$ 
is a continuous map in the induced topology.
By Schauder's fixed point theorem
\cite[Theorem 11.1]{gilbarg} then,
there is $(\underline{u},\underline{\xibold})\in B$
such that $\Jcal(\underline{u},\underline{\xibold})=(\underline{u},\underline{\xibold})$.
It follows by (i) that $\theta_{\fundertilde+v} - w - \underline{w}_\xibold + w'=0$,
and by (l) that
$-w - \underline{w}_\xibold + w_\xibold +w '=w_\xibold $.
Since the smoothness follows by standard regularity theory,
the proof is completed by taking
$\zetabold=(\underline{\xi}_{1,0},\underline{\xi}_{2,0})$
and $f=\fundertilde+v$.
\end{proof}

\appendix

\section{Jacobi elliptic functions and elliptic integrals}
\label{A:elliptic}
\nopagebreak
  
This appendix recalls properties of the Jacobi
elliptic functions and elliptic integrals needed in \S \ref{S:u1:tori}.
For a more leisurely
description of elliptic functions we refer the reader to
\cite{lawden}.

Let $\sn{(t,k)}$, the \textit{Jacobi sn-noidal function with modulus} $k \in [0,1]$, be the unique solution of the equation
\addtocounter{theorem}{1}
\begin{equation}
  \label{sn-noid}
  \dot{z}^2= (1 -z^2)(1-k^2 z^2)
\end{equation}
with $z(0)=0,\  \dot{z}(0)=1$. It follows from this definition that in the case $k=0$ we have
$\sn{(t,0)}= \sin{t}$ and that for $k=1$ we have $\sn(t,1)=\tanh{t}$. By analogy with the trigonometric functions
there is also a \textit{Jacobi cn-noidal function} $\cnd{(t,k)}$ which satisfies
$$ \cnd^2{(t,k)} = {1 - \sn^2{(t,k)}}.$$
There is another Jacobi elliptic function $\dn{}$ satisfying
$$ \dn^2{(k,t)} = {1-k^2 \sn^2{(t,k)}}.$$
Using the definition of $\sn$ given in \ref{sn-noid} and the
relationships between the squares of the other Jacobi elliptic
functions  we find that
$$ \frac{d}{dt} \sn{t}  =  \cnd{t} \dn{t}, \quad  \frac{d}{dt} \cnd{t}  =  -\sn{t} \dn{t}, \quad \frac{d}{dt} \dn{t}  =  -k^2 \sn{t} \cnd{t}.$$
The period of $\sn{(t,k)}$ and $\cnd{(t,k)}$ is $4K(k)$, while
$\dn{(t,k)}$ has period $2K(k)$, where $K$ is the \textit{complete
elliptic integral of the first kind} defined by
$$ K(k) = \int_0^{\pi/2} \frac{dx}{\sqrt{1-k^2 \sin^2{x}}} = \sn^{-1}(1).$$
Similarly, the \textit{complete elliptic integral of the second kind} $E$
is defined by
$$ E(k) = \int_0^{\pi/2} {\sqrt{1-k^2 \sin^2{x}}} \ dx = \int_0^K \dn^2{t}\ dt.$$
For $k\in (0,1),$
$K$ is a positive strictly increasing function of $k$ and
$E$ is a positive strictly decreasing function of $k$.
The derivatives of $K$ and $E$ with respect to $k \in (0,1)$ are \cite[\S 3.8]{lawden}
\begin{eqnarray}
\label{dK:dk}
\frac{dK}{dk} & = & \frac{1}{k{k'}^2}(E-{k'}^2K),\\
\label{dE:dk}
\frac{dE}{dk} & = & \frac{1}{k}(E-K),
\end{eqnarray}
respectively, where ${k'}=\sqrt{1-k^2}$ is the so-called \textit{complementary modulus} to $k$.

$K$ and $E$ can be analytically continued to complex
values of $k$ using the fact that they both  satisfy linear
second-order differential equations, whose solutions exist for
complex values of $k$ \cite[p. 75]{lawden}.
$K(k)$ turns out to be analytic over the whole $k$-plane except
for logarithmic branch points at $k=\pm 1$. In the neighborhood
of $k=1$, it has the following expansion
\cite[p. 244]{lawden}
\addtocounter{theorem}{1}
\begin{equation}
\label{Kexpand}
K(k) = \frac{1}{\pi}K'(k)
\ln{\left(\frac{8}{h}\right)}-\frac{1}{4}h - \frac{7}{32}h^2 +
O(h^3)
\end{equation}
where $h=1-k$ and
\addtocounter{theorem}{1}
\begin{equation}
\label{Kpexpand} K'(k) = \frac{\pi}{2} \left(1 + \frac{1}{2}h +
\frac{5}{16}h^2 + O(h^3) \right).
\end{equation}
$E(k)$ also has logarithmic branch points at $k=\pm1$. In
a neighborhood of $k=1$ its expansion is \cite[p. 244]{lawden},
\addtocounter{theorem}{1}
\begin{equation}
\label{Eexpand}
E(k) = \frac{1}{\pi}J'(k) \ln{\left(\frac{8}{h}\right)} + 1 -
\frac{1}{2}h - \frac{5}{16}h^2 + O(h^3)
\end{equation}
where
\addtocounter{theorem}{1}
\begin{equation}
\label{Jpexpand}
J'(k) = \frac{\pi}{2}\left( h + \frac{1}{4}h^2 + O(h^3)\right).
\end{equation}

There are also \textit{incomplete elliptic integrals of the first and
second kinds} $F$ and $D$ defined by
$$ F(\phi,k) = \int_0^{\phi} \frac{dx}{\sqrt{1-k^2 \sin^2{x}}} = \int_0^u dt, \quad
 D(\phi,k) = \int_0^{\phi} {\sqrt{1-k^2 \sin^2{x}}} \ dx = \int_0^u \dn^2{t}\  dt$$
respectively, where $\sin{\phi}=\sn{u}$. Clearly $F(\pi/2,k)=K(k)$ and $D(\pi/2,k)=E(k)$.
Expanding the previous two integrands in ascending powers of $k^2$
and integrating term by term, we see that
\addtocounter{theorem}{1}
\begin{equation}
\label{Fexpand}
F(\phi,k) = \phi + \frac{1}{4}k^2 (\phi-\sin{\phi}\cos{\phi}) +
\ldots
\end{equation}
and
\addtocounter{theorem}{1}
\begin{equation}
\label{Dexpand}
D(\phi,k) = \phi - \frac{1}{4}k^2 (\phi-\sin{\phi}\cos{\phi}) +
\ldots.
\end{equation}

We define the \textit{incomplete elliptic integral of the third kind} by
\addtocounter{theorem}{1}
\begin{equation}
\label{ellipint3}
\Lambda(\phi,\alpha,k) = \int_0^{\phi} \frac{dx}{(1-\alpha^2\sin^2{x})\sqrt{1-k^2\sin^2{x}}}
= \int_0^u{\frac{dt}{1-\alpha^2\sn^2{t}}} = \Lambda(u,\alpha,k),
\end{equation}
where $\phi$ and $u$ are related by $\sn{u} = \sin{\phi}$.
When $\phi= \frac{\pi}{2}$ or $u=K$, we obtain the \textit{complete elliptic integral of the third kind}, $\Lambda(K,\alpha,k)$.

$\Lambda(K,\alpha,k)$ can be expressed in terms of $K$, $E$, $F$ and $D$.
This expression has four different forms depending
on whether $\alpha^2-k^2$ is positive or negative,
$\alpha$ is bigger than one or $\alpha$ is real.
In the case $0<k<\alpha<1$, the expression takes the form \cite[p. 76]{lawden},
\addtocounter{theorem}{1}
\begin{equation}
\label{compellipint3}
\Lambda(K,\alpha,k)= c(\alpha,k)
\left[ K(k)D(\phi,k')-K(k)F(\phi,k') + E(k)F(\phi,k') \right]
\end{equation}
where $k'$ is the complementary modulus to $k$,
$$c(\alpha,k) = \frac{\alpha}{\sqrt{(\alpha^2-k^2)(1-\alpha^2)}} \quad \textrm{and\quad }
 \sin{\phi}=\frac{\sqrt{\alpha^2-k^2}}{\alpha k'}.$$
For the case $\alpha^2=-\beta^2$ for some $\beta>0$, the analogous expression is \cite[p. 76]{lawden},
\addtocounter{theorem}{1}
\begin{equation}
\label{compellipint3b}
\Lambda(K,i\beta,k) = a(\beta,k)K + b(\beta,k) \left[ KD(\omega,k')-KF(\omega,k')+EF(\omega,k')\right]
\end{equation}
where
$$ a(\beta,k) = \frac{k^2}{\beta^2+k^2}, \quad b(\beta,k)= \frac{\beta}{\sqrt{(\beta^2+k^2)(\beta^2+1)}} \quad \textrm{and\ }
\sin{\omega} = \frac{\beta}{\sqrt{\beta^2+k^2}}.
$$

Finally, the \textit{Heuman Lambda function} $\Lambda_0(\phi,k)$ \cite[150.03]{byrd:friedman} is defined by
\addtocounter{theorem}{1}
\begin{equation}
\label{hlambda:defn}
\Lambda_0 (\phi,k): = \frac{2}{\pi} \left[K(k)D(\phi,k')-K(k)F(\phi,k') + E(k)F(\phi,k') \right].
\end{equation}
This particular combination of elliptic integrals of the first and second kinds appeared in \ref{compellipint3}.
Some basic properties of $\Lambda_0$ are listed in \cite[150.01--154.01]{byrd:friedman}.
According to \cite[710.11,730.04]{byrd:friedman}, the derivatives of $\Lambda_0$ with respect to $k$ and $\phi$ are
\addtocounter{theorem}{1}
\begin{equation}
\label{dk:hlambda}
\frac{\partial}{\partial k} \Lambda_0(\phi,k) = 2\frac{(E-K)\sin{\phi}\cos{\phi}}{\pi k \sqrt{1-{k'}^2\sin^2{\phi}}} < 0,
\end{equation}
and
\addtocounter{theorem}{1}
\begin{equation}
\label{dphi:hlambda}
\frac{\partial}{\partial\phi}\Lambda_0(\phi,k) = 2\frac{E - {k'}^2 \sin^2{\phi} K}{\pi \sqrt{1-{k'}^2\sin^2{\phi}}} > 0
\end{equation}
respectively. 

\section{Interpolation of Legendrian immersions}
\label{A:interpolation}
\nopagebreak

We often need to transit from one Legendrian immersion in $\Sph^5$ to another, where both immersions 
are close to the same totally geodesic Legendrian $\Sph^2$. To make this systematic assume we are given
two Legendrian immersions $X_0: \Omega_0 \rightarrow \Sph^5$ and 
$X_1: \Omega_1 \rightarrow \Sph^5$ where $\Sph^1 \times [a_1,a_2] \subset \Omega_i \subset \Sph^1 \times \R$
for $i=0,1$. We assume that on $\Sph^1 \times [a_1,a_2]$ both immersions are $C^1$-close
to each other and to the set $\Sph^2_{e_1,e_2,e_3} \subset \Sph^5$, 
where $e_1, e_2, e_3$ is a unitary basis of $\C^3$
and $\Sph^2_{e_1,e_2,e_3}:= \Sph^5 \cap P$,
where $P:=\left<e_1,e_2,e_3\right>_\R$.
We define then a Legendrian immersion $X$
$$X= \join{[X_0,X_1;a_1,a_2;e_1,e_2,e_3]}: \Omega \ra \Sph^5$$
as follows, 
where $\Omega$ is the union of  $\Omega_0 \cap \Sph^1 \times (-\infty,a_2)$ with 
$\Omega_1 \cap \Sph^1 \times (a_1,\infty)$. 
First, on $\Omega_0 \cap \Sph^1 \times (-\infty,a_1']$ we define $X$ to be $X_0$
and on $\Omega_1 \cap \Sph^1 \times [a'_2,\infty)$ we define $X$ to be $X_1$, 
where $(a'_1,a'_2)$ is the middle third of $(a_1,a_2)$.

We now proceed to define $X$ on $\Sph^1 \times (a'_1,a'_2)$.
Consider the conical Lagrangian immersions $Y_i: \Sph^1 \times (a_1,a_2) \times \R^+ \ra \C^3$
defined by $Y_i(e^{is},t,r) = rX_i(e^{is},t)$ for $i=0,1$. 
Let $Y_i = Y_i^{\top} + Y_i^{\perp}$ be the orthogonal decomposition where 
$Y_i^{\top}\in P$ and
$Y_i^{\perp}$ is perpendicular to the $3$-plane $P$.
Assuming that $Y_i^{\top}$ is a diffeomorphism 
from its domain to a radial open subset $\Omega_i'$ of the $3$-plane $P$
and using the fact the image of $Y_i$ is a Lagrangian cone close to $\Omega_i'$
there is a function $f_i: \Omega_i' \ra \R$ such that the
image of $Y_i$ coincides with the graph of the $1$-form $df_i$ 
(the closedness of the graphing $1$-form follows from the Lagrangian condition and
exactness then follows using the homogeneity).

We assume that the convex hull $\Omega''$ of $Y_0^{\top}(\Sph^1 \times (a'_1,a'_2)\times \R^+)$
satisfies $\Omega'' \subset \Omega_0'\cap \Omega_1'$. Let $\psi: \Omega_0' \ra \R$ and 
$f: \Omega'_0 \cap \Omega'_1 \ra \R$ be defined
by 
$$ \psi(Y_0^{\top}(e^{is},t,r)) = \psi[a'_1,a'_2](t) 
 \quad \text{and} \quad f = (1-\psi)f_0 + \psi f_1.$$
Since $f_0$ and $f_1$ are homogeneous of degree $2$ and $\psi$ is homogeneous of degree $0$, 
then $f$ is also homogeneous of degree $2$ and so $\graph(df)$ is a Lagrangian cone.
Let $Y^{\top}: \Sph^1 \times (a'_1,a'_2) \times \R^+ \ra P$ be defined by
$$Y^{\top}(e^{is},t,r) = \psi[a_2',a_1'](t) \, Y_0^{\top}(e^{is},t,r) +
\psi[a_1',a_2'](t) \, Y_1^{\top}(e^{is},t,r),$$
and $Y: \Sph^1 \times (a'_1,a'_2) \times \R^+ \ra \C^3$ be defined by
$$ Y = Y^{\top} + (df)^{\perp} \circ Y^{\top}.$$
Clearly, $Y$ is a Lagrangian conical (homogeneous of degree $1$) immersion, and so
we can define the Legendrian immersion $X: \Sph^1 \times (a_1',a_2') \ra \Sph^5$ 
by taking $X = Y/|Y|$.

Note that $X$ depends smoothly on $X_0$ and $X_1$, and when $X_0=X_1$ then $X=X_0$.

\section{Perturbation of Legendrian immersions}
\label{A:perturbation}
\nopagebreak
Given a Legendrian immersion $X:M \ra \Sph^5(R)$ of a surface $M$ and a $C^{2,\beta}$  
function $f:M \ra \R$,  with $\beta\in (0,1)$, 
 we will construct a new Legendrian immersion $X_f: M \ra \Sph^5(R)$.
This construction is an adaptation of the one used in \cite{butscher,butscher:thesis,y:lee} in 
the Lagrangian setting. Since the construction is local we can assume that $X$ is an embedding.
Then we first extend the function $f$ from $M$ to a tubular neighborhood $\Omega$ 
of $X(M)$ in $\Sph^5(R)$ by requiring
that if $v$ is a normal vector to $X(M)$ at $X(p)$ then 
$$f(\exp{(X(p),v\,;\Sph^5(R))} = f(p).$$
We extend  $f$ to the cone $C(\Omega)$ over $\Omega$
by requiring that $f$ is homogeneous of degree $2$. Let $V$ be the $C^{1,\beta}$ Hamiltonian vector
field on $C(\Omega)$ defined by 
$$ V = -J\nabla f.$$
Assuming that $f$ is small enough, we can flow $X$ by $V$ for unit time to an immersion 
$X_f:M \ra \Sph^5(R)$, that is, 
$$ X_f = Y(1,X(p)) \quad \text{for} \quad p\in M,$$
where $Y(t,q)$ is the flow defined by
$$\frac{\partial Y}{\partial t}(t,q) = V(Y(t,q)) \quad \text{and} \quad Y(0,q)=q \quad \text{for} \quad q\in C(\Omega).$$
Since $f$ is homogeneous of degree $2$ and 
the flow of $V$ preserves the symplectic structure on $\C^3$ it follows that
the immersion $X_f:M \ra \Sph^5$ remains Legendrian.

To make a quantitative statement we assume now that we have a Legendrian immersion
$X:D \ra \Sph^5(R)$, ($R\ge 1$), where $D$ is a disk of radius $1$ in the Euclidean plane $\R^2$,
and $X$ satisfies 
\begin{equation}
\addtocounter{theorem}{1}
\label{E:inequality}
\Vert \partial X:C^{2,\beta}(D,g_0)\Vert  \le 1, \quad \text{and} \quad g \ge \epsilon g_0,
\end{equation}
where $\partial X$ are the partial derivatives of the coordinates of $X$, 
$g$ is the metric induced by the immersion $X$ and $g_0$ is the standard
Euclidean metric on $D$. Note that this condition can be arranged by 
first appropriately magnifying the target (see for example \ref{L:quadratic}).

We have the following proposition.
\begin{prop}
\addtocounter{equation}{1}
\label{P:Leg:pert}
There exists a (small) constant $c_1(\epsilon)>0$
such that if $X$ is a Legendrian immersion satisfying \ref{E:inequality}
and the functions $\phi, \varphi:D \ra \R$ satisfy
$$ \Vert \phi:C^{2,\beta}(D,g_0)\Vert  < c_1, \quad \Vert \varphi:C^{2,\beta}(D,g_0)\Vert  < c_1,$$ 
then $X_\phi, X_{\phi+\varphi}: D \ra \Sph^5(R)$ are well-defined by the construction above and satisfy
\begin{multline*}
\Vert X_{\phi+\varphi} - X_\phi + (J\nabla\varphi+ 2R^{-1}\varphi J \partial_r)
: C^{1,\beta}(D,g_0)\Vert 
\le \\
C (\Vert \phi:C^{2,\beta}(D,g_0)\Vert +\Vert \varphi:C^{2,\beta}(D,g_0)\Vert )\,
\Vert \varphi:C^{2,\beta}(D,g_0)\Vert ,
\end{multline*}
and
\begin{multline*}
\Vert \theta_{\phi+\varphi} - \theta_{\phi} - (\Delta + 6R^{-2})\varphi:C^{0,\beta}(D,g_0)\Vert 
\le \\
C (\Vert \phi:C^{2,\beta}(D,g_0)\Vert +\Vert \varphi:C^{2,\beta}(D,g_0)\Vert )\,
\Vert \varphi:C^{2,\beta}(D,g_0)\Vert ,
\end{multline*}
where $\theta_\phi$ and $\theta_{\phi+\varphi}$ are the Lagrangian angles of $X_\phi$ and $X_{\phi+\varphi}$ respectively
and the Laplacian $\Delta$ is taken with respect to the metric $g$ induced by $X$.
\end{prop}

\begin{proof}
That the linear terms are as stated is well known and follows by a straightforward calculation we omit.
The nonlinear terms are given by rational functions of monomials consisting of contractions
of derivatives of $X$ and derivatives of $\varphi$ and $\phi$.
This implies both the existence results and the estimate on the nonlinearity.
\end{proof}

Finally we record the following proposition we use in the paper:

\begin{prop}
\label{P:killing}
\addtocounter{equation}{1}
If $f: \C^3 \ra \R$ is a harmonic Hermitian quadratic, that is a function of the form
$$
f = \sum_{i,j=1}^3 a_{ij}z_i\bar{z}_j, \quad \text{with} \quad a_{ij}=\bar{a}_{ij}, \quad \sum_{i=0}^3{a_{ii}}=0, \quad a_{ij}\in\C,
$$
then $-J\nabla f$ is a Killing field corresponding to an element of $\mathfrak{su}(3)$, 
and its restriction to any minimal Legendrian surface $M$ of $\Sph^5$ satisfies
$$ (\Delta + 6)f = 0,$$
where the Laplacian $\Delta$ is taken with respect to the metric induced on $M$ from $\Sph^5$. 
\end{prop}
\begin{proof}
The proof can be found in \cite[Thm. 3.2]{fu} or in \cite[Lem. 3.4]{joyce:slgcones2}. 
\end{proof}

\bibliographystyle{amsplain}
\bibliography{paper}
\end{document}